\documentclass[a4paper,10pt,DIV=12]{scrartcl}

\usepackage[intlimits]{amsmath}
\usepackage{amssymb}
\usepackage{amsthm}
\usepackage{bm}
\usepackage{paralist}
\usepackage{eurosym}
\usepackage{colonequals}
\usepackage{todonotes}
\usepackage[utf8]{inputenc}
\usepackage[T1]{fontenc} 
\usepackage{textcomp}    
\usepackage{graphicx}
\usepackage{overpic}
\usepackage{xspace}

\usepackage{authblk} 
\usepackage{todonotes}
\usepackage{stmaryrd}
\usepackage{nicefrac}

\usepackage{algorithm2e}
\usepackage{siunitx}
\usepackage{caption}
\usepackage{subcaption}

\usepackage{tikz}
\usetikzlibrary{calc}
\usepackage{pgfplots}
\usepgfplotslibrary{dateplot}

\allowdisplaybreaks

\title{\textbf{Simulation of Deformation and Flow in Fractured,
Poroelastic Materials}} 

\author{
\large
\textsc{Katja K. Hanowski}\thanks{Katja.Hanowski@tu-dresden.de}
}
\author{
\large
\textsc{Oliver Sander}\thanks{Oliver.Sander@tu-dresden.de}
}
\affil{Institut f\"ur Numerische Mathematik, TU Dresden, Germany}
\date{\today}

\providecommand{\R}{\mathbb R} 
\providecommand{\N}{\mathbb N} 

\providecommand{\sprod}[2]{\left( #1\, ,\, #2 \right) }

\providecommand{\ekl}[1]{\left[ #1 \right]} 
\providecommand{\rkl}[1]{\left( #1 \right)} 
\providecommand{\skl}[1]{\left \lbrace #1 \right \rbrace} 
\providecommand{\jkl}[1]{\left \llbracket #1 \right \rrbracket} 

\providecommand{\abs}[1]{\left| #1 \right|} 
\providecommand{\norm}[1]{\left\| #1 \right\|} 

\providecommand{\ie}{i.e.}  
\providecommand{\eg}{e.g.}  
\providecommand{\intd}{\,\operatorname{d}}		

\providecommand{\divsub}[1]{\operatorname{div}_{#1}}		
\renewcommand{\d}{\operatorname{d}}		

\DeclareMathOperator{\dist}{dist}	
\renewcommand{\div}{\text{\upshape div}\,} 
\DeclareMathOperator{\spanx}{span}	
\DeclareMathOperator{\trace}{trace}	

\DeclareMathOperator{\vnabla}{\boldsymbol \nabla}
\DeclareMathOperator{\vdiv}{\textbf {div}}			

\newtheorem{thm}{Theorem}

\newtheorem{rem}[thm]{Remark}
\newtheorem{lem}[thm]{Lemma}
\newtheorem{dfn}[thm]{Definition}

\newtheorem{asus}[thm]{Assumptions}

\graphicspath{{img/}{img/source/}}

\makeatletter
\def\input@path{{sections/}{img/source/}}
\makeatother

\begin{document}

\maketitle 


\begin{abstract}

We introduce a coupled system of partial differential equations for the modeling of the fluid--fluid and fluid--solid
interaction in a poroelastic material with a single static fracture. The fluid flow in the fracture
is modeled by a lower-dimensional Darcy equation, which interacts with the surrounding rock
matrix and the fluid it contains. We explicitly allow the fracture to end within the domain, and the fracture width is an unknown of the problem. The resulting weak problem is nonlinear, elliptic and symmetric,
and can be given the structure of a fixed-point problem.
We show that the coupled fluid--fluid problem has a solution in a specially crafted Sobolev space, even though the fracture width cannot be bounded away from zero near the crack tip.

For numerical simulations,
we combine XFEM discretizations for the rock matrix deformation and
pore pressure with a standard lower-dimensional finite element method for the fracture flow problem. The resulting coupled
discrete system consists of linear subdomain problems coupled by nonlinear coupling conditions.
We solve the coupled system with a substructuring solver and observe very fast convergence.
We also observe optimal mesh dependence of the discretization errors even in the presence of crack tips.

\end{abstract}


\section{Introduction}

Coupled fluid--solid interaction processes in fractured porous media play an important role in engineering applications such as
the design and construction of geothermal power plants, the risk assessment of waste deposits, and the
production of crude oil and gas. Numerical simulation of such processes is an important tool for such
applications. Furthermore it can help researchers in geosciences to gain a better understanding
of intricate subsurface processes.
Simulation remains challenging due to the number of physical processes involved,
the nonlinear coupling, the complex geometries, and the heterogeneous nature of fractured porous rock.
In particular, such systems combine hydrological processes such as fluid flow in the porous
matrix and in the fracture network with mechanical effects: the deformation of the medium under fluid
pressure and external loads.

In this work we focus on the nonlinear coupling between hydrology
and mechanics.
We consider a low-porosity medium containing a single large-scale fracture, which is filled with a porous medium without mechanical stiffness.
Also, we assume that the fracture length scale is much larger than the fracture width.
The pore space of the matrix and the fracture are both assumed to be fully saturated with a fluid.
The mechanical and hydrological equilibrium state of the system is then governed by four coupling processes:
\begin{enumerate} \label{enum:couplings}
  \item \label{enum:coupling-ff} Fluid--fluid coupling: The fluid may diffuse from the fracture into the surrounding medium and vice versa.
  \item \label{enum:coupling-fs} Fluid--solid coupling: The fluid in the fracture exerts a normal force onto the fracture boundaries, which induces a
	deformation of the rock matrix.
  \item \label{enum:coupling-sf} Solid--fluid coupling: The deforming rock matrix changes the fracture domain, which affects the permeability in the fracture.
  \item \label{enum:coupling-p} Poroelasticity: Fluid pressure in the rock matrix influences the matrix stiffness.
\end{enumerate}
The third list item is the main challenge. The crack width enters the fracture flow equation as the inverse permeability,
which renders the overall system nonlinear. The problem becomes even more challenging if the crack ends within the domain, which we explicitly allow. In such a case singularities appear near the crack tip both in the mechanical stresses and in the matrix fluid flow field. These singularities need to be captured by special singularity functions as part of the XFEM discretization. Optimal convergence of the discretization errors can only be observed if these functions are selected correctly. Furthermore, the crack width necessarily tapers off near the crack tip, and therefore cannot be bounded away from zero on the entire domain. Our proof of the existence of solutions to the coupled bulk--fracture fluid problem constructs particular weighted Sobolev spaces to deal with this degeneracy.

\bigskip
Elasticity problems in fractured materials have been widely addressed in the literature during the last century. Analytical solutions have been derived
for several special cases \cite{inglis,rice}, and the singular behavior of the solution near the fracture tip has been investigated \cite{mazia-rossmann-kozlov}.
Various numerical methods have been developed, often in connection with the modeling of fracture propagation \cite{rabczuk-et-al,aliabadi}. A breakthrough in the numerical
modeling of crack growth processes was
the eXtended Finite Element Method (XFEM) \cite{moes-dolbow-belytschko,yazid-abdelkader-abdelmadjid,pommier-et-al}, which enriches the
discrete function spaces by non/polynomial functions locally reproducing the discontinuity along the fracture and the stress singularities at the crack tip.
XFEM overcomes the need to adapt the mesh to the discontinuities and singularities of the solution.

Concerning the modeling and simulation of fluid flow in fractured porous media,
we only mention some recent articles related to the model presented in this paper.
In \cite{martin-jaffre-roberts} a coupled bulk--fracture fluid model is derived from standard single-phase Darcy equations where the fracture is represented
as a lower-dimensional interface. It is assumed that the fracture and surrounding matrix are both filled
with porous media with different material properties, and that the fracture separates the domain into two parts (i.e., there is no crack tip).
Under the assumption that there exists a lower bound for the fracture aperture function,
existence and uniqueness
of solutions in standard Sobolev spaces is proved.
The authors present a discrete domain decomposition formulation of the original transmission problem using Raviart--Thomas Finite Elements. In \cite{angot-boyer-hubert}, the model from~\cite{martin-jaffre-roberts}
is extended to consider fractures with crack tips.
A cell-centered finite volume scheme is applied, and again, a grid adapted to the fracture is required for the discretization of the problem. Existence of a solution is proved by showing the convergence of the finite volume discretization to a function in a subspace of $H_{\text{div}}$.

Since the discretizations proposed in these two papers rely on a grid resolving the fracture, considering multiple interacting fractures or fracture networks becomes computationally expensive. A more flexible discretization is introduced
in \cite{dangelo-scotti}, where the authors use an extended $(RT_0,P_0)$ Raviart--Thomas discretization of the mixed problem. The authors prove consistency, stability and convergence of the proposed numerical scheme in the case of a domain fully cut by a fracture. In \cite{formaggia-et-al}, the ideas of the model from \cite{martin-jaffre-roberts} are generalized to
flows in fracture networks. The fluid pressure and velocity are allowed to jump at the fracture intersections,
but the fluid interaction with the surrounding matrix is not considered. Again, an extended $(RT_0,P_0)$ Raviart--Thomas discretization is applied
to overcome the difficulty of matching grids at the intersections of multiple fractures.
The model from \cite{martin-jaffre-roberts} has been extended to two-phase flow \cite{jaffre-mnejja-roberts,fumagalli-scotti-13} and passive transport \cite{fumagalli-scotti-11}.

In contrast, literature addressing the coupling of flow and deformation in fractured media is scarce. In \cite{watanabe-et-al} the Biot theory is used to describe the matrix flow and deformation, and a reduced model is used for the flow in the network. Both models are coupled
by a crack width function similar to our own approach. The authors use lower-dimensional interface elements to represent the fracture, and an XFEM
discretization to capture the discontinuity of the matrix displacement and the fracture intersections. Nonetheless, to enforce pressure continuity, they require that the interface is resolved by the bulk grid.

An alternative model for the coupling of deformation and flow uses the quasi-static Biot equation, a linearized model for slightly compressible single-phase fluids in the surrounding media, and a lubrication Darcy equation for modeling the flow in the fracture \cite{girault-wheeler-ganis-mear, girault-kumar-wheeler}. Contrary to our own approach, this model includes changes in the porosity of the bulk medium, but requires continuity of the fluid pressure across the interface. Existence and uniqueness in weighted Sobolev spaces is proved for a simplified version of the fracture equation in which the fracture permeability is not influenced by  the fracture width. The problem is discretized using finite elements on matching grids.

The article \cite{garagash-detournay-adachi} considers a fully
saturated porous medium with a semi-infinite fracture filled with a viscous fluid. The fluid flow inside the fracture is modeled using lubrication
theory, and the fluid leak-off into the surrounding medium is modeled using Carter's leak-off law. From known analytical solutions to the displacement of an elastic
solid and lubrication theory, the authors derive asymptotic solutions for the fluid pressure inside the fracture, the fracture aperture function and the fracture opening under the assumption
that the fracture propagates at a constant velocity.

Similarly, in \cite{fries-schaetzer-weber}, the surrounding medium is modeled as an impermeable, homogeneous, isotropic linear elastic solid. The fracture flow is modeled using a standard lubrication equation for incompressible flow between parallel plates. The equations are coupled by using the fluid pressure as Neumann boundary conditions on the fracture. The authors present an efficient XFEM based discretization using a hybrid explicit--implicit crack description. The main focus of the paper is the determination of the stress intensity factors without the influence of rigid body motions of the crack tip.

Finally, a phase-field approach to fracture propagation is presented in \cite{mikelic-wheeler-wick}. There, the variational approach to fracture propagation presented in
\cite{bourdin-francfort-marigo} is extended to fully saturated porous media that contains a fracture. Contrary to \cite{garagash-detournay-adachi},
the fracture propagates only according to an energy minimization principle and may bifurcate. The numerical modeling of this problem remains
a challenging and computationally expensive task (see also \cite{wheeler-wick-wollner} and references therein).

\bigskip

The model and discretization proposed in this paper combine and generalize several of the previous approaches. In particular, we consider the interplay of \emph{three} processes, viz. the fluid flow in the fracture and the bulk matrix, and the elastic deformation of the matrix. As the matrix deformation determines the fracture width, which in turn determines the fracture permeability, the resulting coupled three-field problem is nonlinear in an intricate way.

More specifically, we combine the model from \cite{martin-jaffre-roberts} for the coupled fluid flow in the fracture and in the surrounding medium with the Biot equation for a linear poroelastic material. Both models are coupled by defining the crack width as the normal jump of the displacement and by
prescribing the fracture fluid pressure as a normal force applied to the solid skeleton at the fracture fronts. We derive the weak formulation of the problem based on the primal form of the Darcy equations and the Biot equation. This weak form has a natural formulation as a fixed-point equation. The corresponding fixed-point iteration alternates between solution operators for the elastic deformation and for the coupled matrix-fracture flow problem. Both solution operators are linear, and the overall nonlinearity stems only from the nonlinear coupling of the deformation to the flow problem.
Our model allows cracks to end within the bulk. While our prototype
geometry in Figure~\ref{fig:domains} has only one crack tip, we can easily
handle the case of a fracture with two tips, and of two-dimensional fractures in a three-dimensional bulk,
where the crack tip is a one-dimensional line.
The width function acts as degenerated and singular coefficient in the fluid equations.

As the main theoretical result we prove existence and uniqueness of weak solutions of the coupled matrix--fracture flow problem for fixed fracture width. This is reasonably straightforward if the fracture width is bounded from below away from zero (see \cite{martin-jaffre-roberts, angot-boyer-hubert}). However, in our case the fracture width is an
$H^{\frac 12}_{00}$ function, and hence must go to zero near the crack tip.
From the asymptotic expansion of the deformation field near the tip we deduce a general form of the fracture width function, which is used as weight in the definition of the solution spaces of the fluid problem. With the proper definition of these spaces, we can prove existence and uniqueness of a solution of the fluid problem with fixed fracture width.

Unfortunately, proving existence of solutions of the overall fixed-point problem is not a matter of simply using an appropriate fixed-point theorem. Since the fracture width function enters the definition of the solution spaces for the fluid-problem, the fixed-point iteration actually operates on a sequence of iteration-dependent spaces. We leave this for future work.

 We use an XFEM discretization both for the displacement and the fluid pressure.  As a result we do not have any restriction on how the fracture is positioned relative to the bulk grid (unlike, e.g., \cite{watanabe-et-al, girault-wheeler-ganis-mear}).
The major goal when designing a XFEM discretization is to make the space large enough to obtain optimal discretization error behavior for decreasing mesh size. This is particularly challenging in the presence of fracture tips, where singularities in the solution are to be expected. The optimal error behavior of correctly constructed XFEM discretizations for linear elasticity problems is well-known, and has been proved in \cite{nicaise-renard-chahine}.
Asymptotic analysis for the fluid problem in the bulk cannot be carried out straightforwardly, but known results can be used to conjecture the form of the crack tip singularity. Again, based on this result, the bulk solution of the pressure may admit a discontinuity across the fracture and a velocity singularity at the front.
 However, it is not clear at all whether the same still holds when combining such problems in a nonlinear fashion as we do here. While we have not attempted to prove rigorous bounds for the discretization errors of our XFEM approximation for the coupled
fluid--fluid--elasticity problem, we observe optimal rates in numerical experiments. This leads us to conclude that the XFEM spaces used here are, in a certain sense, the correct ones.

\bigskip
The paper is organized as follows: In Chapter 2 we introduce the governing equations for the coupled processes. As a preparation for the weak formulation and the existence proofs, Chapter 3 introduces several weighted Sobolev spaces and shows some of their relevant properties. In Chapter 4, a weak formulation of the nonlinearly coupled problem is derived and we show existence and uniqueness of weak solutions of the fluid--fluid subproblem. In Chapter 5 we introduce an XFEM discretization of the resulting problem.
Finally, Chapter 6 provides numerical examples showing that the discretization error behaves optimally.


\section{The Coupled Deformation--Flow Model}

We begin by stating the continuous model in its strong form.  We do this in two steps:
First we state the coupled equations for a system with a fully-dimensional
fracture.  In a separate step we then perform a dimension reduction of the fracture.  This approach allows to better distinguish the features introduced by the coupling itself
from those introduced by the dimension reduction.

\subsection{Coupling Flow in a Full-Dimensional Fracture to a Poroelastic Bulk}
\begin{figure}
 \begin{center}
	\begin{tikzpicture}[scale=0.75]
	    \draw[thick,black,fill=lightgray!30!white] (0,0) ellipse (4cm and 3cm);
	    \draw[thick,domain=0:3.5,smooth,variable=\x,dashed,black]  plot ({-\x},{2.0*sqrt(0.05*\x) + 0.05*\x*\x});
	    \draw[thick,domain=0:3.9,smooth,variable=\x,dashed,black]  plot ({-\x},{-2.0*sqrt(0.05*\x) + 0.05*\x*\x});
	    \node[font=\selectfont] at (2,1) {$\Omega_e$};
	    \node[font=\selectfont] at (-3,0.3) {$\Omega_f$};
	    \node[font=\selectfont] at (3,2.6) {$\widetilde{\Gamma}$};
\end{tikzpicture}
 \end{center}
 \caption{Domain $\widetilde \Omega$ with volume fracture $\Omega_f$}
 \label{fig:full_dimensional_fracture}
\end{figure}
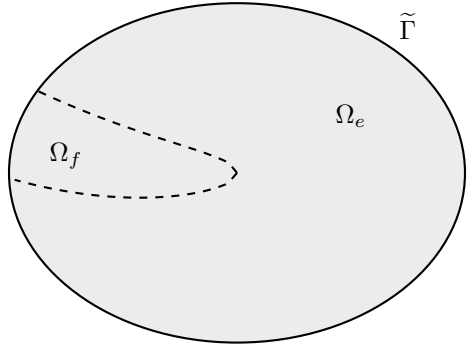

 Let $\widetilde \Omega \subset \R^d$ be a domain of dimension $d=2$
 or $d=3$ with Lipschitz-boundary $\widetilde \Gamma \colonequals \partial \widetilde \Omega$.
 We assume that $\widetilde{\Omega}$ contains a fracture, which, in this first step, we model
 as a subdomain $\Omega_f \subsetneq \widetilde{\Omega}$ (Figure~\ref{fig:full_dimensional_fracture}). Its complement
 $\widetilde \Omega \setminus \Omega_f$ will be denoted by $\Omega_e$.

\subsubsection{Flow Equations in the Fracture and the Matrix}

We consider the matrix and the fracture to both consist of porous media, and to be both
fully saturated with a single-phase incompressible
fluid. These assumptions lead to the well-known Darcy equation
	\begin{subequations} \label{bulk:darcy-all}
	  \begin{alignat}{2}
	   \div \mathbf q &= f_F &\qquad& \text{on $\widetilde \Omega$}, \label{bulk:darcy-mass-all} \\
	   \mathbf q &= -\mathbb K \nabla p && \text{on $\widetilde \Omega$} \label{bulk:darcy-momentum-all},
	  \end{alignat}
	\end{subequations}
  where $\mathbf q$ denotes the seepage velocity, $p$ the pore pressure, and $\mathbb K \in \R^{d \times d}$ the permeability tensor.
  In this formulation, the fracture $\Omega_f$ appears as a region of $\widetilde{\Omega}$ where the
permeability $\mathbb K$ differs considerably from the rest.  To prepare the inclusion into a coupled model, however,
we formulate a substructuring problem.  Is is well known~\cite{quarteroni_valli:1999} that \eqref{bulk:darcy-all} is equivalent
to solving the Darcy equation on the subdomains $\Omega_e$ and $\Omega_f$ separately, and imposing suitable coupling conditions.
We therefore consider the system
	  \begin{alignat*}{2}
	  \div \mathbf q^\Omega &= f_F^\Omega &\qquad& \text{in $\Omega_e$}, \\ 
	  \mathbf q^\Omega &= -\mathbb K^\Omega \nabla p^\Omega && \text{in $\Omega_e$},
	  \end{alignat*}
for the bulk and
\begin{subequations} \label{fracture:darcy-volume}
	\begin{alignat}{2}
	\div \mathbf q^f &= f^f &\qquad& \text{in $\Omega_f$}, \label{fracture:darcy-mass-volume}  \\
	\mathbf q^f &= -\mathbb K^f \nabla p^f &\qquad& \text{in $\Omega_f$}, \label{fracture:darcy-momentum-volume}
	\end{alignat}
\end{subequations}
for the fracture.  The appropriate coupling conditions for this are continuity of the pressure
 \begin{equation} \label{mod:pressurecont}
  p^\Omega = p^f \qquad \text{on } \partial \Omega_e \cap \partial \Omega_f,
 \end{equation}
and continuity of the normal flux
 \begin{equation} \label{mod:fluxcont}
  \mathbf q^\Omega \cdot \nu_f = \mathbf q^f \cdot \nu_f
  \qquad \text{on } \partial \Omega_e \cap \partial \Omega_f,
 \end{equation}
where $\nu_f$ denotes the unit outer normal of $\Omega_f$.
Well-posedness of this problem together with suitable boundary conditions is shown in \cite{quarteroni_valli:1999}.
%
\subsubsection{Poroelastic Behavior of the Matrix}

While the fracture and bulk are modeled as behaving qualitatively the same as far as hydrological processes
are concerned, their mechanical behavior is considered to be different. While the bulk is described as a linear poroelastic material,
we do not assign any stiffness to the fracture at all. Hence mechanically it behaves
as if it was empty.

Let $\mathbf{u} : \Omega_e \to \R^d$ be the displacement of the rock matrix.
We make the small strain assumption, which allows to use the linearized strain
\begin{equation*}
  \mathbf e (\mathbf u) \colonequals \frac 12 \big( \vnabla \mathbf u + (\vnabla \mathbf u )^T  \big)
\end{equation*}
and the St.\,Venant--Kirchhoff material law
  \begin{equation*}
    \sigma (\mathbf u) = \lambda \trace \left( \mathbf e (\mathbf u) \right) \mathbb I +
	  2 \mu \mathbf e (\mathbf u)
  \end{equation*}
for the elastic stress $\sigma$.
The parameters $\lambda$ and $\mu$ are the well-known Lam\'e coefficients.
In a poroelastic medium, internal forces result from the elastic stress $\sigma$ and the isotropic
fluid pressure $p^\Omega \mathbb I$, where $\mathbb I$ is the $d \times d$ identity tensor. Their equilibrium is described by the Biot equation
  \begin{equation*} 
    -\vdiv \rkl{\sigma(\mathbf u) - p^\Omega \mathbb I} = \mathbf f_E,
  \end{equation*}
  where $\mathbf f_E$ denotes a volume force acting on the medium. We omit the possible dependence of the matrix permeability on the deformation field.

We formulate coupling conditions for the bulk--fracture interface based on conservation of momentum. As the fracture is supposed to be empty as far as mechanical behavior is concerned, there is no displacement variable in the fracture. Consequently, we cannot require continuity of the displacement at the interface. We therefore only postulate equality of the normal components of the total stress, which is $p^f \mathbb I$ in the fracture and $\sigma(\mathbf u) - p^\Omega \mathbb I$ in the surrounding matrix:
\begin{equation} \label{mod:stresscont}
 \sigma( \mathbf u) \nu_f - p^\Omega \nu_f = -p^f \nu_f
\end{equation}
on $\partial \Omega_e \cap \partial \Omega_f$. The continuity of pressure \eqref{mod:pressurecont} then yields the  boundary condition $\sigma( \mathbf u) \nu_f = \mathbf 0$ on the interface.

\subsection{Dimension Reduction of the Fracture} \label{sec:dimension_reduction}

Fractures are typically long and thin objects.
Following \cite{martin-jaffre-roberts}, we therefore replace the $d$--dimensional fracture by a $(d-1)$--dimensional hypersurface,
and the equations on $\Omega_f$ by reduced equations obtained by integrating \eqref{fracture:darcy-volume} across the fracture thickness.
The coupling conditions are modified accordingly. Unlike \cite{martin-jaffre-roberts}, we take the curvature of the fracture midsurface into account.

\begin{figure}
    \centering
    \begin{subfigure}{.49\textwidth}
      \centering
      \begin{tikzpicture}[scale=0.75]
	    \draw[thick,black,fill=lightgray!30!white] (0,0) ellipse (4cm and 3cm);

	    \draw[thick,domain=0:3.5,smooth,variable=\x,dashed,black]  plot ({-\x},{2.0*sqrt(0.05*\x) + 0.05*\x*\x});
	    \draw[thick,domain=0:3.9,smooth,variable=\x,dashed,black]  plot ({-\x},{-2.0*sqrt(0.05*\x) + 0.05*\x*\x});
	    \draw[thick,domain=0:3.9,smooth,variable=\x,dashed,solid,red!70!black]  plot ({-\x},{0.05*\x*\x});
		\draw[thick,green!60!black, dashed] (0,0) -- (4,0);
		\draw[thick,<->,blue!70!black] (-2,-0.45) -- (-1.8,0.75);
	    \node[font=\selectfont] at (2,1.6) {$\widetilde{\Omega}$};
	    \node[font=\selectfont,red!70!black] 
		    at (-3.5,0.2) {$\Sigma$};
		\node[font=\selectfont,blue!70!black] 
		    at (-2.1,0.5) {$b$};
	    \node[font=\selectfont,green!60!black] 
		    at (3,-0.4) {$\Sigma_e$};
	    \node[font=\selectfont] at (3,2.6) {$\widetilde{\Gamma}$};
	    \draw[black,fill=black] (0,0) circle (2pt);
\end{tikzpicture}
    \end{subfigure}%
    \begin{subfigure}{.49\textwidth}
      \centering
            \begin{tikzpicture}[scale=0.75]
	    \draw[thick,black,fill=lightgray!30!white] (0,0) ellipse (4cm and 3cm);
 	    \draw[thick,domain=0:3.9,smooth,variable=\x,dashed,solid,black]  plot ({-\x},{0.05*\x*\x});
		\draw[thick,black, dashed] (0,0) -- (4,0);
		\draw[thick,black, ->] (-1,0.05) -- (-0.9,1);
		\draw[thick,black, ->] (0,3) -- (0,4);
		\node[font=\selectfont] at (-0.65,0.5) {$\nu$};
		\node[font=\selectfont] at (0.3,3.5) {$\mathbf n$};
	    \node[font=\selectfont] at (2,1.6) {$\Omega^+$};
	    \node[font=\selectfont] at (2,-1.6) {$\Omega^-$};
	    \node[font=\selectfont] 
		    at (-2.5,0.7) {$\Sigma^+$};
		\node[font=\selectfont] 
			    at (-2.6,0) {$\Sigma^-$};
	    \node[font=\selectfont] at (3,2.6) {$\Gamma^+$};
	    \node[font=\selectfont] at (3,-2.6) {$\Gamma^-$};
	    \draw[black,fill=black] (0,0) circle (2pt);
\end{tikzpicture}
    \end{subfigure}
    \caption{Replacing a thin fracture by a lower--dimensional approximation. Left:
    	Fractured domain with fracture midsurface $\Sigma$ and tangential extension $\Sigma_e$.
    	Right: Dimension--reduced geometry.}
    \label{fig:domains}
\end{figure}
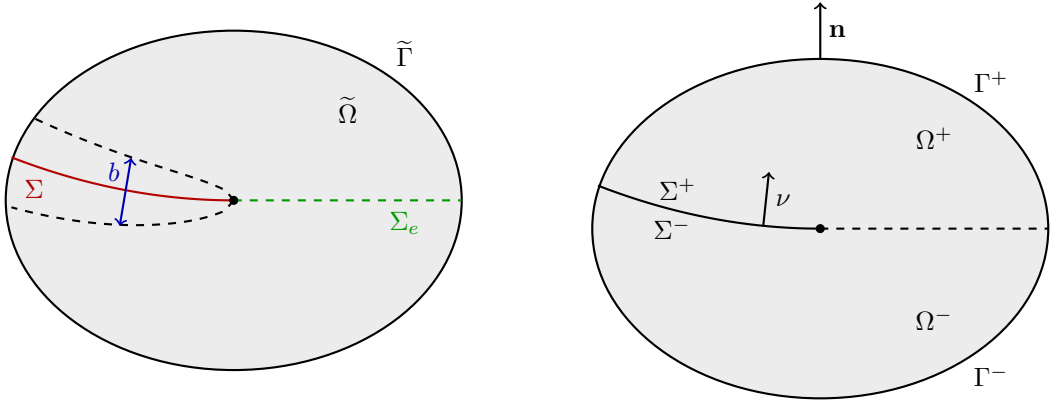

For the rest of this article we suppose that there is
a parametrized hypersurface $\Sigma$ (called the fracture midsurface) such that
the fracture domain $\Omega_f \subset \widetilde \Omega$ can be
represented by
\begin{equation*}
 \Omega_f = \Big \lbrace \mathbf x \in \widetilde \Omega \, \big \vert \,
	\mathbf x = \mathbf s + t\nu, \; \mathbf s \in \Sigma, \; t \in
	\Big( - \frac{b(\mathbf s)}{2}, \frac{b(\mathbf s)}{2} \Big) \Big \rbrace,
\end{equation*}
 where $\nu$ denotes a continuous unit normal vector field to $\Sigma$, and $b:\Sigma \to (0,\infty)$ is the fracture aperture function
 (Figure~\ref{fig:domains}).
	The new bulk domain is $\Omega \colonequals \widetilde \Omega \setminus \Sigma$, with outer boundary $\Gamma \colonequals \partial \Omega \setminus \Sigma$.

 We distinguish two cases. Either, the fracture $\Sigma$ partitions $\widetilde \Omega$ into two disconnected subdomains. In that case, we suppose that the two domains both have Lipschitz boundary, and we label them $\Omega^+$ and $\Omega^-$, respectively. In the other case, $\Omega$ is connected, which implies that at least parts of the fracture boundary $\gamma \colonequals \partial \Sigma$ are contained in $\Omega$. We then suppose that
 there exists a tangential extension $\Sigma_e$ of $\Sigma$ such that
 $\widetilde \Sigma \colonequals \Sigma \cup \Sigma_e$
 subdivides $\widetilde \Omega$ into two disjoint subdomains $\Omega^+$ and $\Omega^-$
 with Lipschitz-boundaries. In either way we denote the boundaries of $\Omega^+$ and $\Omega^-$ by $\Gamma^+$ and $\Gamma^-$, respectively.
  We denote by $\nu^\pm$
  the unit outer normal to $\Sigma^\pm \colonequals \Gamma^\pm \cap \widetilde \Sigma$
  and by $\mathbf n$ the unit outer normal to the outer boundaries
  $\Gamma^\pm \setminus \widetilde \Sigma$. To be specific, we set $\nu \colonequals \nu^- = -\nu^+$.

  \bigskip
  Dimension reduction of the fracture equation \eqref{fracture:darcy-volume} involves splitting up
  the equation into normal and tangential parts.
  The projection operators onto the normal and tangent spaces of the parametrized hypersurface $\Sigma$ are denoted by $\mathbb P_\nu \colonequals \nu \nu^T$ and
  $\mathbb P_\tau \colonequals \mathbb I - \mathbb P_\nu$, respectively.
  For scalar-valued or vector-valued functions $g$ and $\mathbf g$, we
  define the normal derivative by
  $\nabla_\nu g \colonequals \nabla g \mathbb P_\nu$ and $\vnabla_\nu \mathbf g \colonequals \ \vnabla \mathbf g \mathbb P_\nu$, respectively.
  The normal divergence
  operator $\divsub{\nu}$ is
  \[\divsub{\nu}  \mathbf g \colonequals \trace \left( \vnabla_\nu \mathbf g \right)
  = \mathbb P_\nu : \vnabla \mathbf g.\]
  The tangential gradient $\nabla_\tau$ and divergence $\divsub{\tau}$ are defined analogously.
  Finally, we introduce the average and the jump operator on $\Sigma$ by
    \[  \skl{g} \colonequals \frac 12 \rkl{g \big \vert_{\Sigma^+ } + g \big \vert_{\Sigma^- }}
    \qquad \text{and} \qquad
    \jkl{g} \colonequals \rkl{g \big \vert_{\Sigma^+ } - g \big \vert_{\Sigma^- }},\]
  respectively.

Suppose that the fracture permeability tensor $\mathbb K^f$ decomposes additively as $\mathbb K^f = K^\nu \mathbb P^\nu + K^\tau \mathbb P^\tau$ with constants $K^\nu, K^\tau > 0$. Then Equation \eqref{fracture:darcy-momentum-volume}
can be decoupled into tangential and normal parts
  \begin{subequations}
  	 \begin{align}
  	 \mathbb P^\tau \mathbf q^f &= - K^\tau \nabla_\tau p^f, \label{fracture:darcy-momentum-volume-tangential} \\
  	 \mathbb P^\nu \mathbf q^f &= - K^\nu \nabla_\nu p^f. \label{fracture:darcy-momentum-volume-normal}
  	 \end{align}
  \end{subequations}
   Define the averaged fracture pressure $p^\Sigma$ and the averaged tangential seepage
     velocity $\mathbf q^\Sigma$ by
     \[ p^\Sigma (\mathbf s) \colonequals \frac{1}{b(\mathbf s)}\int_{-\frac {b(\mathbf s)}{2}}^{\frac {b(\mathbf s)}{2}} p^f (\mathbf s + t\nu) \intd t \
     \qquad \text{and} \qquad
     \mathbf q^\Sigma (\mathbf s) \colonequals  \frac{1}{b(\mathbf s)}
     \int_{-\frac {b(\mathbf s)}{2}}^{\frac {b(\mathbf s)}{2}} \mathbb P^\tau (\mathbf s) \mathbf q^f (\mathbf s + t\nu) \intd t. \]
    The aim is to express equations \eqref{fracture:darcy-mass-volume}, \eqref{fracture:darcy-momentum-volume-tangential} and
    \eqref{fracture:darcy-momentum-volume-normal} in terms of
    these averaged quantities. Simple calculations show that
    \[ \div \mathbf q^f = \divsub{\nu} \mathbf q^f + \divsub{\tau} \mathbf q^\tau
    + \divsub{\tau} \mathbf q^\nu = \divsub{\nu} \mathbf q^f + \divsub{\tau} \mathbf q^\tau + \kappa \mathbf q^f \cdot \nu, \]
    where $\kappa = \divsub{\tau} \nu$ is the mean curvature of $\Sigma$, $\mathbf q^\nu \colonequals \mathbb P^\nu \mathbf q^f$ and $\mathbf q^\tau \colonequals \mathbb P^\tau \mathbf q^f$.
    Let $\mathbf s \in \Sigma$ be arbitrary.
    In a first step, we integrate the left-hand side of \eqref{fracture:darcy-mass-volume} in normal direction,
    and apply the Gauss and Leibniz integral rules to obtain
    \begin{align} \nonumber
    	\int_{-\frac {b(\mathbf s)}{2}}^{\frac {b(\mathbf s)}{2}} \div \mathbf q^f (\mathbf s + t\nu) \intd t 
    		&= \int_{-\frac {b(\mathbf s)}{2}}^{\frac {b(\mathbf s)}{2}} \divsub{\nu} \mathbf q^f (\mathbf s + t\nu) \intd t 
    			+ \int_{-\frac {b(\mathbf s)}{2}}^{\frac {b(\mathbf s)}{2}} \divsub{\tau} \mathbf q^\tau (\mathbf s + t\nu) \intd t
    			 \\
    		\nonumber & \quad	+ \int_{-\frac {b(\mathbf s)}{2}}^{\frac {b(\mathbf s)}{2}} \kappa(\mathbf s) \mathbf q^f(\mathbf s + t\nu)\cdot \nu(\mathbf s) \intd t  \\
    		&= \ekl{\mathbf q^f \Big(\mathbf s + \frac{b(\mathbf s)}{2} \nu \Big) -
    			\mathbf q^f \Big(\mathbf s - \frac{b(\mathbf s)}{2} \nu\Big)}\cdot \nu(\mathbf s)
	    	 		\label{eq:averageDivergence_step2_1} \\
    		& \quad + \divsub{\tau} \int_{-\frac {b(\mathbf s)}{2}}^{\frac {b(\mathbf s)}{2}}
    			 \mathbf q^\tau (\mathbf s + t\nu) \intd t \label{eq:averageDivergence_step2_2}\\
    		& \quad 
    			 - \frac 12 \rkl{\mathbf q^\tau \Big(\mathbf s 
    			 - \frac{b(\mathbf s)}{2} \nu\Big) + \mathbf q^\tau \Big(\mathbf s + \frac{b(\mathbf s)}{2} \nu\Big)}\cdot \nabla_\tau b(\mathbf s)
    			 \label{eq:averageDivergence_step2_3} \\
	    	& \quad + \kappa(\mathbf s)  
		    	\int_{-\frac {b(\mathbf s)}{2}}^{\frac {b(\mathbf s)}{2}}
			    \mathbf q^\nu(\mathbf s + t\nu) \intd t \cdot \nu(\mathbf s). \label{eq:averageDivergence_step2_4}  
    \end{align}
    To simplify the first term \eqref{eq:averageDivergence_step2_1}, note that
	it only involves values of $\mathbf q^\nu$ on the fracture boundary $\partial \Omega_f \cap \partial \Omega_e$.
    The continuity of fluxes \eqref{mod:fluxcont} then yields
    \begin{equation*} \label{eq:fluxcont_form}
	    \mathbf q^f  \Big(\mathbf s \pm \frac{b(\mathbf s)}{2}\nu(\mathbf s)\Big) \cdot \nu (\mathbf s) = \mathbf q^\Omega  \Big( \mathbf s \pm \frac{b(\mathbf s)}{2}\nu(\mathbf s)\Big) \cdot \nu (\mathbf s). \tag{8b}
	\end{equation*}
    Since $b$ is assumed to be small, we use the continuity of fluxes \eqref{eq:fluxcont_form} to approximate
    $\mathbf q ^\Omega \big( s + \frac{b\rkl{\mathbf s}}2 \nu \mathbf (s)\big) \cdot \nu(\mathbf s)$
    by $\mathbf q^\Omega (\mathbf s)\big\vert_{\Sigma^+} \cdot \nu(\mathbf{s})$
    and $\mathbf q^\Omega \big( \mathbf s - \frac{b(\mathbf s)}{2}\nu(\mathbf s)\big)\cdot \nu (\mathbf s)$
    by $\mathbf q^\Omega (\mathbf s)\big\vert_{\Sigma^-} \cdot \nu(\mathbf{s})$.
    Together, this yields
    \[ \ekl{\mathbf q^f \Big(\mathbf s + \frac{b(\mathbf s)}{2} \nu \Big) -
    	\mathbf q^f \Big(\mathbf s - \frac{b(\mathbf s)}{2} \nu \Big)}\cdot \nu(\mathbf s)
	    \approx \jkl{\mathbf q^\Omega} (\mathbf s) \cdot \nu (\mathbf s). \]
    To rewrite \eqref{eq:averageDivergence_step2_3}, define
    \[ \widetilde{\mathbf q}^\tau (\mathbf s) \colonequals  \frac 12 \rkl{\mathbf q^\tau \Big(\mathbf s - \frac{b(\mathbf s)}{2} \nu\Big) + \mathbf q^\tau \Big( \mathbf s + \frac{b(\mathbf s)}{2} \nu\Big)}. \]
    Applying the trapezoidal rule and the continuity of the normal flux \eqref{eq:fluxcont_form} to the integral in
    \eqref{eq:averageDivergence_step2_4} gives
    \[\int_{-\frac {b(\mathbf s)}2}^{\frac {b(\mathbf s)}2}  \mathbf q^\nu(\mathbf s + t \nu) \intd t \cdot \nu (\mathbf s) \approx b(\mathbf s) \skl{ \mathbf q^\nu}(\mathbf s)\cdot \nu(\mathbf s) = b(\mathbf s)\skl{\mathbf q^\Omega}(\mathbf s) \cdot \nu(\mathbf s).
    \]
    Finally, observe that the integral in \eqref{eq:averageDivergence_step2_2}
    is simply $b\mathbf q^\Sigma$.
    In total, the divergence of $\mathbf q^f$ is represented by
    \begin{align*}
    \int_{-\frac {b(\mathbf s)}{2}}^{\frac {b(\mathbf s)}{2}} \div \mathbf q^f (\mathbf s + t\nu) \intd t &\approx \jkl{\mathbf q^\Omega} (\mathbf s)\cdot \nu (\mathbf s)
    + \divsub{\tau} (b(\mathbf s)\mathbf q^\Sigma(\mathbf s)) - \widetilde{\mathbf q}^\tau (\mathbf s) \cdot \nabla_\tau b (\mathbf s) + b(\mathbf s)\kappa(\mathbf s) \skl{\mathbf q^\Omega}(\mathbf s) \cdot \nu(\mathbf s).
    \end{align*}
    Similarly, since $ K^\tau$ is constant, integrating the right-hand side of equation \eqref{fracture:darcy-momentum-volume-tangential} yields
    \begin{align*}
	    \int_{-\frac {b(\mathbf s)}{2}}^{\frac {b(\mathbf s)}{2}} -  K^\tau \nabla_\tau
	    p^f(\mathbf s + t\nu) \intd t &\approx -  K^\tau \rkl{\nabla_\tau \rkl{bp^\Sigma} -
	    	\skl{p^\Omega} \nabla_\tau b }
		    = -  K^\tau \ekl{b \nabla_\tau \rkl{p^\Sigma} -
		    	\rkl{\skl{p^\Omega} - p^\Sigma} \nabla_\tau b },
    \end{align*}
    where we have used the continuity of the pressure $p^f = p^\Omega$.
    Assuming that $\nabla_\tau b$ is small, this yields
    \begin{align*}
    \divsub{\tau} \rkl{b \mathbf q^\Sigma} &= b f^\Sigma -  \jkl{\mathbf q^\Omega}  \cdot \nu - b \kappa \skl{\mathbf q^\Omega} \cdot \nu , \\ 
    \mathbf q^\Sigma &= -  K^\tau \nabla_\tau p^\Sigma , 
    \end{align*}
    on $\Sigma$, where for any $\mathbf s \in \Sigma$ we define
    \[
     f^\Sigma (\mathbf s) \colonequals \frac{1}{b(\mathbf s)}	\int_{-\frac {b(\mathbf s)}{2}}^{\frac {b(\mathbf s)}{2}} f^f (\mathbf s + t\nu) \intd t.   \]

    Finally, we approximate the equation \eqref{fracture:darcy-momentum-volume-normal} by integrating in normal direction. We use the trapezoidal rule, as before, for the left-hand side and
    apply the fundamental theorem of calculus to the right--hand side, to obtain
    \begin{equation}  \label{fracture:averaged-normal-momentum}
    b \skl{\mathbf q^\Omega}\cdot \nu = - K^\nu \jkl{p^\Omega}. 
    \end{equation}

  To close the system, we need to relate the fracture fluid pressure $p^\Sigma$ to the bulk
  fluid pressure $p^\Omega$ and the bulk flow $\mathbf q^\Omega$.
  To derive the corresponding formula, we will use the composite trapezoidal rule
  \begin{equation} \label{eq:pressure_trapezoidal}
	  \frac{1}{b(\mathbf s)} \int_{- \frac{b(\mathbf s)}{2}}^{\frac{b(\mathbf s)}{2}}
	 p^f (\mathbf s + t\nu) \intd t \approx
	 \frac{1}{4} \ekl{p^f \Big(\mathbf s - \frac{b(\mathbf s)}{2} \nu\Big)
				 	+ 2 p^f (\mathbf s) + p^f \Big(\mathbf s +  \frac{b(\mathbf s)}{2} \nu
				 	\Big)}.
	\tag{12a}
  \end{equation}
 We use the linear approximations
 \begin{align*}
 	p^f(\mathbf s) &\approx p^f \Big(\mathbf s + \frac{b(\mathbf s)}{2}\nu \Big) -
	 	\frac{b(\mathbf s)}{2} \nabla p^f \Big(\mathbf s + \frac{b(\mathbf s)}{2}\nu\Big) \cdot \nu
	 	(\mathbf s)\\
	 	&=  p^f \Big(\mathbf s + \frac{b(\mathbf s)}{2}\nu \Big) +
		 	\frac{b(\mathbf s)}{2K^\nu} \mathbf q^f \Big(\mathbf s + \frac{b(\mathbf s)}{2}\nu\Big)
		 	 \cdot \nu(\mathbf s).\\
	\intertext{By symmetry, we can also write}
	p^f(\mathbf s) &\approx  p^f \Big(\mathbf s - \frac{b(\mathbf s)}{2}\nu\Big) -
	\frac{b(\mathbf s)}{2K^\nu} \mathbf q^f \Big(\mathbf s - \frac{b(\mathbf s)}{2}\nu \Big)
	\cdot \nu(\mathbf s).
 \end{align*}
 Averaging these two expressions yields the approximation
 \begin{align*} 	
	 p^f(\mathbf s) &\approx \frac 12 \ekl{p^f \Big(\mathbf s - \frac{b(\mathbf s)}{2}\nu
	 	\Big)
	 	+ p^f \Big(\mathbf s + \frac{b(\mathbf s)}{2}\nu\Big)} \tag{12b} \\
	  & \quad + \frac{b(\mathbf s)}{2K^\nu}
			 \ekl{\mathbf q^f \Big(\mathbf s + \frac{b(\mathbf s)}{2}\nu \Big)
			  - \mathbf q^f \Big(\mathbf s - \frac{b(\mathbf s)}{2}\nu\Big)
			 	} \cdot \nu(\mathbf s). \tag{12b} \label{eq:averaged_pressure_approximation} 
 \end{align*}
 By the continuity of the pressure and the fluxes on $\partial \Omega_f \cap \partial \Omega_e$ and the approximations $p^\Omega \big \vert_{\partial \Omega_f \cap \partial \Omega^\pm} \approx
	 p^\Omega \big \vert_{\Sigma^\pm}$ and $\mathbf q^\Omega\cdot \nu \big \vert_{\partial \Omega_f \cap \partial \Omega^\pm} \approx
	 \mathbf q^\Omega \cdot \nu \big \vert_{\Sigma^\pm}$ it then follows that
 \begin{align*}
 	p^\Sigma(\mathbf s)
  		&= \frac{1}{b(\mathbf s)} \int_{- \frac{b(\mathbf s)}{2}}^{\frac{b(\mathbf s)}{2}}
  				p^f (\mathbf s + t\nu) \intd t \approx \frac 14 \ekl{4\skl{p^\Omega}(\mathbf s) 
  					+ 2 \frac{b(\mathbf s)}{4K^\nu} \jkl{\mathbf q^\Omega}(\mathbf s) \cdot \nu (\mathbf s) } \\
			 &= \skl{p^\Omega}(\mathbf s) + \frac{b(\mathbf s)}{8K^\nu}
			 \jkl{\mathbf q^\Omega}(\mathbf s) \cdot \nu (\mathbf s).
\end{align*}  
 This equation is a special case of the more general form
 \begin{equation} \label{fracture:average-pressure-approximation}
 p^\Sigma =  \skl{p^\Omega} + (2 \xi - 1) \frac{b}{4K^\nu} \jkl{\mathbf q^\Omega}\cdot \nu
 \end{equation}
 used in \cite{martin-jaffre-roberts}, where $\xi \in \left( \frac 12, 1 \right]$ depends on the choice of
 the discretization of the integral in the definition of $p^\Sigma$. The trapezoidal rule yields $\xi = \nicefrac 12$. With the composite trapezoidal rule \eqref{eq:pressure_trapezoidal} and the averaged expression \eqref{eq:averaged_pressure_approximation}, we have $\xi = \nicefrac 34$.
%
%
\subsection{The Fully Coupled Problem with the Reduced Fracture Equation}
  Combining the dimension-reduced fracture fluid equations with the
  bulk flow and poroelasticity equation of the previous section, we arrive at the following problem:
  Find a fracture fluid pressure $p^\Sigma: \Sigma \rightarrow \R$,
  a bulk fluid pressure $p^\Omega: \Omega \rightarrow \R$, and a bulk
  displacement $\mathbf u: \Omega \rightarrow \R^d$ such that
  \begin{subequations}
  \begin{alignat}{2}
   	-\vdiv \rkl{\sigma(\mathbf u) - p^\Omega \mathbb I} &= \mathbf f_E
   	&\qquad& \text{in }
   	\Omega \label{all:biot} \\
  	\div \mathbf q^\Omega &= f_F^\Omega &\qquad& \text{in }
  	\Omega \label{all:darcy-mass}, \\
  	\mathbf q^\Omega &= -\mathbb K^\Omega \nabla p^\Omega &\qquad& \text{in }
  	\Omega,
  	\label{all:darcy-momentum} \\
  	\divsub{\tau} \rkl{b \mathbf q^\Sigma} &= b f^\Sigma_F -  \jkl{\mathbf q^\Omega}  \cdot \nu
	  	- b \kappa \skl{\mathbf q^\Omega}\cdot \nu &\qquad& \text{on } \Sigma , \label{all:averaged-mass} \\
  	\mathbf q^\Sigma &= - K^\tau \nabla_\tau p^\Sigma
  	&\qquad&\text{on } \Sigma\label{all:averaged-tangential-momentum}.
  \end{alignat}
  \end{subequations}
  The first equation is the momentum balance of poroelasticity, followed by two equations for the matrix Darcy flow and two equations for the reduced fracture Darcy flow. The three processes are coupled partly by source terms in the equations themselves, and partly by explicit coupling conditions. In particular, the matrix fluid pressure $p^\Omega$ appears in the poroelastic momentum balance \eqref{all:biot}.
  The fluid flow $\mathbf q^\Omega$ from the bulk to the fracture appears as a volume term in the fracture Darcy flow equation \eqref{all:averaged-mass}. Conversely,
  by \eqref{fracture:averaged-normal-momentum} and \eqref{fracture:average-pressure-approximation}, the fracture fluid
  pressure $p^\Sigma$ acts as Robin-type boundary condition for the bulk flow $\mathbf q^\Omega$:
  \begin{subequations} \label{all:coupling-conditions}
  	\begin{alignat}{2}
  		 p^\Sigma &=  \skl{p^\Omega} + (2 \xi - 1) \frac{b}{4K^\nu} \jkl{\mathbf q^\Omega}\cdot \nu 	&\qquad&\text{on } \Sigma
  		 \label{all:average-pressure-approximation}, \\
  		  \skl{\mathbf q^\Omega}\cdot \nu &= - \frac{K^\nu}{b} \jkl{p^\Omega}
  		  	&\qquad&\text{on } \Sigma.
  		  \label{all:averaged-normal-momentum}
  	\end{alignat}
  	It acts as a Neumann boundary condition
  	for the poroelastic equation
  	\begin{equation} \label{all:coupling-stress}
  	\sigma(\mathbf u\vert_{\Sigma^\pm}) \cdot \nu^\pm = - p^\Sigma\nu^\pm
  	\end{equation}
  	by the continuity of the total stress \eqref{mod:stresscont}.
  	Finally, the displacement $\mathbf u$ determines the width of the fracture.
  	In the reduced model, this width appears only in form of the fracture aperture function $b:\Sigma \rightarrow \R$. In a linear elastic setting it is reasonable to set the fracture width equal to the normal jump of the displacement field. We therefore impose the
  	coupling condition
  	\begin{equation} \label{all:coupling-width}
  	b = \jkl{\mathbf u}\cdot \nu \qquad \text{on } \Sigma.
  	\end{equation}
  \end{subequations}
  This is the ``difficult'' coupling condition which makes the coupled system nonlinear and possibly degenerate.
  
  The system is closed using appropriate Dirichlet and Neumann boundary
  conditions for all three subsystems. Let $\Gamma_N^E$ and $\Gamma_D^E$ be two disjoint sets
  with $\Gamma_N^E \cup \Gamma_D^E = \Gamma \setminus \gamma$. For the displacement we have the
  boundary conditions
  \begin{subequations}
  	\begin{alignat}{2}
  	\sigma(\mathbf u) \cdot \mathbf n & = \sigma_N  &\qquad & \text{on } \Gamma_N^E,\label{poroelasticity:neumannBC} \\
  	\mathbf u & =  \mathbf u_D  &\qquad & \text{on } \Gamma_D^E, \label{poroelasticity:dirichletBC}
   \end{alignat}
   for given boundary data functions $\sigma_N$ and $\mathbf u_D$.
   Similarly, let $\Gamma_N^F$ and $\Gamma_D^F$ be two disjoint sets with
   $\Gamma_N^F \cup \Gamma_D^F = \Gamma \setminus \gamma$ and let $q_N^\Omega$ and $p_D^\Omega$ be given data. The bulk flow boundary conditions are
   \begin{alignat}{2}
     \mathbf q^\Omega \cdot \mathbf n &= q_N^\Omega &\qquad & \text{on } \Gamma_N^F,\label{darcy:neumannBC} \\
     p^\Omega &= p_D^\Omega && \text{on } \Gamma_D^F. \label{darcy:dirichletBC}
  \end{alignat}
  Finally, the reduced fracture flow problem needs boundary conditions of $\Sigma$. Generally, the boundary of $\Sigma$ is in part a subset of the domain boundary $\partial \overline{\Omega}$ and in part contained in the interior of $\overline{\Omega}$ (the fracture tip). We define the two disjoint sets $\gamma_N^F$ and $\gamma_D^F$ with $\gamma_N^F \cup \gamma_D^F = \gamma$ and
  set the boundary conditions
  \begin{alignat}{2}
    \mathbf q^\Sigma \cdot \tau &= q_N^\Sigma &\qquad & \text{on } \gamma_N^F,\label{fracture:neumannBC} \\
    p^\Sigma &= p_D^\Sigma && \text{on } \gamma_D^F. \label{fracture:dirichletBC}
  \end{alignat}
  \end{subequations}
\section{Sobolev Spaces for the Bulk--Fracture System}
	Standard existence theory for the Darcy equation requires the permeability to be bounded from below away from zero almost everywhere. Unfortunately, this assumption
	does not hold if the fracture ends in the interior of the bulk domain, because coupling condition
	\eqref{all:coupling-width} forces the fracture width $b$ to tend to zero when approaching a crack tip.
	For a rigorous existence theory we therefore have to resort to weighted Sobolev spaces.
	\subsection{Sobolev Spaces on Domains with a Slit} \label{sec:fs_domain_with_slit}
		  The bulk elasticity and fluid problems are posed on a domain with a slit. Solutions of elliptic equations on such domains are usually not first-order Sobolev functions. Instead, we construct certain weighted Sobolev spaces in which elliptic problems become well-posed.

		  Remember that the fracture domain $\Sigma$ (possibly with extension $\Sigma_e$) divides the domain $\Omega$ into two Lipschitz domains $\Omega^+$ and $\Omega^-$. Let $L^2(\Omega^\pm)$ and $H^1(\Omega^\pm)$ be the standard Lebesgue and Sobolev spaces on $\Omega^\pm$, with the norms $\| \cdot \|_{0,\Omega^\pm}$
		  and $\| \cdot \|_{1,\Omega^\pm}$, and scalar products
		  $\rkl{\cdot, \cdot}_{0,\Omega^\pm}$ and $\rkl{\cdot, \cdot}_{1,\Omega^\pm}$, respectively. We then
		  define
		  \[ \widetilde V = \left \{ v \in L^2 (\Omega) \, \Big \vert \, v|_{\Omega^\pm} \in H^1\rkl{\Omega^\pm} \right \}, \]
		  and the broken scalar product
		  \[ \rkl{\cdot,\cdot}_{\widetilde V}:\widetilde V \times \widetilde V \rightarrow \mathbb R, \qquad
			  \rkl{u,v}_{\widetilde V} \colonequals
			   \rkl{u,v}_{H^1\rkl{\Omega^-}} +
			   \rkl{u,v}_{H^1 \rkl{\Omega^+}}.   \]
		  The space $\widetilde V$ is a Hilbert space, and
		  the induced norm is
		  \[ \norm{\cdot}_{1,\Omega}^2 \colonequals \rkl{\cdot, \cdot}_{\widetilde V}
				  =  \| \cdot \|_{1, \Omega^-}^2 +  \| \cdot \|_{1, \Omega^+}^2.\]

		  For any $(d-1)$-dimensional Lipschitz manifold $M$ we denote the
		  standard Sobolev--Slobodeckij space on $M$ by
		  $H^{\frac 12} (M)$ (see, \eg, \cite{wloka}) and the corresponding  norm by
		  $\norm{\cdot}_{\frac 12, M}^2 \colonequals \norm{\cdot}_{0,M}^2 +
			 \abs{\cdot}_{\frac 12, M}^2$,
		  where $\abs{\cdot}_{\frac 12, M}^2$ is the seminorm induced by the symmetric
		  bilinear form
		  \[  \rkl{u,v}_{\frac 12, M} \colonequals
			 \int_{M} \int_M \frac{\abs{u(\mathbf s) - u(\widetilde {\mathbf s})}
			 	\abs{v(\mathbf s) - v(\widetilde{\mathbf s})}}
				 {\abs{\mathbf s - \widetilde{\mathbf s}}^d} \intd \mathbf s \intd \widetilde {\mathbf s}. \]
		  Since $\Omega^+$ and $\Omega^-$ are Lipschitz domains there exist unique linear and continuous trace
		  operators from $H^1(\Omega^\pm)$ onto $H^{\frac 12} ( \Gamma^\pm)$. We
		  need their restrictions to the outer boundary $\Gamma$ and the extended
		  slit $\widetilde \Sigma$.
		  \begin{dfn} \label{lem:elasticity_trace_domainboundary}
		    Denote by $\gamma_\Gamma :
		    \widetilde V \rightarrow H^{\frac 12}(\Gamma)$,
		    $\gamma^+ : H^1\rkl{\Omega^+}
		    \rightarrow H^{\frac 12} ( \widetilde \Sigma)$ and $\gamma^-:
		    H^1\rkl{\Omega^-} \rightarrow H^{\frac 12} (\widetilde\Sigma)$ the
		    restrictions of the global trace operators from $H^1(\Omega^\pm)$ onto $H^{\frac 12} ( \Gamma^\pm)$ to $\Gamma$ and $\widetilde \Sigma$.
		  \end{dfn}
		  In particular the jump and average operators $\jkl{\cdot},\skl{\cdot }: \widetilde V \rightarrow H^{\frac 12} (\widetilde \Sigma)$, given by
		  \[ \jkl{v}_{\widetilde{\Sigma}}
			  \colonequals \left(\gamma^+v - \gamma^-v \right) \qquad \text{and} \qquad   \skl{v}_{\widetilde{\Sigma}} \colonequals \frac 12\left(\gamma^+v + \gamma^-v \right), \]
		  are well defined, linear and continuous. The restrictions of these operators to $\Sigma$ will be denoted
		  by $\jkl{\cdot}$ and $\skl{\cdot}$ (without a subscript), respectively.

          Functions in $\widetilde V$ can have arbitrary jumps across the entire extended fracture $\widetilde \Sigma$. To single out the functions that only jump at the actual interface $\Sigma$ we now define
          the subspace
          \[ V \colonequals \skl{v \in \widetilde V \mid \jkl{v}_{ \Sigma_e} = 0 \text{ almost everywhere}}. \]
          By the continuity of the trace operators
          it follows that $V$ is a closed subspace of $\widetilde V$. This is the usual definition of a first-order Sobolev Space on the non-Lipschitz domain $\Omega \setminus \Sigma$.

		  We now discuss the space of traces on $\Sigma$ of function from $V$. Define the subspace $H^{\frac 12}_0 (\Sigma)$ of $H^{\frac 12}(\Sigma)$ as the completion of $C^{1,1}(\Sigma)$ with compact
		  support in the $H^{\frac 12}$-norm. Let the extension of any $v \in H^{\frac 12}_0 (\Sigma)$ by $0$ be denoted by $\widetilde v$, \ie,
		  \[ \widetilde v = \begin{cases}
		                     v & \text{on } \Sigma, \\ 0 & \text{on } \Sigma_e.
		                    \end{cases}\]
		  This extension $\widetilde v$ is not generally in $H^{\frac 12} (\widetilde \Sigma)$ \cite{khludnev}. The subspace of functions
		  in $H^{\frac 12}_0$ whose extensions are in $H^{\frac 12}(\widetilde \Sigma)$ is called
		  \[ H^{\frac 12}_{00} (\Sigma) := \left \lbrace v \in H^{\frac 12}_0(\Sigma) \mid \widetilde v \in H^{\frac 12} (\widetilde \Sigma) \right \rbrace. \]
		  Let $\d :\Sigma \rightarrow \mathbb R$ denote the geodesic distance of $\mathbf s \in \Sigma$ to $\gamma \cap \widetilde \Omega$.
          We define a scalar product
		  $\rkl{\cdot,\cdot}_{00,\Sigma} :  H^{\frac 12}_{00}
		  (\Sigma) \times H^{\frac 12}_{00} (\Sigma) \rightarrow \mathbb R$ by
		  \[ \rkl{u,v}_{00,\Sigma} \colonequals \rkl{u,v}_{0,\Sigma}
				  + \rkl{u,v}_{\frac 12,\Sigma}
				  + ( u\d^{- \nicefrac 12}, v \d^{- \nicefrac 12} )_{0,\Sigma}. \]
		  The space $ H^{\frac 12}_{00} (\Sigma)$ equipped with
		  this scalar product is a Hilbert space \cite{khludnev}
		  with induced norm $\| \cdot \|_{00, \frac 12, \Sigma}$. Furthermore, a function $v$ is in $H^{\frac 12}_{00}(\Sigma)$ if and only if the norm $\| v \|_{00, \frac 12, \Sigma}$ is finite \cite{khludnev}.
		  
  		  We introduce the space of admissible traces
  		  \[ W_\Sigma \colonequals \skl{\rkl{ v^+, v^-} \in H^{\frac 12}(\Sigma)
  		  	\times   H^{\frac 12}(\Sigma) \mid
  		  	\rkl{v^- - v^+} \in
  		  	H^{\frac 12}_{00}(\Sigma)}. \]
  		  A norm on this space is
  		  \[ \norm{\rkl{ v^-, v^+}}^2 =
  		  \norm{v^-}_{\frac 12, \Sigma}^2
  		  +	 \norm{v^+}_{\frac 12, \Sigma}^2
  		  +  \big\|\rkl{v^+ - v^-}\d^{-\nicefrac 12}\big\|_{0,\Sigma}^2. \]
  		  The following result establishes a relation between the space $V$ of Sobolev functions on the slit domain $\Omega$, and the space $W_\Sigma$
  		  of traces. It shows that $W_\Sigma$ is the correct trace space of $V$ on $\Sigma$.
  		  \begin{lem}[\cite{khludnev}, Theorem 1.25] \label{lem:elasticity_trace_continuity}
  		  	\begin{enumerate}[(i)]~
	  		  	\item There exists a continuous linear trace operator
		  		  	 $\gamma_\Sigma:V \rightarrow W_\Sigma, v \mapsto
		  		  	 \rkl{\gamma^-v, \gamma^+v}$.
		  		\item There exists a continuous linear extension operator
			  		$ E_\Sigma : W_\Sigma \rightarrow V, \rkl{v^-, v^+}
			  		\mapsto v$ such that $\gamma_\Sigma \circ E_\Sigma = \text{id}$.
  		  	\end{enumerate}
  		  \end{lem}
		  In particular, all traces of functions in $V$ 
		  are $H^{\frac 12}_{00} (\Sigma)$-functions.
		  \bigskip
		  All definitions in this section can be made equally well for vector-valued spaces. Such spaces will be written in bold face. For $d$-valued spaces, we
		  mention the following generalized Green's formula, which is proved in  \cite{angot-boyer-hubert}.
		  \begin{lem} \label{lem:elasticity_generalized_green}
		   Let $\Omega$ and $\Sigma$ satisfy the conditions of Section
		   \ref{sec:dimension_reduction}.
		   Let $\mathbf u \in \mathbf L^2(\Omega)$ satisfy $\div \mathbf u \in L^2(\Omega)$ and let $v \in V$. Then
		   \[ \int_\Omega (\div \mathbf u)v + \mathbf u \cdot \nabla v  \normalfont{\intd} \mathbf x =  \rkl{\mathbf u \cdot
		   	\mathbf n, v}_{0,\Gamma}
			  - \rkl{\jkl{\mathbf u \cdot \nu}, \skl{v}}_{0,\Sigma}
		      - \rkl{\skl{u \cdot \nu}, \jkl{v}}_{0,\Sigma}.\]
		  \end{lem}

	\subsection{Weighted Lebesgue Spaces on Parametrized Hypersurfaces}
		 In this section we will assume that the fracture width function $b: \Sigma \rightarrow \R$ is fixed.
		 The functions $b$ and $b^{-1}$ appear as degenerate and singular coefficients in the averaged fracture fluid equation \eqref{all:averaged-mass} and
		 the fluid coupling conditions
		 \eqref{all:average-pressure-approximation},
		 \eqref{all:averaged-normal-momentum}. We must therefore resort to weighted Sobolev spaces on the fracture to obtain well-posed fluid problems. We recall the basic definitions.
		 \begin{dfn} \label{dfn:generalWeights}
		 	Let $\widehat \omega: \R^{N} \rightarrow \mathbb R$ be a function and
		 	let $M \subset \R^N$ be open and bounded.
		 	\begin{enumerate}[(i)]
		 		\item The function $\widehat \omega$ is called a \emph{weight} if
			 		$\widehat \omega \in L^1_{\text{loc}}(\R^{N})$ and $\widehat\omega > 0$ almost everywhere.
			 	\item Let $\widehat \omega$ be a weight. The
				 	\emph{weighted $L^2$-space} $L^2_{\widehat\omega}\rkl{M}$
				 	is defined as
			 		\[ L^2_{\widehat\omega}\rkl{M} \colonequals \skl{v:M \rightarrow \mathbb R
			 			\mid v \text{ is measurable},
			 			\norm{v}_{\widehat \omega,M} < \infty}, \]
			 		where the norm
			 		\[ \norm{v}_{\widehat \omega,M}^2 \colonequals
				 		\int_M \abs{v}^2 \widehat \omega \intd \mathbf x \]
				 	is induced by the weighted scalar product
				 	\[ \rkl{u,v}_{\widehat \omega,M} \colonequals
					 	\rkl{\widehat{\omega}^{\frac 12}u, \widehat{\omega}^{\frac 12} v}_{0,M}. \]
		 		\item The space $L^2_{\widehat \omega,0}(M)$ is defined as the
			 		closure of $C^\infty_0(M)$ in $L^2_{\widehat \omega}\rkl{M}$.
		 		\item A weight $\widehat \omega$
			 		 is called an \emph{$A_2$-weight}
			 		if there exists a constant $A > 0$ such that for all balls $B$ in $\R^{N}$ we have
			 		\[ \Bigg(\frac{1}{|B|} \int_B \widehat \omega(\mathbf x) \intd \mathbf x \Bigg)\Bigg(\frac{1}{|B|}
			 		\int_B \widehat \omega(\mathbf x)^{-1} \intd \mathbf x \Bigg) \le A, \]
			 		independent of $B$.  \label{dfn:A2weight}
		 	\end{enumerate}
		 \end{dfn}
		 If $\widehat \omega$ is an $A_2$-weight, then $L^2_{\widehat \omega}(M)$ and $L^2_{\widehat \omega,0}(M)$
		 are Banach spaces (\cite{turesson}, Proposition 2.1.2). Furthermore,
		 $C^\infty(M) \cap L^2_{\widehat \omega} (M)$ is then dense in $L^2_{\widehat \omega}(M)$ (\cite{miller}, Lemma 2.4).

		 We extend these definitions to parametrized hypersurfaces.
		 \begin{dfn} \label{dfn:weightedLebesgue}
		 	Let $\Sigma \subset \R^d$ be a hypersurface defined by the homeomorphism
		 	$\alpha:\Sigma \rightarrow \widehat \Sigma \subset \R^{d-1}$.
		 	\begin{enumerate}[(i)]
		 		\item We call a function $\omega:\Sigma \rightarrow \mathbb R$ a \emph{weight on $\Sigma$} if there exists 
		 		a weight $\widehat \omega:\R^{d-1} \rightarrow \R$ such that 
		 		$\widehat \omega \vert_{\widehat \Sigma} = \omega \circ \alpha^{-1}$. 
					\label{dfn:weightOnSurface}
				\item  Let $\omega:\Sigma \rightarrow \mathbb R$
				 	be a weight on $\Sigma$. We say that a function $v:\Sigma \rightarrow \R$ belongs to $L^2_\omega(\Sigma)$, if
				 	$v \circ \alpha^{-1} \in L^2_{\widehat \omega} (\widehat \Sigma)$. \label{dfn:weightedLbesgueOnSurface}
				\item Let $\omega:\Sigma \rightarrow \mathbb R$
				be a weight on $\Sigma$. We say that $\omega$ is an 
				$A_2(\Sigma)$-weight, if $\widehat \omega$ is an $A_2$-weight.
				\label{dfn:A2weightOnSurface}
		 	\end{enumerate}
		 \end{dfn}
		 Spaces of functions on parametrized hypersurfaces with $A_2$-weighted norms are Hilbert spaces.
		 \begin{thm} \label{thm:wlebesgue-hilbert}
		 	Let $\Sigma$ be a hypersurface defined by the homeomorphism
		 	$\alpha:\Sigma \rightarrow \widehat \Sigma \subset \R^{d-1}$ and $\omega$ an $A_2(\Sigma)$-weight. The space $L^2_\omega(\Sigma)$ equipped with the norm
		 	\[ \norm{v}_{0,\omega,\Sigma}^2 \colonequals
			 	\norm{v \circ \alpha^{-1}}_{0,\widehat \omega,\widehat \Sigma}^2 \]
		 	is a separable Hilbert space. The norm is induced by the scalar product
		 	\[ \rkl{u,v}_{\omega, \Sigma} \colonequals
			 	\sprod{\widehat{\omega}^{\frac 12} \rkl{u \circ \alpha^{-1}}}{ \widehat{\omega}^{\frac 12}
		 		 \rkl{v \circ \alpha^{-1}}}_{0,\alpha [\widehat \Sigma]}.
		 	\]
		 	Furthermore, $C^\infty(\Sigma) \cap L^2_\omega (\Sigma)$ is dense in
		 	$L^2_\omega (\Sigma)$.
		 \end{thm}
		 \begin{proof}
		 	It is easy to see that the norm $\norm{\cdot}_{0,\omega,\Sigma}$ is induced by the scalar product $\rkl{\cdot,\cdot}_{\omega,\Sigma}$. 

		 	To show the completeness of $L^2_{\omega}(\Sigma)$, let $\rkl{v_n}_{n \in \N}$ be a Cauchy sequence in $L^2_\omega (\Sigma)$. Then the sequence
		 	 $\rkl{v_n \circ \alpha^{-1}}_{n \in N}$ is a Cauchy sequence in $L^2_{\widehat \omega} (\widehat \Sigma)$. 
		 	 Since $L^2_{\widehat \omega} (\widehat \Sigma)$
		 	is complete, there exists a $\widehat v \in L^2_{\widehat \omega} (\widehat \Sigma)$ with $v_n \circ \alpha^{-1} \longrightarrow \widehat v$ for $n \longrightarrow \infty$. Define $v\colonequals
			\widehat v \circ \alpha$, then
			\[ \norm{v_n - v}_{0,\omega,\Sigma}
				= \norm{\rkl{v_n - v}\circ \alpha^{-1}}_{0,\widehat \omega,
					 \widehat \Sigma}
				= \norm{v_n \circ \alpha^{-1} - \widehat v} _{0,\widehat \omega, \widehat \Sigma} \longrightarrow 0 \qquad \text{for} \qquad
				n \rightarrow \infty. \]
			This proves completeness.

			Now, let $v \in L^2_\omega (\Sigma)$ be a function. Then
			$\rkl{v \circ \alpha^{-1}} \in L^2_{\widehat \omega} (\widehat \Sigma)$. Since $L^2_{\widehat \omega}(\widehat \Sigma) \cap C^\infty (\widehat \Sigma)$ is dense in $L^2_{\widehat \omega}(\widehat \Sigma)$, there exists a sequence $\rkl{\widehat v_n}_{n \in \N}$
			in $C^\infty(\widehat \Sigma) \cap L^2_{\widehat \omega} (\widehat \Sigma)$ converging towards $\rkl{v \circ \alpha^{-1}}$ in
			$ L^2_{\widehat \omega} (\widehat \Sigma)$. For $\rkl{v_n}_{n \in \N}$ defined by $v_n \colonequals \widehat v_n \circ \alpha$ it follows that
			\[ \norm{v_n - v}_{0,\omega,\Sigma}
			= \norm{\rkl{v_n - v}\circ \alpha^{-1}}_{0,\widehat \omega,
				\widehat \Sigma}
			= \norm{\widehat v_n - v\circ \alpha^{-1}} _{0,\widehat \omega, \widehat \Sigma}  \longrightarrow 0 \qquad \text{for} \qquad
			n \rightarrow \infty. \]
		 \end{proof}

		 We now construct weighted Lebesgue spaces using the fracture width function as
		 our weight.
		 \begin{asus} \label{hyp:weight}
		 	Let $\Sigma$ be a hypersurface with boundary $\gamma \colonequals \partial \Sigma$ defined by the homeomorphism
		 	$\alpha:\Sigma \rightarrow \widehat \Sigma$. We assume that $b:\Sigma \rightarrow \mathbb R$ satisfies the following
		 	conditions:
		 	\begin{enumerate}[{A}1]
		 		\item The function $b$ is a weight on $\Sigma$. \label{hyp:weight1}
		 		\item Let $\chi$ be a smooth cutoff function with $\chi = 1$ in a neighborhood of
			 		$\gamma$. There exists a constant $c>0$ such that the function $b$ is given by
			 		\[b = c (\chi \dist^{\frac 12} + (1-\chi) f), \]
			 		where $f$ denotes an integrable function and $\dist: \R^d \rightarrow \R$ denotes the distance in $\R^d$ to $\gamma$.  \label{hyp:weight2}
			 	\item The function $b$ is essentially bounded on $\Sigma$ by a constant
				 	 $b_{\text{max}} > 0$.  \label{hyp:weight3}
				\item The function $f$ is essentially 
					bounded from below away from zero,
					\ie, there exists a constant $a > 0$ such that
					$a \le f$ almost everywhere on $\Sigma$.  \label{hyp:weight4}
		 	\end{enumerate}
		 \end{asus}
		 The crucial assumption A\ref{hyp:weight2} is motivated in Chapter \ref{sec:weak_formulation} by asymptotic expansion of the matrix displacement field near the crack tip.
		 \begin{lem}[] \label{lem:ap-weight}
		 	Let $\Sigma$ be a bounded hypersurface defined by the homeomorphism
		 	$\alpha:\Sigma \rightarrow \widehat \Sigma$ with boundary $\gamma = \partial \Sigma$, and let $b:\Sigma \rightarrow \R$ satisfy Assumptions A\ref{hyp:weight1}-A\ref{hyp:weight4}. Then $b$ is an $A_2(\Sigma)$-weight.
		 \end{lem}
		 \begin{proof}
		 	  We define $F = \{ \mathbf s \in \Sigma \mid \chi(\mathbf s) = 1 \}$.
		 	  Since $\dist > 0$ away from the crack tip, there exists a constant
		 	  $d_\text{min} >0$ such that $d_\text{min} \le \dist^{\frac 12}$ on $\Sigma \setminus F$.
		 	  With Assumption A\ref{hyp:weight3} it thus follows that
		 	  \[ b \le b_{\text{max}}
		 	  \le \frac{b_{\text{max}}}
		 	  {d_\text{min}}  \dist^{\frac 12} \qquad \text{on }
		 	  \Sigma \setminus F. \]
		 	  Since $b = c\dist^{\frac 12}$ on $F$
		 	  we conclude that $b \le M \dist^{\frac 12}$ on $\Sigma$, where
		 	  $M \colonequals \max \skl{c, \frac{b_{\text{max}}}{d_\text{min}}}$.
		 	  With Assumptions A\ref{hyp:weight2} and A\ref{hyp:weight4} it furthermore follows that
		 	  \[ b = \chi \dist^{\frac 12} + \rkl{1- \chi} f \ge
			 	  \chi d_\text{min} + \rkl{1 - \chi} a
			 	  \ge \min \skl{d_\text{min}, a} \]
		 	 on $\Sigma \setminus F$. Since $\Sigma$ is bounded, we can find a constant $d_\text{max} > 0$, such that $\dist^{\frac 12} \le d_\text{max}$ on $\Sigma$. This
		 	 yields the estimate
		 	 \[ b^{-1} \le \frac{1}{\min \skl{d_\text{min}, a}} =
		 	 \frac{1}{\min \skl{d_\text{min}, a}} \frac{ d_\text{max}}{ d_\text{max}} \le
			 	 \frac{d_\text{max}}{\min \skl{d_\text{min}, a}} \dist^{-\frac 12}.  \]
		 	 Since $b^{-1} = c^{-1}\dist^{-\frac 12}$ on $F$, we conclude that
		 	 $b^{-1} \le m\dist^{-\frac 12}$, where $m = \max
		 	 \skl{\frac{d_\text{max}}{\min \skl{d_\text{min}, a}}, c^{-1}}$.

		 	 Let $\widehat{\d}: \R^{d-1} \rightarrow \R$ denote the distance of $\mathbf s \in \R^{d-1}$ to the set $\widehat \gamma \colonequals \partial \widehat \Sigma$. The boundedness of $\Sigma$ yields that there exist
		 	 constants $d,D >0$, such that $d \widehat{\d }\rkl{\mathbf s}^{\frac 12} \le
		 	 \dist \rkl{\alpha^{-1}( \mathbf s), \gamma}^{\frac 12} \le  D \widehat{\d} \rkl{\mathbf s}^{\frac 12}$ for all $ \mathbf s \in \widehat \Sigma$. Together with the previous estimates we get
		 	 $ dm \widehat \d \le b\circ \alpha \le DM \widehat \d^{\frac 12}$
		 	 on $\widehat \Sigma$. Finally, by Lemma 3.3 in \cite{duran-lopez}, it follows that the function $\widehat d^{\frac 12}$ is an $A_2$ weight, \ie,
		 	 there exists a constant $\widehat A$, such that
		 	 \[\rkl{\frac{1}{\vert \widehat B \vert}\int_{\widehat B}
		 	 	\widehat d \intd \mathbf s} ^{\frac 12} \rkl{\frac{1}{\vert \widehat B \vert}\int_{\widehat B}
		 	 	\widehat d^{-\frac 12} \intd \mathbf s} \le \widehat A \]
		 	 independent of the choice of the ball $\widehat B \subset \R^{d-1}$.

		 	 Then the function $\widehat b$ defined by 
		 	 \[ \widehat b \colonequals 
			 	 \begin{cases}
				 	 b \circ \alpha^{-1} & \text{on } \widehat \Sigma \\
				 	 \widehat d^{\frac 12} & \text{otherwise}
			 	 \end{cases} \]
		 	 satisfies for arbitrary balls $\widehat B$ in $\R^{d-1}$
		 	 \begin{align*}
		 	 	\rkl{\frac{1}{\vert \widehat B \vert}\int_{\widehat B}
			 	 		\widehat b \intd \mathbf s}\rkl{\frac{1}{\vert \widehat B \vert}\int_{\widehat B}
			 	 		\widehat b^{-1} \intd \mathbf s}
				 	 &= \frac{1}{\vert \widehat B \vert^2}
				 	 \rkl{\int_{\widehat B\cap \widehat \Sigma}
				 	 	\widehat b \intd \mathbf s
				 	 	+ \int_{\widehat B\setminus \widehat \Sigma}
				 	 	\widehat b \intd \mathbf s}
				 	 \rkl{\int_{\widehat B\cap \widehat \Sigma}
				 	 	\widehat b^{-1} \intd \mathbf s +
				 	 	\int_{\widehat B\setminus \widehat \Sigma}
				 	 	\widehat b^{-1} \intd \mathbf s} \\
				 	&\le \frac{\max\skl{1,dm} \max\skl{1,DM}}{\vert \widehat B \vert^2}
				 	\rkl{\int_{\widehat B}
				 		\widehat d^{\frac 12} \intd \mathbf s}
				 	\rkl{\int_{\widehat B}
				 		\widehat d^{-\frac 12} \intd \mathbf s} \\
					 & \le \max\skl{1,dm} \max\skl{1,DM} \widehat A.
		 	 \end{align*}
		 	 Hence, by Definition \ref{dfn:weightedLebesgue}\eqref{dfn:A2weightOnSurface}, $b$ is an 
		 	 $A_2(\Sigma)$-weight.
		 \end{proof}
		 As in the Euclidean case, the fact that $b$ is an $A_2$-weight allows to conclude that the corresponding Lebesgue spaces have desirable properties.
		 \begin{lem}
	 		Let $\Sigma$ be a bounded hypersurface defined by the homeomorphism
	 		$\alpha:\Sigma \rightarrow \widehat \Sigma$ with boundary $\gamma = \partial \Sigma$, and let $b:\Sigma \rightarrow \R$ satisfy Assumptions A\ref{hyp:weight1}-A\ref{hyp:weight4}.
		 	Define the $L^2_{b}(\Sigma)$ and
		 	$L^2_{b^{-1}}(\Sigma)$ with the corresponding norms
		 	$\norm{\cdot}_{0,b,\Sigma}$
		 	and $\norm{\cdot}_{0,b^{-1},\Sigma}$
		 	as in Definition \ref{dfn:weightedLebesgue} and
		 	Theorem \ref{thm:wlebesgue-hilbert}.
		 	Then the spaces $L^2_{b}(\Sigma)$ and $L^2_{b^{-1}}(\Sigma)$ are
		 	separable Hilbert spaces. We have the continuous embeddings
	 		 \[  L^2_{b^{-1}}(\Sigma) \hookrightarrow L^2(\Sigma) \hookrightarrow L^2_{b} (\Sigma). \]
		 \end{lem}
		 \begin{proof}
		 	The assumptions A\ref{hyp:weight1}-A\ref{hyp:weight4} yield that $b$ is an $A_2$-weight and $L^2_b(\Sigma)$ is a separable Hilbert space by Theorem \ref{thm:wlebesgue-hilbert}. By the symmetry of definition \ref{dfn:generalWeights}\eqref{dfn:A2weight}, $b^{-1}$ is an $A_2$-weight as well; which proves the assertion about $L^2_{b^{-1}}(\Sigma)$. The continuous embeddings exist because $b$ is bounded from above.
		 \end{proof}

	\subsection{A Sobolev Space with Weighted Trace Space} \label{sec:fs_weightedTrace}
		In this section we will construct a Sobolev space representing the volume
		pressure $p^\Omega$. Due to the coupling conditions \eqref{all:coupling-conditions},
		we need such a space to have its trace in
		\[  H^{\frac 12}_{b^{-1}} (\Sigma) \colonequals L^2_{b^{-1}}(\Sigma) \cap
			H^{\frac 12}(\Sigma).\]
		To this end, we the norm
		\[
		\| v \|_{\frac 12, b^{-1},\Sigma}^2 \colonequals   \| v \|_{0,b^{-1},\Sigma}^2  +
		\abs{v}_{\frac 12, \Sigma}^2. \]
		From
		$L^2_{b^{-1}}(\Sigma) \hookrightarrow L^2(\Sigma)$ and the definition
		of $H^{\frac 12}_{b^{-1}}(\Sigma)$ it follows that
		$H^{\frac 12}_{b^{-1}}(\Sigma) \hookrightarrow H^{\frac 12}(\Sigma)$.
		 \begin{thm} \label{thm:fluid_volume_hilbert}
			 Let $\Sigma$ be a bounded hypersurface defined by the homeomorphism
			 $\alpha:\Sigma \rightarrow \widehat \Sigma$ with boundary $\gamma = \partial \Sigma$, and let $b:\Sigma \rightarrow \R$ satisfy Assumptions A\ref{hyp:weight1}-A\ref{hyp:weight4}. Then the space $H^{\frac 12}_{b^{-1}}$ equipped with the norm $\| \cdot \|_{\frac 12, b^{-1},\Sigma}$ is a separable Hilbert space.
		 \end{thm}
		 \begin{proof}
		 	The square of the norm $\| \cdot \|_{\frac 12, b^{-1},\Sigma}$ is the sum of squared norms that are induced by scalar products. Therefore
		 	$\norm{\cdot}_{\frac 12, b^{-1},\Sigma}$ is induced by a scalar product
		 	as well.

		 	 Let $(s_n)_{n \in \mathbb N}$ be a Cauchy sequence in $H^{\frac 12}_{b^{-1}}(\Sigma)$. Then  $(s_n)_{n \in \mathbb N}$ is a Cauchy sequence
		 	 in the Hilbert space $L^2_{b^{-1}}(\Sigma)$. Hence there exists
		 	 a $s \in L^2_{b^{-1}}(\Sigma)$ with $s_n \rightarrow s$ in $L^2_{b^{-1}}(\Sigma)$.

		 	 To prove completeness it remains to show that $s \in H^{\frac 12}_{b^{-1}}(\Sigma)$, i.e., that $\abs{s}_{\frac 12, \Sigma}$
		 	 is bounded and that $\abs{s_n - s}_{\frac 12, \Sigma} \rightarrow 0$. But since $L^2_{b^{-1}}(\Sigma)
		 	 \hookrightarrow L^2(\Sigma)$ we can conclude from the corresponding result
		 	 for $H^{\frac 12}(\Sigma)$ (\cite{wloka}, Theorem 3.1) that both conditions are satisfied.

		 	 Finally, the space $H^{\frac 12}_{b^{-1}}(\Sigma)$ can be identified with a closed subspace
		 	 of the separable product space $L^2_{b^{-1}}(\Sigma) \times L^2(\Sigma \times \Sigma)$,
		 	 which implies the separability of $H^{\frac 12}_{b^{-1}}(\Sigma)$
		 	 (see \cite{wloka}, Theorem 3.1).
		 \end{proof}

		The coupling constraints
		\eqref{all:average-pressure-approximation} and
		\eqref{all:averaged-normal-momentum} require that the trace space on $\Sigma$ of the volume pressure space should be $H^{\frac 12}_{b^{-1}}(\Sigma)$. Therefore we define
		\[ V_{b^{-1}} \colonequals \skl{ v \in V \mid \gamma^\pm v \in H^{\frac{1}{2}}_{b^{-1}} (\Sigma) }, \]
		and equip this space with the norm
		\[ \| v \|_{1,b^{-1},\Omega}^2 = \| \gamma^+v \|_{0,b^{-1},\Sigma}^2 + \| \gamma^-v \|_{0,b^{-1},\Sigma}^2 + \|v \|_{1,\Omega}^2. \]
		It is easy
		to see that $V_{b^{-1}} \hookrightarrow V$.

		\begin{lem} \label{lem:fluid_bulk_trace}
			The restriction of $\gamma_\Gamma:V \rightarrow H^{\frac 12}(\Gamma)$
			from Definition \ref{lem:elasticity_trace_domainboundary}
			onto an operator from $V_{b^{-1}}$ to $H^{\frac 12}(\Gamma)$,
			and the restrictions of $\gamma^\pm : V \rightarrow H^{\frac 12}(\Sigma)$
			from Definition \ref{lem:elasticity_trace_domainboundary}
			to operators from $V_{b^{-1}}$ to $H^{\frac 12}(\Sigma)$
			 are well defined,
			linear and continuous.
		\end{lem}
		\begin{proof}
			The embedding
			$V_{b^{-1}} \hookrightarrow V$
			yields that the restriction of $\gamma_\Gamma:V \rightarrow H^{\frac 12}(\Gamma)$
			onto an operator from $V_{b^{-1}}$ is well defined, linear and continuous.

			Analogously,
			the restrictions of the trace operators $\gamma^\pm : V \rightarrow
			H^{\frac 12}(\Sigma)$
			onto operators from
			$V_{b^{-1}}$ to $H^{\frac 12}(\Sigma)$ are well defined and linear.
			By the definition of $V_{b^{-1}}$ it follows that the restricted operators $\gamma^\pm$ map onto $H^{\frac 12}_{b^{-1}}(\Sigma)$.
			Furthermore, for any $v \in V_{b^{-1}}$ the following estimate is satisfied:
			\begin{align*}
			\| \gamma^+ v \|_{\frac 12, b^{-1},\Sigma}^2 &= \| \gamma^+ v \|_{0, b^{-1},\Sigma}^2
			+ \abs{\gamma^+ v}_{\frac 12,\Sigma}^2 \\
			& \le  \| \gamma^+ v \|_{0, b^{-1},\Sigma}^2 +  \| \gamma^+ v \|_{\frac 12,\Sigma}^2  \\
			& \le C \rkl{ \| \gamma^+ v \|_{0, b^{-1},\Sigma}^2 + \|  v \|_{1,\Omega^+}^2 }.
			\end{align*}
			Analogously one can show that $\| \gamma^- v \|_{\frac 12,b^{-1},\Sigma}^2 \le C \rkl{ \| \gamma^- v \|_{0, b^{-1},\Sigma}^2 + \|  v \|_{1,\Omega^-}^2  }$. This proves
			continuity.
		\end{proof}
		As in Section \ref{sec:fs_domain_with_slit}, we introduce the space of admissible traces. Define
		\[W_{b^{-1}} \colonequals \skl{\rkl{v^-,v^+} \in H^{\frac 12}_{b^{-1}}(\Sigma)
			\times   H^{\frac 12}_{b^{-1}}(\Sigma) \mid
			\rkl {v^+ - v^-} \in
			H^{\frac 12}_{00}(\Sigma)} \]
		with norm
		\[ \norm{\rkl{v^-,v^+}}^2_{W_{b^{-1}}} \colonequals
		\norm{\gamma^- v}_{\frac 12,b^{-1}, \Sigma}^2
		+	 \norm{\gamma^+ v}_{\frac 12,b^{-1}, \Sigma}^2. \]
		This norm is equivalent to the more natural choice
			\[ \norm{\rkl{v^+,v^-}}^2 =
			\norm{v^-}_{\frac 12,b^{-1}, \Sigma}^2
			+	 \norm{v^-}_{\frac 12,b^{-1}, \Sigma}^2
			+  \big\|(v^+ - v^-)\d^{-\nicefrac 12}\big\|_{0,\Sigma}^2, \]
			since $b^{-1}$ behaves like $\dist^{-\frac 12} \ge \d^{-\frac 12}$ near the crack tip.
			This yields that there exists a constant $C$ such that
			\begin{align*}
				\big \|(v^+ - v^-)\d^{-\nicefrac 12} \big\|_{0,\Sigma}
					&\le \big\|v^+ \d^{-\nicefrac 12}\big\|_{0,\Sigma}
						+ \big\|v^- \d^{-\nicefrac 12}\big\|_{0,\Sigma} \\
					&\le C \rkl{\norm{v^+}_{0,b^{-1},\Sigma}
						+ \norm{v^-}_{0,b^{-1},\Sigma}}.
			\end{align*}

		The next Lemma follows directly by the definitions of the
		spaces $V_{b^{-1}}$ and $W_{b^{-1}}$.
		\begin{lem}~ \label{lem:fluid_bulk_trace_fracture}
			\begin{enumerate}[(i)]
				\item The restriction of the trace operator
					$\gamma_\Sigma:V \rightarrow W_{\Sigma}$
					from Lemma \ref{lem:elasticity_trace_continuity}
					onto an operator from $V_{b^{-1}}$ to $W_{b^{-1}}$
					is well defined, linear and continuous. \label{lem:fluid_bulk_trace_fracture_i}
				\item The restriction of the extension operator
					$E_\Sigma: W_\Sigma \rightarrow V$ from Lemma \ref{lem:elasticity_trace_continuity} onto an operator
					from $W_{b^{-1}}$ to $V_{b^{-1}}$ is well defined, linear and continuous.
			\end{enumerate}
		\end{lem}
		\begin{proof}~
			\begin{enumerate}[(i)]
				\item By Lemma \ref{lem:fluid_bulk_trace} it follows that
				the Operator $\gamma_\Sigma$ is linear and continuous.
				By the definition of $V_{b^{-1}}$ (or more precisely: by
				the definition of $V$) it follows that $\jkl{v} \in H^{\frac 12}_{00}(\Sigma)$.
				\item By definition of $W_{b^{-1}}$ it follows that for
				$\widetilde v = \rkl{v^-, v^+} \in W_{b^{-1}}$ we have
				$v^\pm \in H^{\frac 12}_{b^{-1}}(\Sigma) \subset
				H^{\frac 12}(\Sigma)$. Then the extension of $\widetilde v$ from Lemma \ref{lem:elasticity_trace_continuity}
				satisfies $E_\Sigma \widetilde v \in V$. To show that
				$E_\Sigma \widetilde v \in V_{b^{-1}}$ it is required that
				$\gamma_\Sigma \rkl{E_\Sigma \widetilde v} \in H^{\frac 12}_{b^{-1}}\rkl{\Sigma}$. But this a direct consequence of
				Lemma \ref{lem:elasticity_trace_continuity}, since
				$\gamma_\Sigma \circ E_\Sigma = \text{id}$.
			\end{enumerate}
		\end{proof}

		\begin{thm} \label{thm:fluid_bulk_hilbert}
			The space $V_{b^{-1}}$ equipped with the norm $\norm{\cdot}_{1,b^{-1},\Omega}$ is a separable Hilbert space.
		\end{thm}
		\begin{proof}
			Obviously the norm $\norm{\cdot}_{1,b^{-1},\Omega}$ is induced by a scalar product.

			By definition $V$ is a Hilbert space and $V_{b^{-1}}$ is a subspace of $V$.
			The space $W_{b^{-1}}$ is a closed subspace of the space
			$H^{\frac 12}_{b^{-1}}{\Sigma} \times H^{\frac 12}_{b^{-1}}(\Sigma)$, which,
			by Theorem \ref{thm:fluid_volume_hilbert}, is a separable Hilbert space.
			The linearity and continuity of the trace operator $\gamma_\Sigma$ (Lemma \ref{lem:fluid_bulk_trace_fracture}) yields that
			$V_{b^{-1}}$ is a Banach space with respect to the graph norm
			\[ \norm{v}_G^2 \colonequals \norm{v}_{1,\Omega}^2
				+ \norm{\rkl{\gamma^+ v,\gamma^- v}}_{W_{b^{-1}}}^2. \]
			The continuity of the trace operators (Lemma \ref{lem:fluid_bulk_trace}) yields
			the equivalence of the norm $\norm{\cdot}_G$ and $\norm{\cdot}_{1,b^{-1},\Omega}$,
			and hence completeness of $V_{b^{-1}}$ with respect
			to the norm $\norm{\cdot}_{1,b^{-1},\Omega}$.
			Hence, the space $V_{b^{-1}}$ is a closed subspace of
			the separable Hilbert space $V$, which again implies the separability of $V_{b^{-1}}$.
		\end{proof}

		Note that from $\gamma^+v^2 + \gamma^-v^2 = 2 \skl{p}^2 + \frac 12 \jkl{p}^2$ an equivalent formulation of the norm $\| v \|_{1,b^{-1}, \Omega}^2$ is
		\[ \| v \|_{1,b^{-1}, \Omega}^2 =  \frac 12 \| \jkl{v} \|_{0,b^{-1},\Sigma}^2 +
			2\| \skl{v} \|_{0,b^{-1},\Sigma}^2 + \|v \|_{1,\Omega}^2. \]
	\subsection{A Weighted Sobolev Space on the Fracture}
		We now introduce a space representing the averaged fracture pressure $p^\Sigma$,
		namely
		\[
		 H^1_b (\Sigma) \colonequals \skl{ v \in L^2_{b^{-1}} (\Sigma) \mid \nabla_\tau v \in \mathbf L^2_{b}(\Sigma) }. \]
		This space is equipped with the norm
		 \begin{align*}
		  \| v \|_{1,b,\Sigma}^2 &:=   \| v \|_{0,b^{-1},\Sigma}^2  + \| \nabla_\tau v \|_{0,b,\Sigma}^2.
		 \end{align*}
		 \begin{thm} \label{thm:fluid_fracture_hilbert}
			 Let $\Sigma$ be a bounded hypersurface defined by the diffeomorphism
			 $\alpha:\Sigma \rightarrow \widehat \Sigma$ with boundary $\gamma = \partial \Sigma$, and let $b:\Sigma \rightarrow \R$ satisfy Assumptions A\ref{hyp:weight1}-A\ref{hyp:weight4}.
		   Then, the space $H^1_b(\Sigma)$ equipped with the norm $\| \cdot \|_{1,b,\Sigma}$ is a Hilbert space.
		 \end{thm}
		 \begin{proof}
		 	It is easy to see that the norm is induced by the scalar products for
		 	$L^2_{b^{-1}}\rkl{\Sigma}$ and $\mathbf L^2_b \rkl{\Sigma}$.
		 	We therefore only need to show that $H^1_{b}(\Sigma)$ is complete. We do this by showing that $H^1_{\widehat b}(\widehat \Sigma)$ with $\widehat b$ defined as in the proof of Lemma \ref{lem:ap-weight} is complete, and proceeding as in Theorem
		 	\ref{thm:wlebesgue-hilbert}. But this is a direct consequence of the fact that $\widehat b$ is a $A_2$ weight, since this yields that $\widehat b,\widehat b^{-1} \in L^1_{\text{loc}}(\R^{d-1})$ and hence that the space $H^1_{\widehat b}(\widehat \Sigma)$ is a Hilbert space (see, \ie,  \cite{kufner-opic-84}).
		 \end{proof}
 		 \begin{thm} \label{thm:fluid_fracture_density}
	 		 Let $\Sigma$ be a bounded hypersurface defined by the diffeomorphism
	 		 $\alpha:\Sigma \rightarrow \widehat \Sigma$ with boundary $\gamma = \partial \Sigma$, and let $b:\Sigma \rightarrow \R$ satisfy Assumptions A\ref{hyp:weight1}-A\ref{hyp:weight4}.
	 		 	The space $C^\infty_0(\overline \Sigma) \cap H^1_b(\Sigma)$ is dense in $H^1_b(\Sigma)$.
 		 \end{thm}
		 \begin{proof}
		 	We only need to show that $C_0^\infty (\widehat \Sigma) \cap H^1_{\widehat b}(\widehat \Sigma)$ is dense in $H^1_{\widehat b}(\widehat \Sigma)$. The assertion then follows directly as in the proof of Theorem \ref{thm:wlebesgue-hilbert}.

			 Our proof follows \cite{wloka} (Section 3.2).
			 It is easy to see that for any  $\Psi \in C^k(\widehat \Sigma)$, $k \in \mathbb N$ and $v \in H^1_{\widehat b}(\widehat \Sigma)$ we have $\Psi  v \in H^1_{\widehat b}(\widehat \Sigma)$.
			 Let $\eta \in C_0^\infty(\R^{d-1})$ be radial, decreasing and positive with $\int \eta = 1$. Define for $\varepsilon >0$ a 
			 function $\eta_\varepsilon (x) \colonequals \varepsilon^{-d} \eta (\nicefrac x \varepsilon)$ and
			 $v_\varepsilon = v \ast \eta_\varepsilon$.
			 For all $v \in L^2_{\widehat b^{\pm 1}} (\Sigma)$ we have $v_\varepsilon \in C_0^\infty(\R^{d-1})$ and $v_\varepsilon \rightarrow v$ in $L^2_{\widehat b^{\pm 1}} (\widehat \Sigma)$ \cite{miller}. Consequently, for any open
			  $\widehat \Sigma' \Subset \widehat \Sigma$ with $\text{dist}(\widehat \Sigma',\widehat \Sigma) > 0$ and $v \in H^1_{\widehat b}(\widehat \Sigma)$ we have $v_\varepsilon \rightarrow v$ in $H^1_{\widehat b}(\widehat \Sigma')$ (\cite{wloka}, Lemma 3.3).

			  We define $H^1_{\widehat b,0}(\widehat \Sigma)$ as the completion of $C^\infty_0(\widehat \Sigma)$-functions
			  with respect to the $H^1_{\widehat b}(\widehat \Sigma)$ norm.
			  Then $H^1_{\widehat b,0}(\widehat \Sigma)$ is a closed subspace of $H^1_{\widehat b}(\widehat \Sigma)$ and
			  therefore a separable Hilbert space.
			  Together with the density property on compact subsets this definition yields that all functions with bounded support in $H^1_{\widehat b}(\widehat \Sigma)$ are dense in $H^1_{\widehat b}(\widehat \Sigma)$ (\cite{wloka}, Theorem 3.3).
			  With that we can show that
			  the set $H^1_{\widehat b}(\widehat \Sigma) \cap C^\infty (\overline{\widehat \Sigma})$ is dense in
			  $H^1_{\widehat b}(\widehat \Sigma)$ (\cite{wloka}, Theorem 3.4 and 3.5).
		 \end{proof}
		\begin{thm}[Poincar\'e--Friedrichs inequality] \label{thm:fluid_fracture_poincare}
			Let $\Sigma$ be a bounded hypersurface defined by the diffeomorphism
			$\alpha:\Sigma \rightarrow \widehat \Sigma$ with boundary $\gamma = \partial \Sigma$, and let $b:\Sigma \rightarrow \R$ satisfy Assumptions A\ref{hyp:weight1}-A\ref{hyp:weight4}.
			Then there exists a positive constant $C$ such that
			for any $u \in C_0^\infty\rkl{\Sigma}$ we have
			\[ \norm{u}_{0,b^{-1},\Sigma} \le C \norm{\nabla_\tau
				u }_{0,b, \Sigma}^2. \]
		\end{thm}
		\begin{proof}
			We will show the corresponding result on $\widehat \Sigma$. A simple coordinate transformation then yields the result. Let $\widehat b$ be
			defined as in the proof of Lemma \ref{lem:ap-weight}.

			Lemma 4.1 in \cite{edmunds-opic} says that if
			 $H^1_{\widehat b,0}(\widehat \Sigma)$ is compactly embedded in
			 $L^1_{\widehat b^{-1}}(\widehat \Sigma)$ and if there exist $x_0 \in \partial \widehat \Sigma$
			 and $R>0$ such that
			 $\widehat \Sigma(x_0,R) \colonequals B_R(x_0) \cap M$ is a Lipschitz domain and $\widehat b, \widehat b^{-1} \in L^1 \rkl{\Sigma(x_0,R)}$, then
			 the weighted Friedrichs inequality
			 \[ \int_{\widehat \Sigma} |v|^2 \widehat{b}^{-1} \intd \mathbf s \le K  \int_{\widehat \Sigma} |  \nabla v|^2 \widehat b \intd \mathbf s  \]
			 holds for all $v \in C_0^\infty (\widehat \Sigma)$. Herein 
			 $\nabla$ denotes the gradient operator on $\R^{d-1}$.

			It it easy to see that the embedding of $H^1_{\widehat b,0}(\widehat \Sigma)$ into
			$L^1_{\widehat{b}^{-1}}(\widehat \Sigma)$
			is continuous. To show that the embedding is also compact we introduce the
			 space
			\[ H_{\widehat b \widehat b}^1(\widehat \Sigma) \colonequals \skl{ v \in L^2_{\widehat b}(\widehat \Sigma) \mid \widehat\nabla v \in \mathbf L^2_{\widehat b}(\widehat \Sigma) }, \]
			with norm
			\[  \norm{ v }_{1,\widehat b \widehat b,\widehat \Sigma}^2
			\colonequals \norm{v}_{0,\widehat b,\widehat \Sigma}^2
			+ \| \nabla v \|_{0,\widehat b,\widehat \Sigma}^2. \]
			Again, it is easy to see that this norm is induced by a scalar product.
			The space $H_{\widehat b \widehat b}^1(\widehat \Sigma)$ equipped with this norm
			is a Hilbert space \cite{kufner-opic-84}. Furthermore, by
			Theorem 19.11 in \cite{kufner-opic-90} it follows that
			$ H_{\widehat b \widehat b}^1(\widehat \Sigma)$ is compactly embedded in $L^2_{\widehat{b}^{-1}}(\widehat \Sigma)$.

			Let $(v_k)_{k \in \mathbb N}$ be a bounded sequence in $H_{\widehat b,0}^1(\widehat \Sigma)$.
			Since $\norm{v_k}_{0,\widehat b,\widehat \Sigma}^2  \le C\norm{v_k}_{0,\widehat{b}^{-1},\widehat \Sigma}$, the sequence $(v_k)_{k \in \mathbb N}$ is bounded in
			$H_{\widehat b \widehat b}^1(\widehat \Sigma)$.
			Hence there exists a subsequence $\rkl{v_{k_j}}_{j \in \mathbb N}$ and a function $v$ in $L^2_{\widehat{b}^{-1}}(\widehat \Sigma)$ with
			$s_{k_j} \rightarrow s$ in $L^2_{\widehat{b}^{-1}}(\widehat \Sigma)$.
			This proves that $H_{\widehat b}^1(\widehat \Sigma)$ is compactly
			embedded in $L^2_{\widehat{b}^{-1}}(\widehat \Sigma)$.

			Since $\widehat \Sigma$ is a bounded domain with Lipschitz boundary, its
			boundary is locally convex, \ie there exists a neighborhood
			$\mathcal U$, where $\partial \widehat \Sigma$ can be represented as convex
			graph of a Lipschitz continuous function. Take $x_0 \in \mathcal U \cap \partial \widehat \Sigma$ and choose $R$ such that
			$B_R(x_0) \subset \mathcal U$, then $B_R(x_0) \cap \Sigma$ is
			convex and therefore a Lipschitz domain.
			By Lemma \ref{lem:ap-weight} it follows that $\widehat b$ is an $A_2$-weight.
			This yields $\widehat b, \widehat{b}^{-1} \in L^1_{\text{loc}}(B)$ for all balls $B$ contained in $\mathbb R^{d-1}$. Hence we have $\widehat b, \widehat{b}^{-1} \in  L^1_{\text{loc}}(B(x_0,R))$.
			This proves the assertion, since
			\begin{align*}
				\norm{v}_{0,b^{-1},\Sigma} &=
					\norm{v \circ \alpha^{-1}}_{0,\widehat b^{-1},\widehat \Sigma}\\
				&\le K \| \nabla \rkl{v \circ \alpha^{-1}}\|_{0,\widehat b,\Sigma}  \\
				&= K\| \rkl{\nabla_\tau (v)\circ \alpha^{-1}}  \nabla \alpha^{-1} \|_{0,\widehat b,\widehat \Sigma} \\
				&\le K\|  \nabla \alpha^{-1} \|_{L^\infty(\widehat \Sigma)}
					\norm{\nabla_\tau v}_{0,b,\Sigma}\\
				&= C \norm{\nabla_\tau v}_{0,b,\Sigma}.
			\end{align*}
		\end{proof}
		\begin{rem} \label{rem:weighted_trace_fracture}
			We did not prove any trace and extension theorems concerning
			the space $H^1_b(\Sigma)$, but want to mention the following result from
			\cite{brewster-mitrea}. Define the weighted Sobolev space
			\[ H^1_{b^{-1}b^{-1}} (\Sigma) \colonequals
			\skl{v \in L^2_{b^{-1}}(\Sigma) \mid \nabla_\tau v \in \mathbf L^2_{b^{-1}}(\Sigma) }^2 \]
			with norm $\norm{v}_{b^{-1},b^{-1},\Sigma}^2
			\colonequals \norm{v}_{0,b^{-1},\Sigma}^2 + \norm{\nabla_\tau v}_{0,b^{-1},\Sigma}$.
			Then $H^{1}_{b^{-1}b^{-1}}(\Sigma) \hookrightarrow H^1_b(\Sigma)$. Furthermore denote by $H^{\frac 34}(\gamma)$ the usual Sobolev-Slobodeckij space on $\gamma$. Then there exists a linear and continuous trace operator
			$\gamma_{b^{-1}}: H^1_{b^{-1}b^{-1}} (\Sigma) \rightarrow H^{\frac 34} (\gamma)$ and a linear and continuous extension operator
			$E_{b^{-1}} : H^{\frac 34}(\gamma) \rightarrow  H^1_{b^{-1}b^{-1}} (\Sigma)$,
			such that $\rkl{\gamma_{b^{-1}} \circ E_{b^{-1}}} f = f$ for all
			$f \in  H^{\frac 34}(\gamma)$.
		\end{rem}
\section{Weak Formulation} \label{sec:weak_formulation}
  Using the Sobolev spaces introduced in the previous section, we will derive weak formulations of the fluid problem in the bulk, the averaged fluid problem in the fracture and the poroelasticity equation.
  For fixed crack width functions $b$
  we prove existence and uniqueness of solutions to the coupled fluid--fluid problem. Finally,
  the weak formulation of the coupled problem is introduced.

\subsection{Weak Formulations of the Subdomain Problems}
We first derive the weak formulation of the three subproblems individually.
\subsubsection{The Weak Elastic Problem} \label{subsubsec:weak-elasticity}

Consider the problem defined by the equation \eqref{all:biot} together with
boundary conditions \eqref{poroelasticity:neumannBC}, \eqref{poroelasticity:dirichletBC} and \eqref{all:coupling-stress}.
Suppose that $p^\Omega$ and $p^\Sigma$ are fixed. Then \eqref{all:biot} is the linear elasticity equation with particular volume terms and boundary conditions.
Multiplication of equation \eqref{all:biot} with a test function $\mathbf v \in \mathbf V_0$, integration over $\Omega$ and
the componentwise application of the generalized Green formula (Lemma \ref{lem:elasticity_generalized_green}) yields
\begin{align*}
	 \int_{\Omega} \sigma\rkl{\mathbf u} : \mathbf e \rkl{\mathbf v} \intd x
		&= \rkl{\mathbf f_E, \mathbf v}_{0,\Omega}
			+ \rkl{\sigma \rkl{\mathbf u}\cdot \mathbf n, \mathbf v}_{0,\Gamma_N^E}
			- \rkl{\nabla p^\Omega, \mathbf v}_{0,\Omega}
			- \rkl{\jkl{\sigma\rkl{\mathbf u}}\cdot \nu, \skl{\mathbf v}}_{0,\Sigma} \\
		& \quad - \rkl{\skl{\sigma\rkl{\mathbf u}}\cdot \nu, \jkl{\mathbf v}}_{0,\Sigma}. 
\end{align*} 
Inserting the coupling and boundary conditions \eqref{all:coupling-stress} and
\eqref{poroelasticity:neumannBC} then yields the equation
\[ \int_{\Omega} \sigma\rkl{\mathbf u} : \mathbf e \rkl{\mathbf v} \intd x
=\rkl{\mathbf f_E, \mathbf v}_{0,\Omega}
+ \rkl{\sigma_N, \mathbf v}_{0,\Gamma_N^E} - \rkl{\nabla p^\Omega, \mathbf v}_{0,\Omega}
 + \rkl{p^\Sigma, \skl{\mathbf v}}_{0,\Sigma}. \]
Define the affine space
\[	\mathbf V_E  \colonequals \big \{ \mathbf u \in \mathbf V \mid
			\mathbf u \vert_{\Gamma_D^E \setminus \widetilde \Sigma} =  \mathbf u_D \big\} \]
and the bilinear form $a_E: \mathbf V_E \times \mathbf V_E \rightarrow \R$ by
\[ a_E\rkl{\mathbf u, \mathbf v} \colonequals \int_{\Omega} \sigma(\mathbf u) : \mathbf e(\mathbf v) \intd x. \]
Furthermore, introduce the linear forms $c_{E,p^\Omega}: \mathbf V_E  \rightarrow \R$ and $c_{E,p^\Sigma}: \mathbf V_E  \rightarrow \R$ by
\begin{equation*}
c_{E,p^\Omega} (\mathbf v) \colonequals \int_\Omega \nabla p^\Omega \cdot \mathbf v \intd x
\quad \text{and} \quad
c_{E,p^\Sigma} (\mathbf v) \colonequals \int_\Sigma p^\Sigma \jkl{\mathbf v} \cdot \nu \intd a,
\end{equation*}
and $l_E^\Omega: \mathbf V_E \to \mathbb R$ by
\[
l_E^\Omega (\mathbf v) = \int_\Omega \mathbf f_E \cdot \mathbf v\intd x
+ \int_{\Gamma_N^E} \sigma_N \cdot \mathbf v \intd s.
\]
Define $l_{E,p^\Sigma,p^\Omega}: \mathbf V_E \rightarrow \mathbb R$ by
\[
l_{E,p^\Sigma,p^\Omega} (\mathbf v) \colonequals l_E^\Omega(\mathbf v)  -c_{E,p^\Omega}(\mathbf v) + c_{E,p^\Sigma}(\mathbf v).
\]

Then the elasticity problem is formally equivalent to the following weak formulation:
For given $p^\Omega \in V_b$ and $p^\Sigma \in H^1_b(\Sigma)$, find $\mathbf u \in \mathbf V_E$ such that
\begin{equation}
\label{eq:weak_elasticity_problem}
a_E(\mathbf u, \mathbf v) = l_{E,p^\Sigma,p^\Omega}(\mathbf v)
\qquad \forall \, \mathbf v \in \mathbf V_0.
\end{equation}

 This problem is a standard linear elasticity problem in fractured domains.
 It has a unique solution, which depends continuously on the data $p^\Omega$ and
 $p^\Sigma$ from the fluid problems.
 A proof can be found in \cite{khludnev}.
 \begin{thm} \label{thm:existence-solid}
 	Let $\Sigma$ be of class $C^2$, and
 	suppose that $\mu \ge 0$ and $\lambda + 2\mu \ge 0$. Assume that the Dirichlet boundary
 	$\Gamma_D^E$ has a positive $(d-1)$-dimensional measure. Let $p^\Sigma \in H^1_b(\Sigma)$
 	and $p^\Omega \in V_b$, and assume that $\mathbf u_D = \mathbf 0$ and
 	$\sigma_N = \mathbf 0$. Then
 	there exists a unique solution $\mathbf u \in \mathbf V_E$ to problem \eqref{eq:weak_elasticity_problem}.
 	     Furthermore the estimate
 	     \[  \norm{\mathbf{u}}_{\mathbf V_E} \le \frac{C}{2\mu} \rkl{\norm{\mathbf f_E}_{0,\Omega} + \norm{\nabla p^\Omega}_{0,\Omega} + \norm{p^\Sigma}_{0,\Sigma}} \]
 	     is satisfied for a constant $C$ depending only on the domain $\Omega$.
 \end{thm}
We now want to gain a better understanding of the behavior of the solution near the crack tip.
It is well known in the literature that the weak solution of a linear elastic problem
in a fractured domain can be decomposed into a singular term $\mathbf u_s$ and a regular term $\mathbf u_r$, where the singular term describes the discontinuous behavior of the solution across the
fracture and the singular behavior of the stress nearby the crack front \cite{mazia-rossmann-kozlov,costabel-dauge-edge-asymptotics}.
To understand the behavior of the solution near a single the crack front we introduce a polar coordinate frame $\rkl{r,\Theta}$ with $\Theta = 0$ along the tangential extension of the fracture and $r \colonequals \text{dist}(\mathbf x, \gamma)$ for all $x \in \R^d$. In three dimensions we extend this frame to a cylindrical coordinate frame by introducing a parametrization of $\gamma$ and using the arc length parameter as third coordinate.

Rewriting equations \eqref{all:biot}, \eqref{all:coupling-stress} in polar or cylindrical coordinates, and applying
the Mellin transform
\[ \widehat u (\alpha;\Theta) \colonequals \mathcal M [u(r,\Theta)](\alpha) = \int_{\mathbb R^+} r^{-\alpha - 1} u(r,\Theta)
\intd r \]
to the resulting problem, generates a differential equation with parameter $\alpha \in \mathbb C$. All $\alpha$ for which the resulting problem with zero right-hand-side admits nontrivial solutions $\mathbf v_\alpha(\Theta)$ are called eigenvalues, and the corresponding functions $\mathbf v_\alpha(\Theta)$ are called eigenvectors. These eigenvectors determine the singularity of the problem.

\begin{thm} Let $\Omega \subset \R^2$ or $\Omega \subset \R^3$.
	\begin{enumerate}[(i)] \label{cor:asymptotics}
		\item \label{cor:asymptotics_1} The eigenvalues of the elasticity problem with symmetric Neumann boundary conditions on the fracture are given by $\Lambda_E = \{ \nicefrac k2 \mid k \in \mathbb Z \}$.
		\item \label{cor:asymptotics_2} The Mellin transform is well defined and invertible for all $\alpha$ with $\text{Re}(\alpha) \notin \Lambda_E$. Hence $\widehat u(\alpha,\Theta)$ is holomorphic
		for all $\alpha$ with $\text{Re}(\alpha) \in \R \setminus \Lambda_E$.
		\item The singular part $\mathbf u_s$ of the solution to the elasticity problem with 
		symmetric Neumann boundary conditions on the fracture 
		is of the form
			\[  \mathbf u_{s} = \chi \rkl{ K_{\frac 12}^{(1)} \mathbf v_{\frac 12}^{(1)}
				+ K_{\frac 12}^{(2)} \mathbf v_{\frac 12}^{(2)}
				+ K_{\frac 12}^{(3)} \mathbf v_{\frac 12}^{(3)}} \]
			near the crack front.
			The real constants $K_{\frac 12}^{(i)}$ are called
			\emph{stress intensity factors} and $\chi$ is a smooth cutoff-function, which is equal to one in a neighborhood of the crack front.
			The components of the eigenfunctions $\mathbf v_{\frac 12}^{(i)}$ are in the span of the set
			\[ \widetilde{\mathcal S} =  \big\{  \sin(\nicefrac \Theta 2),
			 \cos(\nicefrac \Theta 2),  \sin(\Theta)\sin(\nicefrac \Theta 2),  \sin(\Theta)\cos(\nicefrac \Theta 2) \big\}.\]
	\end{enumerate}
\end{thm}
The important conclusion from the latter Theorem is, that the fracture width function $b$ can be approximated nearby the crack tip by
\[ b = \jkl{\mathbf u_s} \cdot \nu \simeq r^{\frac 12}
	= \text{dist}(\cdot, \gamma)^{\frac 12}. \]
We want to mention, that we only considered edge asymptotics in the three-dimensional case here, since $\Sigma$ is smooth enough. If $\Sigma$ only has a Lipschitz boundary, then vertex singularities have to be considered as well.

\subsubsection{The Fluid Bulk Problem}
Assume next that the width $b \in H^{\frac 12}_{00}\rkl{\Sigma}$ and the fracture fluid pressure $p^\Sigma \in H^1_b(\Sigma)$
are fixed.
We consider the fluid bulk problem introduced by the equations \eqref{all:darcy-mass}, \eqref{all:darcy-momentum} together with the
  boundary conditions \eqref{darcy:neumannBC}, \eqref{darcy:dirichletBC} and the coupling conditions \eqref{all:average-pressure-approximation},
 \eqref{all:averaged-normal-momentum} and \eqref{all:coupling-width}.
  Let
  \[ V_F^\Omega \colonequals \skl{ p^\Omega \in V_b\mid
  	p^\Omega\vert_{\Gamma_D^F \setminus \widetilde \Sigma} =  p_D^\Omega }
   \qquad \text{and} \qquad
   V_{b,0} =  \skl{ r^\Omega \in V_b \mid r^\Omega\vert_{\Gamma_D^F \setminus \widetilde \Sigma} =  0}. \]

 We multiply equation \eqref{all:darcy-mass} with a test function and apply the generalized Green's formula (see Lemma \ref{lem:elasticity_generalized_green}) to derive the equation
 \begin{align*}
 	\rkl{f_F^\Omega, r^\Omega}_{0,\Omega} -
 	\rkl{\mathbf q^\Omega \cdot \mathbf n, r^\Omega}_{0,\Gamma_N^F} =
 	 - \rkl{\mathbf q^\Omega, \nabla r^\Omega}_{0,\Omega} - \rkl{\jkl{\mathbf q^\Omega} \cdot \nu, \skl{r^\Omega}}_{0,\Sigma}
 	   -  \rkl{\skl{\mathbf q^\Omega}\cdot \nu,
 	 	\jkl{r^\Omega}}_{0,\Sigma}.
 \end{align*}
	Replacing the bulk fluid velocity jump $\jkl{\mathbf q^\Omega} \cdot \nu$
	and average $\skl{\mathbf q^\Omega} \cdot \nu$ by the
  corresponding pressure values in the coupling conditions \eqref{all:average-pressure-approximation} and
  \eqref{all:averaged-normal-momentum} yields the equation
  \begin{align*}
  	\rkl{f_F^\Omega, r^\Omega}_{0,\Omega} -
  	\rkl{\mathbf q^\Omega \cdot \mathbf n, r^\Omega}_{0,\Gamma_N^F} & =
  	 - \rkl{\mathbf q^\Omega, \nabla r^\Omega}_{0,\Omega}
  	 - \frac{1}{2 \xi - 1} \rkl{\frac{4 K^\nu}{b} \rkl{p^\Sigma - \skl{p^\Omega}}, \skl{r^\Omega}}_{0,\Sigma} \\
  	 & \quad 	+  \rkl{\frac{K^\nu}{b} \jkl{p^\Omega}, \jkl{r^\Omega}}_{0,\Sigma}.
  \end{align*}
  Inserting the momentum equation \eqref{all:darcy-momentum} and the Neumann boundary conditions
  \eqref{darcy:neumannBC}, we derive the weak problem: For given $b \in H^{\frac 12}_{00}(\Sigma)$ and $p^\Sigma \in H^1_b(\Sigma)$, find
  $p^\Omega \in V_F^\Omega$ such that
\begin{equation} \label{weak:fluid_bulk}
 a_{F,b}^\Omega(p^\Omega, r^\Omega) - c_{F,b,p^\Sigma}(r^\Omega) = l_F^\Omega(r^\Omega) \qquad \forall \, r^\Omega \in V_{b,0}.
\end{equation}
Here, the bilinear form $a_F^\Omega :  V_F^\Omega \times V_F^\Omega \rightarrow \R$ is defined by
\begin{equation*}
a_{F,b}^\Omega (p^\Omega, r^\Omega)
\colonequals
\sprod{\mathbb{K} \nabla p^\Omega}{\nabla r^\Omega}_{0,\Omega}
	+ \frac{1}{2 \xi - 1} \sprod{\frac{4K^\nu}{b} \skl{p^\Omega}}{\skl{r^\Omega}}_{0,\Sigma}
	 + \sprod{\frac{K^\nu}{b}\jkl{p^\Omega}}{\jkl{r^\Omega}}_{0,\Sigma}.
\end{equation*}
The linear form $c_{F,b,p^\Sigma} : V_F^\Omega \rightarrow \R$,
\begin{equation*}
 c_{F,b,p^\Sigma} (r^\Omega) \colonequals \frac{1}{2\xi - 1} \sprod{\frac{4K^\nu}{b} p^\Sigma}{\skl{r^\Omega}}_{0,\Sigma}
\end{equation*}
represents the coupling to the fracture fluid pressure $p^\Sigma$.
Finally, the linear form  $l_F^\Omega: V_F^\Omega \rightarrow \R$
\begin{equation*}
  l_F^\Omega (r^\Omega) \colonequals \sprod{f_F^\Omega}{r^\Omega}_{0,\Omega} - \sprod{ q_N^\Omega}{r^\Omega}_{0,\Gamma_N^F},
\end{equation*}
represents the external source terms.

  The map $a_{F,b}^\Omega$ is symmetric and bilinear in $p^\Omega$ and $r^\Omega$, but depends nonlinearly on the fracture width
   $b = \jkl{\mathbf u}\cdot \nu$. Consequently, it depends nonlinearly on the displacement field $\mathbf u$. 
   We postpone the discussion on the existence of weak solutions to Section
   \ref{subsec:fluid-existence}.
\subsubsection{Fracture Fluid Problem}
Suppose now that the fracture width $b \in H^{\frac 12}_{00}\rkl{\Sigma}$ and the bulk fluid pressure $p^\Omega \in V_F^\Omega$
are fixed.
Define
\[	V_F^\Sigma  \colonequals \skl{ p^\Sigma \in H^1_b(\Sigma) \mid p^\Sigma \vert_{\gamma_D^F} = p_D^\Sigma}
   \qquad \text{and} \qquad
 H^1_{b,0} (\Sigma) \colonequals
\skl{ r^\Sigma \in H^1_{b}(\Sigma) \mid r^\Sigma \vert_{\gamma_D^F} = 0}.  \]
Consider the averaged fluid problem on the fracture, given by Equations
  \eqref{all:averaged-mass} and \eqref{all:averaged-tangential-momentum}, together with the boundary conditions \eqref{fracture:neumannBC}, and
  \eqref{fracture:dirichletBC} and the coupling conditions \eqref{all:coupling-width} and \eqref{all:average-pressure-approximation}.
  Insert \eqref{all:averaged-tangential-momentum}
  into \eqref{all:averaged-mass}, and multiply it with a test function $r^\Sigma \in V_F^\Sigma$. We integrate the resulting
  Laplace--Beltrami-like problem over $\Sigma$ and apply integration by parts, to obtain
  \[  \rkl{ b K^\tau \nabla p^\Sigma, \nabla r^\Sigma}_{0,\Sigma}
  		= \rkl{b f^\Sigma, r^\Sigma}_{0,\Sigma} - \rkl{\jkl{\mathbf q^\Omega} \cdot \nu, r^\Sigma}_{0,\Sigma} - \rkl{\kappa b \skl{\mathbf q^\Omega}\cdot \nu, r^\Sigma}_{0,\Sigma}
  			- \rkl{b q_N^\Sigma, r^\Sigma}_{0,\gamma_N^F}. \]
  Using equation  \eqref{all:average-pressure-approximation} to replace the
  fluid velocity normal jump by a term only depending on $p^\Sigma$ and $p^\Omega$ yields the following problem:
  For given $b \in H^{\frac 12}_{00}(\Sigma)$ and $p^\Omega \in V_F^\Omega$, find $p^\Sigma \in V_F^\Sigma$, such that
\begin{equation} \label{weak:fluid_fracture}
 a_{F,b}^\Sigma(p^\Sigma, r^\Sigma) - c_{F,b,p^\Omega}(r^\Sigma) = l_{F,b}^\Sigma(r^\Sigma) \qquad \forall \, r^\Sigma \in H^1_{b,0}(\Sigma).
\end{equation}
Here, the bilinear form $a_{F,b}^\Sigma : V_F^\Sigma \times V_F^\Sigma$ is defined by
\begin{equation*}
 a_{F,b}^\Sigma (p^\Sigma, r^\Sigma)
  \colonequals
  \sprod{ b K^\tau \nabla_\tau p^\Sigma}{\nabla_\tau r^\Sigma}_\Sigma
    + \frac{1}{2\xi -1}\sprod{\frac{4 K^\nu}{b} p^\Sigma}{r^\Sigma}_{0,\Sigma},
\end{equation*}
the linear form $c_{F,b,p^\Omega} : V_F^\Sigma \rightarrow \R$ by
 \begin{equation*}
 	c_{F,b,p^\Omega} (r^\Sigma) \colonequals \frac{1}{2\xi - 1} \sprod{\frac{4K^\nu}{b} r^\Sigma}{\skl{p^\Omega}}_{0,\Sigma}
	 	+ \rkl{\kappa K^\nu \jkl{p^\Omega}, r^\Sigma}_{0,\Sigma},
 \end{equation*}
 and the source term $l_F^\Sigma: \mathbf V_F^\Sigma \to \R$ is defined by
  \[ l_{F,b}^\Sigma (r^\Sigma)
      = \sprod{b f^\Sigma}{r^\Sigma}_{0,\Sigma}
      	- \sprod{b q_N^\Sigma}{r^\Sigma}_{0,\gamma_N^F}. \]

  The map $a_{F,b}^\Sigma$ is symmetric and bilinear in $p^\Sigma$ and $r^\Sigma$, but depends nonlinearly on $\mathbf u$ through
  the fracture width $b = \jkl{\mathbf u}\cdot \nu$. The map $l_{F,b}^\Sigma$
  is linear in $r^\Sigma$, but also depends linearly on $\mathbf u$. 
  Note that the curvature $\kappa$ only plays a role if $p^\Omega$ does jump across the fracture.
  Existence and uniqueness of solutions to this problem is discussed
  in the next section.
\subsection{Existence of Solutions of the Weak Fluid--Fluid Problem} \label{subsec:fluid-existence}
For a fixed crack width function $b$, the coupled fluid problems
\eqref{weak:fluid_bulk} and \eqref{weak:fluid_fracture} form a joint linear variational problem.
From the asymptotic expansion of the displacement field $\mathbf u$
(Section \ref{subsubsec:weak-elasticity}) we know that $b$ behaves like $\dist^\frac 12(\cdot, \gamma)$ near the crack tip. This allows us to show existence of unique solutions to the coupled fluid--fluid problem in suitable weighted Sobolev spaces, under the assumption that the crack is open.

To this end, we define the bilinear forms
$c_{F,b}:V_F^\Sigma \times V_F^\Omega \rightarrow \R$ by
\[ c_{F,b}(r^\Sigma, r^\Omega) \colonequals c_{F,b,r^\Sigma}(r^\Omega) = c_{F,b,r^\Omega}(r^\Sigma) - \sprod{\kappa K^\nu \jkl{r^\Omega}}{r^\Sigma}_{0,\Sigma} \]
and $c_{F,b,\kappa} :V_F^\Sigma \times V_F^\Omega \rightarrow \R$ by
\[c_{F,\kappa}(r^\Sigma, r^\Omega) \colonequals  \sprod{\kappa K^\nu \jkl{r^\Omega}}{r^\Sigma}_{0,\Sigma}.  \]
For given $b \in H^{\frac 12}_{00}(\Sigma)$ combining the weak formulations 
\eqref{weak:fluid_bulk} and \eqref{weak:fluid_fracture} yields: find $p^\Sigma \in V_F^\Sigma$ and $p^\Omega \in V_F^\Omega$ such that
\begin{subequations}
	\label{weak:coupled-fluid}
	\begin{alignat}{2}
	\label{eq:weak-coupled-fluid-bulk}
	a_{F,b}^\Omega(p^\Omega, r^\Omega) - c_{F,b}(p^\Sigma, r^\Omega) &= l_F^\Omega(r^\Omega) & \qquad & \forall \, r^\Omega \in V_{b,0},\\
	\label{eq:weak-coupled-fluid-fracture}
	a_{F,b}^\Sigma(p^\Sigma, r^\Sigma) - c_{F,b}(r^\Sigma, p^\Omega)
	- c_{F,b,\kappa}(r^\Sigma,p^\Omega) &= l_{F,b}^\Sigma(r^\Sigma) && \forall \, r^\Sigma \in H^1_{b,0}(\Sigma).
	\end{alignat}
\end{subequations}
 We now prove existence and uniqueness of a
solution of this weak coupled fluid--fluid problem.

\begin{thm} \label{thm:existence_fluid}
	Assume that $\abs{\kappa} \le \kappa_{\text{max}} < \infty$ and that $\mathbb K$ is symmetric, bounded and uniformly elliptic, \ie, there exists
	a constant $K > 0$ such that
	\[ x^T \mathbb K x > K \norm{x}^2 \qquad \forall \, x \in \R^d \quad \text{and} \quad
	\text{a.e. on } \Omega. \]
	Furthermore, assume that $K^\tau, K^\nu$ are positive constants.
	Let $b$ satisfy Assumptions \ref{hyp:weight} and let $\xi \in \big (\frac 12,1 \big]$. Suppose that the curvature of the fracture is bounded in the sense that 
	\[ b_{\text{max}}^2 \kappa_{\text{max}}^2 < \frac{4 K^\tau}{K^\nu} C_\Sigma, \]
	where $C_\Sigma$ denotes the Poincar\'e constant of $\Sigma$ from Theorem \ref{thm:fluid_fracture_poincare}.
	Let $\Gamma_D^F \neq \emptyset$, $f_F^\Omega \in L^2(\Omega)$, $f_F^\Sigma \in L^2_+(\Sigma)$, $q_N^\Omega \in L^2(\Gamma_N)$, and $q_N^\Sigma \in L^2_+(\Sigma)$. Finally, assume that
	\[ p_D^\Omega \in W_D^\Omega \colonequals \skl{s \in H^{\frac 12}(\Gamma_D^F) \mid
		E_\Omega s \in V_b},
	\qquad
	p_D^\Sigma \in W_D^\Sigma \colonequals
	\skl{s \in H^{\frac 12}(\gamma_D^F) \mid
		E_\Sigma s \in H^1_b(\Sigma)}, \]
	where $E_\Omega : H^{\frac 12}(\Gamma_D^F) \rightarrow V$ and
	$E_\Sigma: H^{\frac 12}(\gamma_D^F) \rightarrow H^1(\Sigma)$
	are the standard extension operators.
	Then
	there exists a unique solution $(p^\Sigma,p^\Omega) \in V_F^\Sigma \times V_F^\Omega$ of the weak coupled problem \eqref{weak:coupled-fluid}.
\end{thm}
Before presenting the proof, note first that the trace space $W_D^\Omega$ is not empty. From the trace
theorem \ref{lem:fluid_bulk_trace}
we know that the trace of a
$V_b$-function restricted to $\Gamma_D^F$ is in $H^{\frac 12}(\Gamma_D^F)$. By Lemma \ref{lem:fluid_bulk_trace_fracture}
we conclude that the trace of this $V_b$-function
restricted to $\Sigma$ is in $W_b$. And thus, again by Lemma
\ref{lem:fluid_bulk_trace_fracture}, it follows that the application
of the global extension operator yields the identity and thus a
$V_b$-function.

Similarly, the space $W_D^\Sigma$ is not empty. By Remark \ref{rem:weighted_trace_fracture} each function in $H^{\frac 34}(\gamma_D)$ has an extension in $H^1_{-b,-b}(\Sigma)$ and thus
in $H^1_b(\Sigma)$.
\begin{proof}
	Without loss of generality we assume that $p_D^\Omega =0$ and $p_D^\Sigma = 0$.
	Then $V_F \colonequals  V_F^\Sigma \times V_F^\Omega = H^1_{b,0}(\Sigma) \times
	V_{b,0}$, and we equip this space
	with the norm $\norm{(p^\Sigma,p^\Omega)}_{V_F}^2 \colonequals \norm{p^\Sigma}_{1,b,\Sigma}^2 + \norm{p^\Omega}_{1,b^{-1},\Omega}^2$.
	Define the bilinear form $k_b: V_F \times V_F \rightarrow \R$ by adding the left hand sides of \eqref{eq:weak-coupled-fluid-bulk} and \eqref{eq:weak-coupled-fluid-fracture}
	\[ k_b(p,r) \colonequals a_{F,b}^\Omega(p^\Omega, r^\Omega) + a_{F,b}^\Sigma (p^\Sigma, r^\Sigma) - c_{F,b} (p^\Sigma, r^\Omega)- c_{F,b} (r^\Sigma, p^\Omega) - c_{F,b,\kappa}(r^\Sigma,p^\Omega).   \]
	Likewise, define the linear form $l_b:V_F \rightarrow \R$ by adding the right hand sides of
	 \eqref{eq:weak-coupled-fluid-bulk} and \eqref{eq:weak-coupled-fluid-fracture}
	\[ l_b(r) = l_F^\Omega (r^\Omega) + l_{F,b}^\Sigma (r^\Sigma), \]
	where $p = \rkl{ p^\Sigma, p^\Omega} \in V_F$ and $r = \rkl{r^\Sigma, r^\Omega} \in V_F$.
	It is easy to check that each solution $p \in V_F$ of the problem
	\begin{equation} \label{weak:summed_problem}
	k_b(p,r) = l(r) \qquad \forall \, r \in V_{F}
	\end{equation}
	is a weak solution of problem \eqref{weak:coupled-fluid} and vice versa.

	We will now prove existence and uniqueness of a solution to \eqref{weak:summed_problem} using the Lax--Milgram Lemma. The continuity of the linear form $l$ can be easily shown using the
	Cauchy--Schwarz inequality.
	Similarly, boundedness of the bilinear form $k_b$ can be shown using the Cauchy--Schwarz inequality, the boundedness of the permeability tensors, and the continuity of the trace operators.

	To prove that $k$ is coercive, we introduce the map $g_b:V_F \to \R$
	\begin{align*}
	g_b(p) &\colonequals \frac{4}{2 \xi - 1} \bigg[ \sprod{K^\nu b^{-1} \skl{p^\Omega}}{\skl{p^\Omega}}_{0,\Sigma}
	+  \sprod{ K^\nu b^{-1} p^\Sigma}{p^\Sigma}_{0,\Sigma}
	- 2\sprod{K^\nu b^{-1} p^\Sigma}{\skl{p^\Omega}}_{0,\Sigma}  \, \bigg] \\
	& \qquad + \sprod{ K^\nu b^{-1} \jkl{p^\Omega}}{\jkl{p^\Omega}}_{0,\Sigma}
	- \sprod{\kappa K^\nu \jkl{p^\Omega}}{p^\Sigma}_{0,\Sigma},
	\end{align*}
	and we note that
	\[ g_b(p) = k_b(p,p) - \rkl{\mathbb K \nabla p^\Omega, \nabla p^\Omega}_{0,\Omega} -
		\rkl{bK^\tau \nabla_\tau p^\Sigma, \nabla_\tau p^\Sigma}_{0,\Sigma}. \]
	From the Cauchy--Schwarz inequality and the $\varepsilon$-weighted Young inequality
	\[  \abs{ab} = \big \vert \varepsilon^{\frac 12} a \big \vert
		\big \vert \varepsilon^{-\frac 12}b \big \vert
		\le\frac 12\rkl{\varepsilon a^2 + \varepsilon^{-1}b^2} \qquad \forall a,b \in \R,
		\quad \varepsilon > 0, \]
	we deduce that for any $\varepsilon_1 > 0$
	\begin{align*}
	\rkl{ b^{-1} K^\nu p^\Sigma, \skl{p^\Omega} }_{0,\Sigma}
	& \le K^\nu \norm{  p^\Sigma}_{0,b^{-1},\Sigma} K^\nu \norm{  \skl{p^\Omega}}_{0,b^{-1},\Sigma} \\
	&\le \frac 12 \rkl{  \varepsilon_1 K^\nu \norm{  p^\Sigma}_{0,b^{-1},\Sigma}^2
		+ \frac {K^\nu}{\varepsilon_1} \norm{  \skl{p^\Omega}}_{0,b^{-1},\Sigma}^2},
	\end{align*}
	and that for any $\varepsilon_2 > 0$
	\begin{align*}
		 \sprod{\kappa K^\nu \jkl{p^\Omega}}{p^\Sigma}_{0,\Sigma}
			 & \le K^\nu b_{\text{max}} \kappa_{\text{max}} 
				 \sprod{\frac 1b \jkl{p^\Omega}}{p^\Sigma}_{0,\Sigma} \\
			 & \le \frac{1}{2}
				 \rkl{\varepsilon_2  K^\nu b_{\text{max}} \kappa_{\text{max}} 
				 	\norm{\jkl{p^\Omega}}_{0,b^{-1},\Sigma}^2 
				 	+ \frac{K^\nu b_{\text{max}} 
				 		\kappa_{\text{max}} }{\varepsilon_2} 
						 	\norm{p^\Sigma}_{0,b^{-1},\Sigma}^2  }. 
	\end{align*}
	Using these inequalities we can find a lower bound for $g_b(p)$
	\begin{align*}
	g_b(p) &= \frac{4K^\nu}{2\xi - 1} \Bigg[ \norm{ \skl{p^\Omega}}_{0,b^{-1},\Sigma}^2
	+ \norm{  p^\Sigma}_{0,b^{-1},\Sigma}^2 	
	- 2 \rkl{ b^{-1} p^\Sigma , \skl{p^\Omega}  }_{0,\Sigma} \Bigg]\\
	& \quad
	+  K^\nu\norm{ \jkl{ p^\Omega}}_{0,b^{-1},\Sigma}^2
	- \sprod{\kappa K^\nu \jkl{p^\Omega}}{p^\Sigma}_{0,\Sigma} \\
	& \ge  \frac{4 K^\nu}{2\xi - 1} \ekl{ \rkl{1 - \varepsilon_1}
		\norm{ \skl{p^\Omega}}_{0,b^{-1},\Sigma}^2
		+ \frac{\varepsilon_1 - 1}{\varepsilon_1} \norm{   p^\Sigma}_{0,b^{-1},\Sigma}^2 } \\
	& \qquad + \frac{2 -  b_{\text{max}} \kappa_{\text{max}}\varepsilon_2}{2} K^\nu
		 \norm{\jkl{ p^\Omega}}_{0,b^{-1},\Sigma}^2 
		- \frac{b_{\text{max}} \kappa_{\text{max}}}{2\varepsilon_2} K^\nu 
			\norm{p^\Sigma}_{0,b^{-1},\Sigma}^2\\
	& =\frac{4 K^\nu (1-\varepsilon_1)}{2\xi - 1}  
		\norm{  \skl{p^\Omega}}_{0,b^{-1},\Sigma}^2
		+  \frac{2 -  b_{\text{max}} \kappa_{\text{max}}\varepsilon_2}{2}
		 K^\nu \norm{\jkl{ p^\Omega}}_{0,b^{-1},\Sigma}^2  \\
	& \qquad + \frac{\rkl{8 \varepsilon_1 - 8}\varepsilon_2 - 
		b_{\text{max}} \kappa_{\text{max}}\rkl{2\xi - 1} \varepsilon_1}{2\varepsilon_1 \varepsilon_2 \rkl{2 \xi -1 }} K^\nu 
		\norm{p^\Sigma}_{0,b^{-1},\Sigma}^2.
	\end{align*}
	We use this result to find a lower bound for $k_b(p,p)$. 
	Introducing the constant $\eta \in (0,1)$, we have
	\begin{align*}
	k_b(p,p) & = \big \| \mathbb K^{\frac 12} \nabla p^\Omega \big \|_{0,\Omega}^2 +  \big \|\rkl{K^\tau}^{\frac 12} \nabla_\tau p^\Sigma \big \|_{0,b,\Sigma}^2 + g_b(p) \\
		& \ge  K \norm{ \nabla p^\Omega }_{0,\Omega}^2
			+ K^\tau \norm{\nabla_\tau p^\Sigma}_{0,b,\Sigma}^2
			+ g_b(p) \\
		& \ge \frac{K}{2}\min\{1,C_\Omega\} \norm{p^\Omega}_{1,\Omega}^2 
			+ \eta K^\tau C_\Sigma \norm{p^\Sigma}_{0,b^{-1},\Sigma}^2
			+ (1-\eta)K^\tau \norm{\nabla_\tau p^\Sigma }_{0,b,\Sigma}^2
			+ g_b (p) \\
		& \ge \frac{K}{2}\min\{1,C_\Omega\} \norm{p^\Omega}_{1,\Omega}^2 
		+ \eta K^\tau C_\Sigma \norm{p^\Sigma}_{0,b^{-1},\Sigma}^2
		+ (1-\eta)K^\tau \norm{\nabla_\tau p^\Sigma }_{0,b,\Sigma}^2 \\
		& \qquad + \frac{4 K^\nu (1-\varepsilon_1)}{2\xi - 1}  
		\norm{  \skl{p^\Omega}}_{0,b^{-1},\Sigma}^2
		+  \frac{2 -  b_{\text{max}} \kappa_{\text{max}}\varepsilon_2}{2}
		K^\nu \norm{\jkl{ p^\Omega}}_{0,b^{-1},\Sigma}^2  \\
		& \qquad + \frac{\rkl{8 \varepsilon_1 - 8}\varepsilon_2 - 
			b_{\text{max}} \kappa_{\text{max}}\rkl{2\xi - 1} \varepsilon_1}{2\varepsilon_1 \varepsilon_2} K^\nu 
		\norm{p^\Sigma}_{0,b^{-1},\Sigma}^2 \\
		& = \frac{K}{2}\min\{1,C_\Omega\} \norm{p^\Omega}_{1,\Omega}^2 
			+  (1-\eta)K^\tau \norm{\nabla_\tau p^\Sigma }_{0,b,\Sigma}^2 \\
		& \qquad + \frac{4 K^\nu (1-\varepsilon_1)}{2\xi - 1}  
		\norm{  \skl{p^\Omega}}_{0,b^{-1},\Sigma}^2
		+  \frac{2 -  b_{\text{max}} \kappa_{\text{max}}\varepsilon_2}{2}
		K^\nu \norm{\jkl{ p^\Omega}}_{0,b^{-1},\Sigma}^2  \\
		& \qquad + \frac{\ekl{\rkl{ 2\eta K^\tau C_\Sigma (2 \xi - 1) +  8
					K^\nu }\varepsilon_1 - 8K^\nu}\varepsilon_2 - 
			b_{\text{max}} \kappa_{\text{max}}\rkl{2\xi - 1} K^\nu \varepsilon_1}{2\varepsilon_1 \varepsilon_2 \rkl{2 \xi -1 }} 
		\norm{p^\Sigma}_{0,b^{-1},\Sigma}^2, 
	\end{align*}
	where we have used the Poincar\'e inequalities on $\Sigma$ and $\Omega$ with positive constants $C_\Sigma$ and $C_\Omega$.

	To ensure coercitivity all coefficients in the last estimate need to be positive.  We therefore have to find constants $\varepsilon_1 \in (0,1)$, 
	$\varepsilon_2 \in \rkl{0, \frac{2}{b_{\text{max}} \kappa_{\text{max}}}}$,
	and $\eta \in (0,1)$, such that 
	\[ \frac{\ekl{\rkl{ 2\eta K^\tau C_\Sigma (2 \xi - 1) +  8
				K^\nu }\varepsilon_1 - 8K^\nu}\varepsilon_2 - 
		b_{\text{max}} \kappa_{\text{max}}\rkl{2\xi - 1} K^\nu \varepsilon_1}{2\varepsilon_1 \varepsilon_2 \rkl{2 \xi -1 }} \ge 0. \]
	This inequality holds if and only if
	\begin{equation} \label{in:coercivity}
	 \ekl{\rkl{ 2\eta K^\tau C_\Sigma (2 \xi - 1) +  8
					K^\nu }\varepsilon_1 - 8K^\nu}\varepsilon_2 \ge
			b_{\text{max}} \kappa_{\text{max}}\rkl{2\xi - 1} K^\nu 
			\varepsilon_1. 
	\end{equation}
	
	Since $b_{\text{max}}^2 \kappa_{\text{max}}^2 < 
		\frac{4K^\tau C_\Sigma}{K^\nu}$ by assumption, we have 
	$0 \le \frac{b_{\text{max}}^2 \kappa_{\text{max}}^2 K^\nu}{4 K^\tau C_\Sigma}
	< 1$. Hence choose any $\eta > \frac{b_{\text{max}}^2 \kappa_{\text{max}}^2 K^\nu}{4 K^\tau C_\Sigma}$ less than one. This yields $4 K^\tau C_\Sigma \eta - 
	b_{\text{max}}^2 \kappa_{\text{max}}^2 K^\nu > 0$ and thus 
	\[ 0 < C_1 \colonequals \frac{16 K^\nu}{\ekl{4 K^\tau C_\Sigma \eta - 
		b_{\text{max}}^2 \kappa_{\text{max}}^2 K^ \nu} \rkl{2 \xi - 1} 
		+ 16 K^\nu} < 1.  \]
	Now choose $\varepsilon_1 > C_1$. Then we have 
	\begin{align*}
		2\ekl{\rkl{2 K^\tau C_\Sigma \eta \rkl{2 \xi - 1} + 8 K^\nu} 
		\varepsilon_1- 8 K^\nu} 
	    &=
		\rkl{4 K^\tau C_\Sigma \eta \rkl{2 \xi - 1} + 16 K^\nu} \varepsilon_1
		- 16 K^\nu \\
		&> b_{\text{max}}^2 \kappa_{\text{max}}^2 K^ \nu \rkl{2 \xi - 1} 
			\varepsilon_1,
	\end{align*}
	or equivalently
	\[	\frac{2}{b_{\text{max}} \kappa_{\text{max}} }  > C_2 \colonequals
		\frac{b_{\text{max}} \kappa_{\text{max}} K^\nu \rkl{2\xi - 1}\varepsilon_1}{\rkl{ 2\eta K^\tau C_\Sigma (2 \xi - 1) +  8
				K^\nu }\varepsilon_1 - 8K^\nu} \ge 0.
	\]
	Finally, choose $\varepsilon_2 > C_2$. Then estimate \eqref{in:coercivity} 
	is satisfied.
	
	This choice for $\varepsilon_1$, $\varepsilon_2$ and $eta$, implies that there exists a positive constant $C$, which only depends on $\Omega$, $\Sigma$, $K^\nu$, $K^\tau$, $K$, $\xi$, $\varepsilon_1$, $\varepsilon_2$ and $\eta$, such that
	\[ k_b(p,p) \ge C \norm{p}_{V_F}^2. \]
	Hence $k_b$ is coercive on $V_F$, and from the Lax--Milgram Lemma follows that \eqref{weak:summed_problem} has a
	unique solution in $V_F$, which in turn implies the original assertion.
\end{proof}
\begin{rem}\label{rem:asymptotics_fluid}
	Asymptotic analysis for the bulk fluid problem with fixed $p^\Sigma$, as it was done for the elasticity problem in Section \ref{subsubsec:weak-elasticity}, is difficult due to the
	$b$-dependent and non-symmetric coupling conditions on $\Sigma$. It is not clear to the authors whether these coupling conditions can be addressed by the standard spectral theory for elliptic problems with variable coefficients. Nonetheless, the primal form of the bulk fluid problem is a simple Laplace-type problem. The asymptotic expansion of this problem for various types of boundary conditions is well known \cite{mazia-rossmann-kozlov}. Eigenfunctions are of the form $c_1r^\alpha \sin(\alpha \Theta) + c_2 r^\alpha \cos(\alpha)$, and for plane problems with Neumann or Dirichlet boundary conditions on $\Sigma$ we have that $\alpha = \nicefrac 12$ is the lowest order term of the asymptotic expansion.
\end{rem}

\subsection{The Coupled Weak Problem}
Combining the three individual problems and coupling terms we obtain the coupled weak problem:
Find $(\mathbf u, p^\Omega, p^\Sigma) \in \mathbf V_E \times V_F^\Omega \times V_F^\Sigma$
  such that
\begin{subequations}
  \label{eq:coupledWeak}
  \begin{alignat}{2}
    a_F^\Omega (\mathbf u, p^\Omega, r^\Omega) - c_F(\mathbf u, p^\Sigma, r^\Omega) & = l_F^\Omega(r^\Omega) &\qquad& \forall \, r^\Omega \in V_{b,0}, \\
    a_F^\Sigma (\mathbf u, p^\Sigma, r^\Sigma) - c_F(\mathbf u, r^\Sigma, p^\Omega) - c_{F,\kappa}(\mathbf u, r^\Sigma, p^\Omega) & = l_F^\Sigma(\mathbf u,r^\Sigma) && \forall \, r^\Sigma \in H^1_{b,0}(\Sigma), \\
    a_E (\mathbf u, \mathbf v) - c_E^\Omega(\mathbf v, p^\Omega) + c_E^\Sigma(\mathbf v, p^\Sigma)  & = l_E^\Omega(\mathbf v) && \forall \, \mathbf v \in \mathbf V_0.
  \end{alignat}
\end{subequations}
  All terms in these three equations are as defined in the previous sections,  
  with the only difference that the dependencies on $b = \jkl{\mathbf u} \cdot \nu$ are replaced by dependencies on $\mathbf u$.

  The complete problem~\eqref{eq:coupledWeak} is nonlinear through the intricate influence of the changing crack width $b$
  on the fluid flow in the fracture.  Attempting to solve the entire system monolithically using a Newton-type method
  is certainly an option.  However, as both the fluid--fluid problem and the elasticity problem are linear when
  regarded separately, it is much more convenient to use a substructuring solver that iterates between the two.
  As the numerical tests in Section~\ref{sec:numerical_results} show, such a substructuring solver converges in
  very few iterations.

  We therefore proceed to write the coupled system as a fixed-point equation.  Since the fluid problems depend
  on the displacement only through the crack width $b = \jkl{\mathbf u}\cdot \nu$, we write the fixed-point equation
  in this variable.  By Theorem~\ref{thm:existence_fluid}, the coupled fluid--fluid problem~\eqref{weak:coupled-fluid}
  has a unique solution $(p^\Sigma, p^\Omega) \in V_F^\Sigma \times V_F^\Omega$ for any $b \in H^{1/2}_{00}(\Sigma)$ satisfying Assumptions \ref{hyp:weight}. We therefore obtain that the fluid solution operator
  \begin{equation*}
  	S_f : b \mapsto (p^\Omega, p^\Sigma)
  \end{equation*}
  is well-defined. Showing continuity of this operator is problematic, as the spaces $V_F^\Sigma$ and $V_F^\Omega$ depend on the argument $b$.

  Likewise, by Theorem~\ref{thm:existence-solid}, the elasticity problem~\eqref{eq:weak_elasticity_problem}
  has a unique solution for each $(p^\Omega, p^\Sigma) \in  V_F^\Sigma \times V_F^\Omega$.  Consequently,
  the elasticity solution operator
  \begin{equation*}
  	S_e : (p^\Omega, p^\Sigma) \mapsto \mathbf u
  \end{equation*}
  is well-defined, and Theorem~\ref{thm:existence-solid} additionally shows that tit is even continuous.

  Finally, the normal jump operator
  \begin{equation*}
  	j : \mathbf V_E \to H^{1/2}_{00}(\Sigma),
  	\qquad
  	j : \mathbf{u} \mapsto b \colonequals \jkl{ \mathbf u} \cdot \nu
  \end{equation*}
  is well-defined and continuous.
  Hence, we can write the weak coupled system as the fixed-point problem: Find $b \in H^{1/2}_{00}(\Sigma)$ satisfying Assumptions \ref{hyp:weight}
  such that
  \begin{equation}
  	\label{eq:fixed-point-equation}
  	b = (j \circ S_e \circ S_f)b.
  \end{equation}
  Showing existence of solutions to this equation is beyond the scope of this work.

\section{Discretization and Solver} \label{sec:discretization}

For the discretization of the coupled problem~\eqref{eq:coupledWeak} we use two unrelated grids: A $d$-dimensional one
for the bulk fluid and elasticity problems, and a $(d-1)$-dimensional one for the fracture fluid equation.
We use first-order Lagrange finite elements for the fluid equation on the fracture. Discretizing the bulk equations
is more challenging: Both solution fields are discontinuous at the fracture, and both fields
develop singularities at fracture tips.  These problems are overcome by an appropriate XFEM discretization.

\subsection{Finite Element Discretization of the Fracture Flow Equation}
 
	Denote by $\mathcal S_{h}$ a conforming and shape regular triangulation approximating $\Sigma$, with $h^\Sigma$ the maximum element diameter. Let $N^\Sigma$ denote the number of nodes in $\mathcal S_h$ and let $\mathcal N^\Sigma \colonequals \{ 1,\dots,N^\Sigma \}$ be the corresponding index set.

The averaged fracture pressure is discretized using first-order Lagrangian finite elements on $\Sigma$. We denote this finite element space by
\[
  V_{F,h}^\Sigma =  \spanx \{\varphi_i^\Sigma \}_{i \in \mathcal N^\Sigma},
\]
with nodal basis functions
   $\varphi_i^\Sigma : \Sigma \rightarrow \R$ associated with the fracture grid nodes $\mathbf s_i$, $i \in \mathcal N^\Sigma$.
  For a function $p_h^\Sigma \in V_{F,h}^\Sigma$ we denote by $p_i^\Sigma$ its coefficient with respect to $\varphi_i^\Sigma$ for all $i \in \mathcal N^\Sigma$.
\subsection{XFEM Spaces for the Bulk Problems}
  We assume that $\widetilde \Omega$ is a polygon and denote by $\mathcal T_{h}$ a conforming and shape regular
  triangulation of $\widetilde \Omega$. Let $h^\Omega$ denote the maximal diameter of an element of $\mathcal T_h$, and $N^\Omega$ the number of nodes in $\mathcal T_h$.
  The corresponding index set is $\mathcal N^\Omega \colonequals
  \{1, \dots, N^\Omega\}$.

  For a fixed $R \ge 0$, denote by $\mathcal J_R$ the set of indices of
  nodes that are within a cylindrical region around the crack front with radius $R$, or are contained in elements that
  intersect the crack front
  \[ \mathcal J_R \colonequals   \big\{ i \in\mathcal N^\Omega \mid \exists \, T \in \mathcal T_h \text{ with } \mathbf x_i \in T \text { and } \gamma \cap T \neq \emptyset \text { or }
    \dist(\mathbf x_i, \partial \Sigma) \le R \big\}. \]
  Further, introduce the set of indices of nodes contained in elements that are cut by the crack, but not already contained in $J_R$,
   \[\mathcal  K_R \colonequals \big\{ i \in \mathcal N^\Omega \setminus \mathcal  J_R \mid \exists \, T \in \mathcal T_h \text{ with } \Sigma \cap T \neq \emptyset  \big\}.\]
  Furthermore we introduce the standard finite element function  $\varphi_i^\Omega : \widetilde \Omega \rightarrow \R$ associated with the bulk grid node
  $i$, and the Heaviside function $H:\widetilde \Omega \rightarrow \R$,
  \begin{equation} \label{eq:def-heavyside}
    H(\mathbf x) =
    \begin{cases}
    -1 & \text{if $\mathbf x \in \Omega^-$}, \\
    1 & \text{else}.
    \end{cases}
  \end{equation}
	\begin{figure}[htb]
		\centering
		\begin{minipage}[htb]{0.49\linewidth}
			\centering
		     \begin{overpic}[width=0.6\linewidth,natwidth=110,natheight=87,clip, trim = 0 3 0 6]{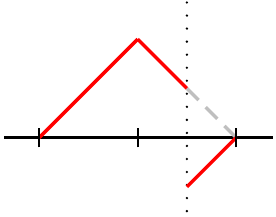}
				      	\put(8,15){$\mathbf x_{i-1}$}
				      	\put(47,15){$\mathbf x_i$}
				      	\put(80,15){$\mathbf x_{i+1}$}
		     \end{overpic}
		     \caption{1D Heaviside shape function $H\varphi_i^\Omega$}
		     \label{fig:bulkGrid3D}
		\end{minipage}
		\hfill
		\begin{minipage}[htb]{0.49\linewidth}
			\centering
			\begin{tikzpicture}[scale=0.7]
\pgfmathsetmacro{\gridx}{7}
\pgfmathsetmacro{\gridy}{5}

\draw[green!50!white,fill=green!10!white] (0.5*\gridx,0.5*\gridy) circle (0.25*\gridx);
\draw[step=1cm, very thin] (0,0) grid (\gridx,\gridy);
\draw[thick] (0,0.5*\gridy) --++ (0.5*\gridx,0);
\draw[thick, dashed] (0.5*\gridx,0.5*\gridy) --++ (0.5*\gridx,0);
\node[circle, fill=black,draw=black,inner sep=0pt,minimum size=4pt] at (0.5*\gridx,0.5*\gridy) {};

\foreach \x in {0,1} 
	\foreach \y in {2,3}
		\node[circle,fill=red!60!black,draw=black, inner sep=0pt,minimum size=5pt] at (\x,\y) {};
\foreach \x in {2,3,4,5} 
	\foreach \y in {2,3}
		\node[rectangle, fill=green!70!black,draw=black,inner sep=0pt,minimum size=5pt] at (\x,\y) {};
\foreach \x in {3,4} 
	\foreach \y in {1,4}
		\node[rectangle, fill=green!70!black,draw=black,inner sep=0pt,minimum size=5pt] at (\x,\y) {};
\end{tikzpicture}
			\caption{Node types for 2D enrichment -- red circles: Heaviside enriched nodes; green squares: crack tip function enriched nodes}
		\end{minipage}
	\end{figure}
\subsubsection{Discrete Displacement Space}

   We define $\mathbf e_\alpha$ as the $\alpha$-th canonical basis vector of $\R^d$, and
  \[\mathbf V_{E,h} \colonequals \spanx \left[ \bigcup_{\alpha=1}^d \Big(  \{\varphi_i^\Omega \mathbf e_\alpha \}_{i \in \mathcal N^\Omega} \cup  \{H \varphi_i^\Omega \mathbf e_\alpha \}_{i \in \mathcal K_R}
   \cup  \skl{F_j \varphi_i^\Omega \mathbf e_\alpha \mid j \in \{1,2,3,4\}, i \in \mathcal J_R} \Big) \right],
  \]
  and
  \begin{equation} \label{enr:solid}
  	\skl{F_j(r,\Theta)}_{j=1}^4 = \skl{\sqrt r \sin \rkl{\frac{\Theta}{2}}, \sqrt r \cos \rkl{\frac{\Theta}{2}}, \sqrt r \sin \rkl{\frac{\Theta}{2}} \sin\rkl{\Theta},
  		\sqrt r \cos \rkl{\frac{\Theta}{2}} \sin\rkl{\Theta} }.
  \end{equation}  
   Here $(r,\Theta)$ are the polar coordinates around the crack front introduced in Section \ref{subsubsec:weak-elasticity}. 

  The tip enrichment functions $F_j$, $j=1,\dots,4$, are the standard enrichment functions in the XFEM theory for linear elasticity problems. As shown in Theorem \ref{cor:asymptotics} these functions span the first order asymptotic expansion in a neighborhood of the crack tip.

  We can represent each function $\mathbf u _h = \sum_{\alpha=1}^{d} u_{h,\alpha}^\Omega \mathbf e_\alpha \in \mathbf V_{E,h}$ by
  \[ u_{h,\alpha}^\Omega(\mathbf x) =  \sum_{i \in \mathcal N^\Omega} u_{i,\alpha} \varphi_i^\Omega (\mathbf x)
	      + \sum_{i \in \mathcal J_R} v_{i,\alpha} H(\mathbf x)  \varphi_i^\Omega (\mathbf x)
	      + \sum_{i \in \mathcal K_R} \sum_{j=1}^4 c_{i,\alpha}^{(j)} F_j(\mathbf x)  \varphi_i^\Omega (\mathbf x), \qquad \alpha=1,\dots,d.  \]
  The average and jump of $\mathbf u_h$ at the fracture can be determined easily using only coefficient values and the standard Lagrangian hat function
  \begin{align*}
      \jkl{u_{h,\alpha}^\Omega(\mathbf s)} & =  \sum_{i \in \mathcal J_R} 2v_{i,\alpha}^\Omega \varphi_i^\Omega (\mathbf s)
	      + \sum_{i \in \mathcal  K_R} 2 c_{i,\alpha}^{(1)} \sqrt{r(\mathbf s)}  \varphi_i^\Omega (\mathbf s), \\
      \skl{u_{h,\alpha}^\Omega(\mathbf s)} & =   \sum_{i \in \mathcal N^\Omega} u_{i,\alpha}^\Omega \varphi_i^\Omega (\mathbf s),
  \end{align*}
  for all $\mathbf s \in \Sigma$. These relations make this set of enrichment functions particularly easy to work with. 
  \subsubsection{Discrete Pressure Space}
  The discretization space of the bulk pore pressure is defined by
  \[ V_{F,h}^\Omega \colonequals \spanx \{\varphi_i^\Omega \}_{i \in \mathcal N^\Omega} \cup \spanx \{H \varphi_i^\Omega \}_{i \in \mathcal K_R}  \cup \spanx \{G_j \varphi_i^\Omega \mid j = 1,2 \}_{i \in \mathcal J_R}, \]
  where $H$ is Heaviside function \eqref{eq:def-heavyside}
  and
  \begin{equation} \label{enr:fluid}
	G_1(r,\Theta) = F_1(r,\Theta) = \sqrt{r} \sin \rkl{\frac \Theta 2}, \quad
  	G_2(r,\Theta) = F_2(r,\Theta) =\sqrt{r} \cos \rkl{\frac \Theta 2}. 
  \end{equation}
  As mentioned in Remark \ref{rem:asymptotics_fluid}, the functions $G_1$ and $G_2$ span
  the first order asymptotic expansion in the near-tip approximation of the solution of the Laplace equation.
  Each $p_h^\Omega \in V_{F,h}^\Omega$ can be represented by
  \[       p_h^\Omega(\mathbf x) = \sum_{i \in \mathcal N^\Omega} p_i^\Omega \varphi_i^\Omega (\mathbf x) + \sum_{i \in K_r} q_i^\Omega H(\mathbf x)  \varphi_i^\Omega (\mathbf x)
  + \sum_{i \in J_r} \sum_{j=1}^2 r_{ij}^\Omega G_j(r,\Theta)  \varphi_i^\Omega (\mathbf x) . \]
  As for the displacement enrichment functions, we have for all $\mathbf s \in \Sigma$ that
  \begin{align*}
  \jkl{p_h^\Omega(\mathbf s)} & =  \sum_{i \in J_r} 2q_i^\Omega \varphi_i^\Omega (\mathbf s) +\sum_{i \in K_r}
  2r_{i1}^\Omega \sqrt{r} \varphi_i^\Omega (\mathbf s), \\
  \skl{p_h^\Omega(\mathbf s)} & =  \sum_{i \in \mathcal{N}^\Omega} p_i^\Omega \varphi_i^\Omega (\mathbf s).
  \end{align*}
  Hence the average and the jump of a pressure finite element function can be easily evaluated.
\subsection{Discrete Coupled Problems}
  The finite element spaces defined in the previous sections are conforming, and we can therefore obtain the discrete problem formulation by restricting the weak problem 
  \eqref{eq:coupledWeak} to the finite element spaces. The result reads:
  Find $(\mathbf u_h, p^\Omega_h, p^\Sigma_h) \in \mathbf V_{E,h} \times V_{F,h}^\Omega \times V_{F,h}^\Sigma$
  such that
  \begin{subequations}
  	\label{eq:coupledWeakDiscrete}
	  \begin{alignat}{2}
	    a_F^\Omega (\mathbf u_h, p^\Omega_h, r^\Omega_h) - c_F(\mathbf u_h, p^\Sigma_h, r^\Omega_h) & = l_F^\Omega(r^\Omega_h) &\qquad& \forall \, r^\Omega_h \in V_{F,h}^\Omega, \\
	    a_F^\Sigma (\mathbf u_h, p^\Sigma_h, r^\Sigma_h) - c_F(\mathbf u_h, r^\Sigma_h, p^\Omega_h) & = l_F^\Sigma(\mathbf u_h, r^\Sigma_h) && \forall \, r^\Sigma_h \in V_{F,h}^\Sigma, \\
	    a_E (\mathbf u_h, \mathbf v_h) - c_E^\Omega(\mathbf v_h, p^\Omega_h) + c_E^\Sigma(\mathbf v_h, p^\Sigma_h)  & = l_E^\Omega(\mathbf v_h) && \forall \, \mathbf v_h \in \mathbf V_{E,h}.
	  \end{alignat}
  \end{subequations}
  In view of numerically solving this system with a substructuring method, we write it as two separate linear problems,
  connected by a nonlinear coupling condition. The first subproblem is the
  discrete coupled fluid--fluid problem: For given $b_h \in H^{\frac 12}_{00}(\Sigma)$,
  find $(p_h^\Omega, p_h^\Sigma)  \in V_{F,h}^\Omega \times V_{F,h}^\Sigma$ such that
  \begin{subequations}
   \label{weak:discrete-fluid}
  \begin{alignat}{2}
      a_{F,b_h}^\Omega(p_h^\Omega, r_h^\Omega) - c_{F,b_h}(p_h^\Sigma, r_h^\Omega) &= l_F^\Omega(r_h^\Omega) &\qquad& \forall \, r_h^\Omega \in V_{F,h}^\Omega, \\
      a_{F,b_h}^\Sigma(p_h^\Sigma, r_h^\Sigma) - c_{F,b_h}(r_h^\Sigma, p_h^\Omega)
      - c_{F,\kappa, b_h}(r_h^\Sigma, p_h^\Omega) &= l_{F,b_h}^\Sigma(r_h^\Sigma) && \forall \, r_h^\Sigma \in V_{F,h}^\Sigma. 
  \end{alignat}
  \end{subequations}
 The second subproblem is the discrete weak elasticity problem: For given
  $(p_h^\Omega, p_h^\Sigma)  \in V_{F,h}^\Omega \times V_{F,h}^\Sigma$, find $\mathbf u_h \in \mathbf V_{E,h}$ such that
 \begin{equation} \label{weak:discrete-solid}
    a_E(\mathbf u_h, \mathbf v_h) = l_{E,p_h^\Sigma, p_h^\Omega} (\mathbf v_h) \qquad \forall \, \mathbf v_h \in \mathbf V_{E,h}.
 \end{equation}
  The two subproblems are coupled nonlinearly via the discrete crack width function
  \[ b_h (\mathbf s) \colonequals \jkl{\mathbf u_h (\mathbf s)}\cdot \nu =
	\sum_{\alpha=1}^d \left[ \sum_{i \in J_r} 2 v_{i,\alpha} \varphi_i (\mathbf s) \nu_\alpha
 	+ \sum_{i \in K_r} 2 c_{i,\alpha}^{(1)} \sqrt{r(\mathbf s)} \varphi_i (\mathbf s) \nu_\alpha \right], \qquad \mathbf s \in \Sigma.  \]
  The space of all discrete crack width functions is denoted by
  \[ H^{\frac 12}_h \colonequals \spanx \{  \varphi_i \nu_\alpha \mid i \in K_r, \alpha = 1,\dots, d \} \cup \spanx \{ \sqrt{r} \varphi_i \nu_\alpha \mid i \in J_r, \alpha = 1,\dots, d \}
   \]
   and we remark that this is a subspace of $H^{\frac 12}_{00}(\Sigma)$. 
   
  Existence and uniqueness of solutions to the subproblems are direct consequences of Theorems \ref{thm:existence-solid} and \ref{thm:existence_fluid}.

  In \cite{nicaise-renard-chahine} the following optimal convergence result for the elasticity subproblem \eqref{weak:discrete-solid} was proved for two-dimensional domains.
    \begin{thm}
     Assume that the displacement solution $\mathbf u$ of Problem \eqref{eq:weak_elasticity_problem} satisfies
      \[ \mathbf u - \mathbf u_s \in H^2(\Omega), \]
     where $\mathbf u_s$ denotes the singular part of $\mathbf u$ near the crack tip. Denote by $\mathbf u_h$ the
     solution of problem \eqref{weak:discrete-solid} and let $\chi$ be a smooth  cutoff function at the crack tip.
     Then
     \[ \norm{\mathbf u - \mathbf u_h}_{1,\Omega} \lesssim h \norm{\mathbf u - \chi \mathbf u_s}_{2,\Omega}. \]
    \end{thm}
    No corresponding result for the pressure subproblem is known.
\subsection{Substructuring Solver}
It is convenient to solve the nonlinear coupled problem \eqref{eq:coupledWeakDiscrete}
by iterating between the linear subproblems \eqref{weak:discrete-fluid} and 
\eqref{weak:discrete-solid}.

We suppose in the following that the discrete crack width function remains greater than  zero away from the crack tip.
Let $k=1,2,\dots$ be the iteration number, $b_{h,0} \in H^{\frac 12}_h$ be an initial discrete fracture width function, and $\beta \in (0,1]$ a damping parameter.
The substructuring solver proceeds in four steps:
\begin{enumerate}
	\item Solve the coupled fluid--fluid problem: Find $(p^\Sigma_{k}, p^\Omega_{k}) \in V_{F,h}^\Sigma \times V_{F,h}^\Omega$ such that
	\begin{alignat*}{2}
	a_{F,b_{k-1}}^\Omega \rkl{p^{\Omega}_k, r^\Omega} - c_{F,b_{k-1}} \rkl{p_k^{\Sigma}, r^\Omega}
	&= l_F^\Omega\rkl{r^\Omega}
	&\qquad& \forall \, r^\Omega \in V_{F,h}^\Omega, \\
	a_{F,b_{k-1}}^\Sigma \rkl{p_{k}^{\Sigma}, r^\Sigma} - c_{F,b_{k-1}}\rkl{r^\Sigma,  p_k^{\Omega}} &= l_{F,b_{k-1}}^\Sigma(r^\Sigma)
	&&  \forall r^\Sigma \in V_{F,h}^\Sigma,
	\end{alignat*}
	\item Solve the elasticity problem:  Find $ \widetilde{\mathbf u}_k \in \mathbf V_{E,h}$ such that
	\begin{equation*}
	a_{E,p^\Sigma_{k}, p^\Omega_{k}}^\Omega \rkl{  \widetilde{\mathbf{u}}_k, \mathbf v} = l_{E,p^\Sigma_{k}, p^\Omega_{k}} \rkl{\mathbf v}
	\qquad
	\forall \, \mathbf v \in \mathbf V_{E,h},
	\end{equation*}
	\item Damped update: \begin{equation*}
	\mathbf{u}_k \colonequals \rkl{1 - \beta} \mathbf{u}_{k-1} + \beta \widetilde{\mathbf{u}}_k.
	\end{equation*}
	\item Compute the normal jump at the fracture:
	\begin{equation*}
	b_k \colonequals \jkl{\mathbf u_k}\cdot \nu.
	\end{equation*}
	\end{enumerate}
	A rigorous proof of convergence of this iteration is left for future work. In numerical experiments we observe a very fast convergence (Section \ref{subsection:solverproperties}).

\section{Numerical Results}
\label{sec:numerical_results}
We close the article by giving a few numerical results. In particular, we show that the discretization error of the XFEM discretization proposed in the previous section behaves optimally when the mesh is refined. This justifies our choice of enrichment functions. Besides that, we show that the substructuring method converges very fast. Finally, we give a three-dimensional simulation. Our implementation is based on the \textsc{Dune} libraries\footnote{www.dune-project.org}, with the \texttt{dune-grid-glue}\footnote{www.dune-project.org/modules/dune-grid-glue} module to couple the bulk and fracture grids \cite{bastian-et-al}.
To solve systems of linear equations we use the UMFPACK direct solver \cite{suitesparse-project}.

\subsection{Substructuring Solver Convergence} \label{subsection:solverproperties}
To investigate the convergence rate of the iterative solver we consider a two-dimensional problem on the domain $\Omega = [0,1] \times \left[ -\frac 12, \frac 12 \right]$ (lengths are in kilometer).
The midsurface of the fracture is given by $\Sigma = \left[0, \frac 12 \right ] \times \{ 0 \}$.
The permeability tensors in the bulk and in the fracture
domains are homogeneous and isotropic, with $K = \SI{0.1}{\milli\text{D}}$ for the bulk and $K^\nu = K^\tau = \SI{100}{\text{D}}$ for the fracture. With these values the fluid can flow easily along and across the
fracture, whereas the rock matrix is much less permeable. Mechanically, the solid skeleton behaves according to the St.\,Venant--Kirchhoff
material law with Young's modulus $E=\SI{1}{\giga\pascal}$ and Poisson ratio $\nu = 0.3$.

\begin{figure}
\centering
\begin{subfigure}{.49\textwidth}
  \centering
  \begin{tikzpicture}
\pgfmathsetmacro{\squarex}{5}
\pgfmathsetmacro{\squarey}{5}

\draw[thick,red!70!black] (\squarex,0) --++ (0,\squarey) --++(-\squarex, 0) --++ (0,-\squarey);
\draw[thick,blue!70!black] (0,0) --++(\squarex,0);
\draw[thick,black] (0,0.5*\squarey) --++(0.5*\squarex,0);

\draw[black,fill=green!60!black,draw=none] (0,0.5*\squarey) circle (2pt);
\draw[black,fill=yellow!60!black,draw=none] (0.5*\squarex,0.5*\squarey) circle (2pt);

\node[font=\selectfont,red!70!black] at (0.5*\squarex,1.05*\squarey)
	{$\Gamma_N$};
\node[font=\selectfont,red!70!black] at (-0.1*\squarex,0.75*\squarey)
{$\Gamma_N$};
\node[font=\selectfont,red!70!black] at (-0.1*\squarex,0.25*\squarey)
{$\Gamma_N$};
\node[font=\selectfont,red!70!black] at (1.1*\squarex,0.5*\squarey)
{$\Gamma_N$};
\node[font=\selectfont,blue!70!black] at (0.5*\squarex,-0.075*\squarey)
	{$\Gamma_D$};

\node[font=\selectfont,green!60!black] at (-0.1*\squarex,0.5*\squarey)
	{$\gamma_D$};
\node[font=\selectfont,yellow!60!black] at (0.6*\squarex,0.5*\squarey)
{$\gamma_N$};
\end{tikzpicture}
  \caption{Boundary conditions}
  \label{fig:boundaryConditions}
\end{subfigure}
\begin{subfigure}{.46\textwidth}
  \centering
  \includegraphics[trim=0 -70 0 0,clip,width=0.76\textwidth,natwidth=918,natheight=986]{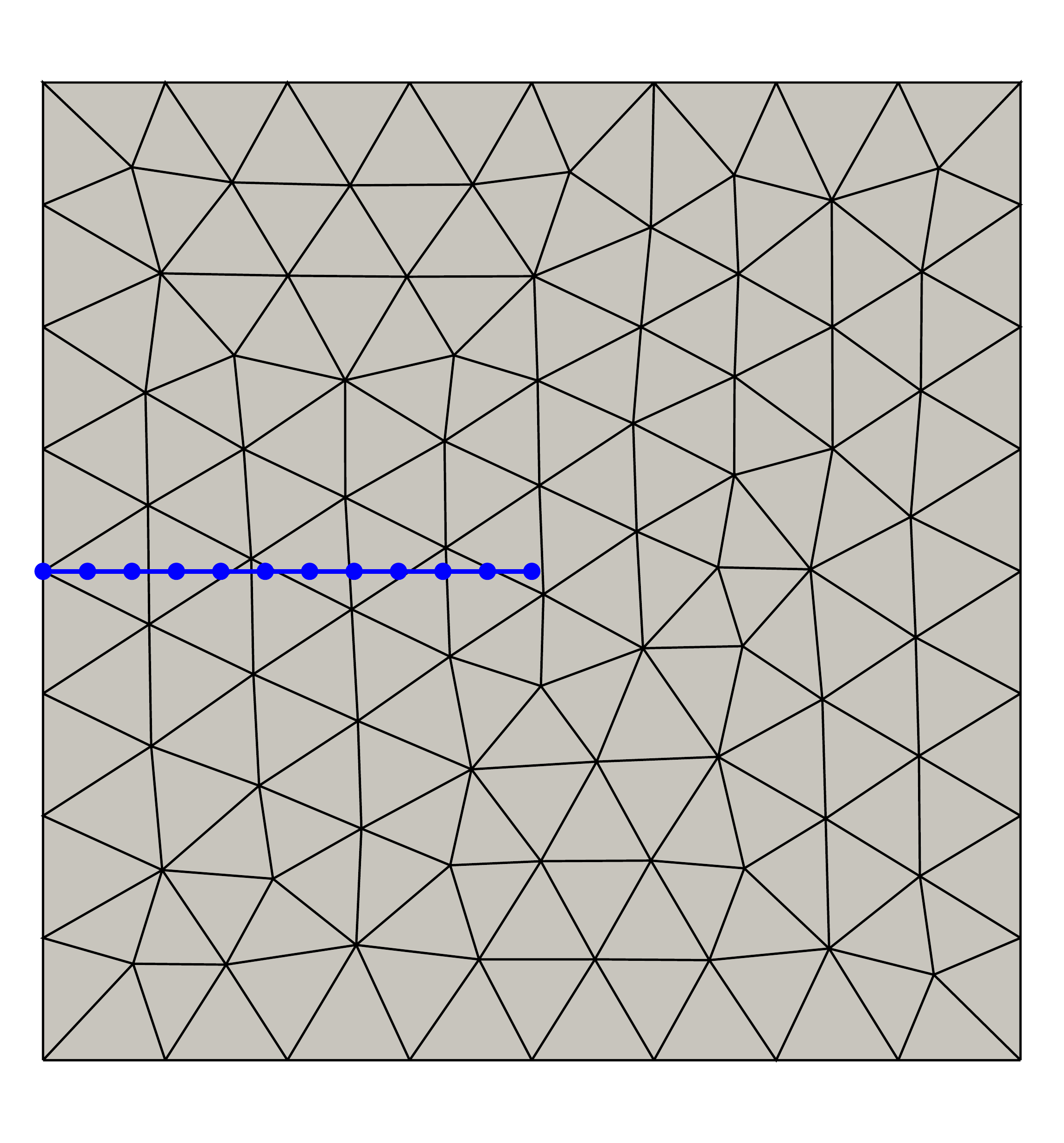}
  \caption{Grid}
  \label{fig:coarsestGrid}
\end{subfigure}
\caption{Problem setting}
\label{fig:test}
\end{figure}
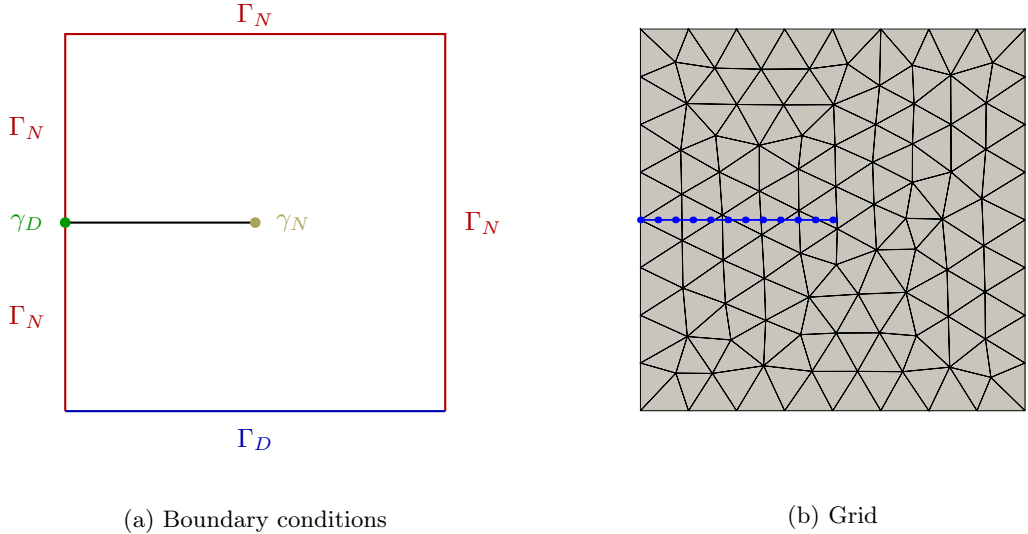
We prescribe zero Dirichlet boundary conditions for both the fluid problem and the
elasticity problem on the lower boundary of the bulk domain. Furthermore, a prescribed Dirichlet pressure of $p_0^\Sigma = \SI{0.5}{\mega\pascal}$ is applied
to the left boundary of the fracture $\gamma_D$. Zero Neumann boundary conditions are applied to the remaining parts of $\partial \Omega$, and to the crack tip $\gamma_N$
(see Figure \ref{fig:boundaryConditions}).
We set the XFEM enrichment radius to $R =0.125$, the solver damping parameter $\beta$ to one, and prescribe the initial crack width function $b_{h,0} = \sqrt{r} \cdot \SI{e-2}{\meter}$.

We discretize the bulk domain with an unstructured triangle grid, and the fracture domain with a uniform one-dimensional grid. Both grids are shown in Figure \ref{fig:coarsestGrid}. We create hierarchies of grids of different mesh size by refining both of them uniformly. For the test for the solver convergence speed we use up to $4$ steps of uniform refinement.

To measure the solver speed for a given pair of bulk and fracture grids, we compute a reference solution $\rkl{p_{h,\ast}^\Omega,p_{h,\ast}^\Sigma,\mathbf u_{h,\ast}}$ by performing $20$ substructuring iterations. After this many iterations, the solver is well beyond the limit of machine accuracy. We then compute the algebraic error of any iteration
$ (p_{h,k}^\Omega,p_{h,k}^\Sigma,\mathbf u_{h,k})$, $k=1,2,\dots$, as the relative $H^1$-error against this reference solution

\[ \text{err}_k = \frac{\big \| \big(p_{h,\ast}^\Omega,p_{h,\ast}^\Sigma, \textbf u_{h,\ast}\big) - \big(p_{h,k}^\Omega,p_{h,k}^\Sigma, \textbf u_{h,k}\big)
\big \|_1}{\big \| \big(p_{h,\ast}^\Omega,p_{h,\ast}^\Sigma, \textbf u_{h,\ast}\big) \big \|_1}.  \]

\begin{figure}[!ht]
	\begin{center}
\begin{tikzpicture}
\begin{semilogyaxis}[
 xlabel=Iteration Step,
 ylabel=Relative $H^1$-Error,
 legend pos=south west,
 width=0.55\linewidth
 ]
	\addplot[line width=1pt,red,mark=square*, mark size=2pt] table [col sep=comma,trim cells=true,x=it,y=r0] {errH};
	\addplot[line width=1pt,gray,mark=otimes*, mark size=2pt] table [col sep=comma,trim cells=true,x=it,y=r1] {errH};
	\addplot[line width=1pt,green,mark=triangle*, mark size=3pt] table [col sep=comma,trim cells=true,x=it,y=r2] {errH};
	\addplot[line width=1pt,cyan,mark=x, mark size=3pt] table [col sep=comma,trim cells=true,x=it,y=r3] {errH};
	\addplot[line width=1pt,blue,mark=+, mark size=3pt] table [col sep=comma,trim cells=true,x=it,y=r4] {errH};

	\legend{$0$ ref\\$1$ ref\\$2$ ref\\$3$ ref\\$4$ ref\\} 
\end{semilogyaxis}
\end{tikzpicture} 
	\end{center}
	\caption{Solver convergence}
	\label{fig:convergenceRate_Unstructured}
\end{figure}
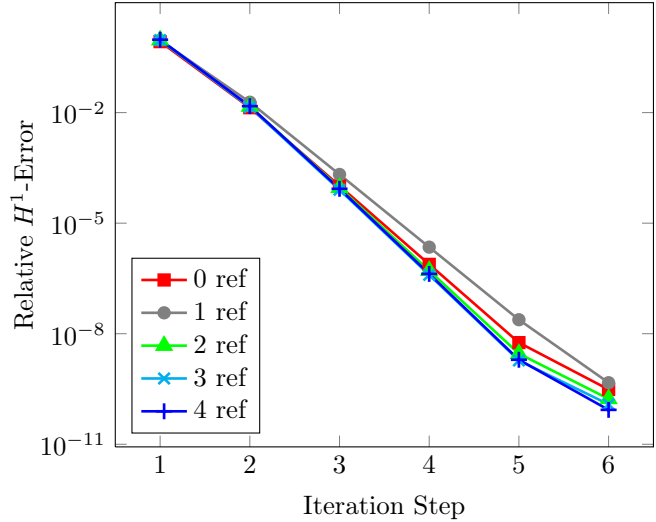

The result of the test can be seen in Figure \ref{fig:convergenceRate_Unstructured}, where we have plotted the algebraic error as a function of the iteration number $k$ for the different grid sizes. The error per iteration decays linearly on a logarithmic scale at a very fast rate. Indeed, a reduction of the relative error by a factor of $10^{-8}$ is achieved within only 5 iteration steps, and the convergence rate appears to be almost completely independent of the mesh size. We conclude that the substructuring method is a very competitive way to solve the coupled bulk--fracture system.

\subsection{Discretization Error Measurements}
We now measure the discretization error of the XFEM discretization proposed in Section \ref{sec:discretization}. This is a crucial experiment: The choice of the special tip enrichment functions \eqref{enr:solid} and \eqref{enr:fluid} can only be justified if the resulting discretization error behaves optimally as a function of the mesh size.

\begin{figure}[!ht]
  	\begin{center}
	\begin{tikzpicture}
	 \begin{loglogaxis}[
		 xlabel=Mesh Size $h^\Omega$,
		 ylabel=Relative Error $\mathbf u^\Omega$,
			legend pos=outer north east,
		 width=0.38\linewidth
		 ]
		\addplot[line width=1pt,red,mark=square*, mark size=2pt] 
			table [col sep=comma,trim cells=true,x=uoh,y=uoe1] {errUO};
		\addplot[line width=1pt,blue,mark=otimes*, mark size=2pt] 
			table [col sep=comma,trim cells=true,x=uoh,y=uoe2] {errUO};
		\addplot[line width=1pt,style=dotted,mark=+, mark size = 3pt,mark options={solid}]  
			table [col sep=comma,trim cells=true,x=uoh,y=uoh1] {errUO};
		\addplot[line width=1pt,style=dashed,mark=+, mark size = 3pt,mark options={solid}] 
			table [col sep=comma,trim cells=true,x=uoh,y=uoh2,mark=x] {errUO};
		
		\legend{$L^2$\\$H^1$\\$h$\\$h^2$\\} 
	\end{loglogaxis}
	\end{tikzpicture} \quad
	\begin{tikzpicture}
	\begin{loglogaxis}[
			xlabel=Mesh Size $p^\Omega$,
			ylabel=Relative Error $p^\Omega$,
			legend pos=outer north east,
			width=0.38\linewidth
		]
		\addplot[line width=1pt,red,mark=square*, mark size=2pt] 
			table [col sep=comma,trim cells=true,x=poh,y=poe1] {errPO};
		\addplot[line width=1pt,blue,mark=otimes*, mark size=2pt] 
			table [col sep=comma,trim cells=true,x=poh,y=poe2] {errPO};
		\addplot[line width=1pt,style=dotted,mark=+, mark size = 3pt,mark options={solid}]  
			table [col sep=comma,trim cells=true,x=poh,y=poh1] {errPO};
		\addplot[line width=1pt,style=dashed,mark=+, mark size = 3pt,mark options={solid}]  
			table [col sep=comma,trim cells=true,x=poh,y=poh2] {errPO};
		\addplot[line width=1pt,style=dotted,gray] 
			table [col sep=comma,trim cells=true,x=poh,y=por1] {errPO};
		\addplot[line width=1pt,style=dashed,gray] 
			table [col sep=comma,trim cells=true,x=poh,y=por2] {errPO};
		
		\legend{$L^2$\\$H^1$\\$h^{1.5}$\\$h^3$\\$h$\\$h^2$\\} 
	\end{loglogaxis}
	\end{tikzpicture}\\
	\begin{tikzpicture}
	\begin{loglogaxis}[
			xlabel=Mesh Size $h^\Sigma$,
			ylabel=Relative Error $p^\Sigma$,
			legend pos=outer north east,
			width=0.38\linewidth
		]
		\addplot[line width=1pt,red,mark=square*, mark size=2pt] 
			table [col sep=comma,trim cells=true,x=psh,y=pse1] {errPS};
		\addplot[line width=1pt,blue,mark=otimes*, mark size=2pt] 
			table [col sep=comma,trim cells=true,x=psh,y=pse2] {errPS};
		\addplot[line width=1pt,style=dotted,mark=+, mark size = 3pt,mark options={solid}]  
			table [col sep=comma,trim cells=true,x=psh,y=psh1] {errPS};
		\addplot[line width=1pt,style=dashed,mark=+, mark size = 3pt,mark options={solid}] 
			table [col sep=comma,trim cells=true,x=psh,y=psh2] {errPS};
		\addplot[line width=1pt,style=dotted,gray] 
			table [col sep=comma,trim cells=true,x=psh,y=psr1] {errPS};
		
		\legend{$L^2$\\$H^1$\\$h^{2}$\\$h^2$\\$h$\\} 
	\end{loglogaxis}
	\end{tikzpicture}
  	\end{center}
  	\caption{Discretization error}
  	\label{fig:discretizationError}
\end{figure}
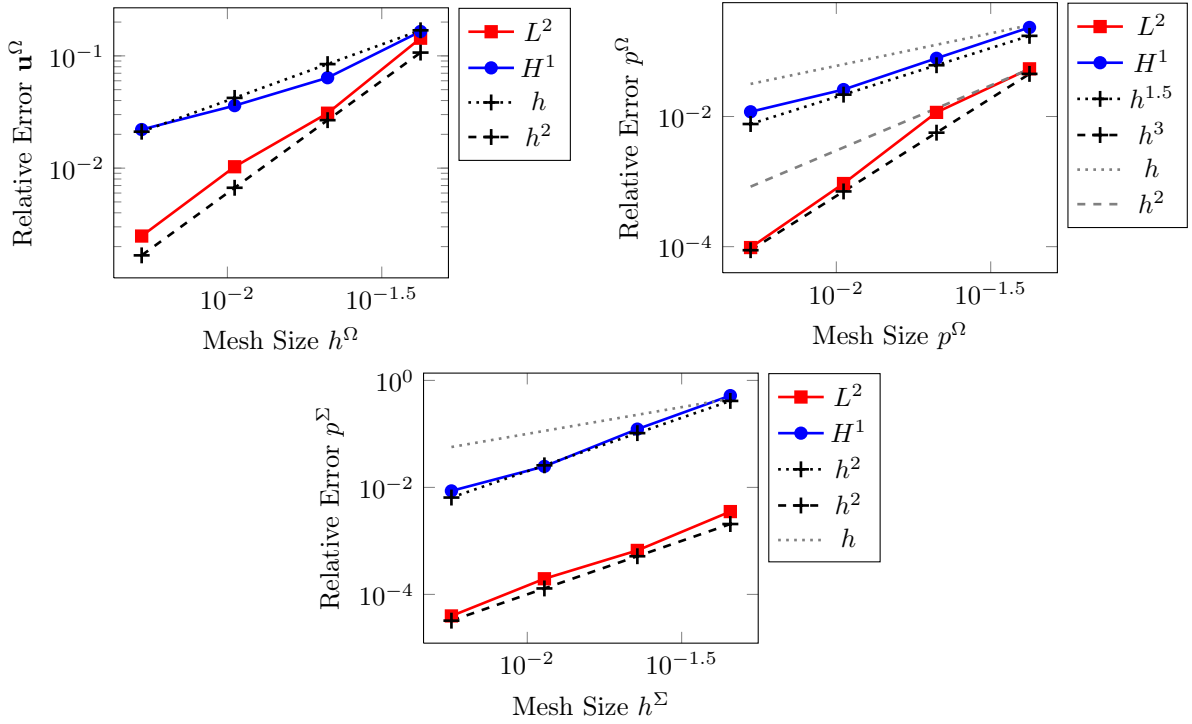

We use the same benchmark problem as in the previous section.
For each refinement level $j$ we have computed reference solutions $\big(p_{h,\ast}^{\Omega,j},p_{h,\ast}^{\Sigma,j},\mathbf u_{h,\ast}^j\big)$, for which the algebraic error is below $10^{-9}$. Of these five triples of functions, we pick
$\big(p_{h,\ast}^{\Omega,4},p_{h,\ast}^{\Sigma,4},\mathbf u_{h,\ast}^4\big)$, the one on the finest grid, to be the reference solution, and we compute $L^2$ and $H^1$ errors for the
coarser four solutions.
The results of this test can be seen in Figure \ref{fig:discretizationError}.
The errors are split up into the three components, that is: the matrix displacement $\mathbf u$, the bulk fluid pressure $p^\Omega$, and the fracture fluid pressure $p^\Sigma$. All three components show optimal error behavior, \ie, a decay of at least $\mathcal O(h^2)$ for the $L^2$ error and a decay of $\mathcal O(h)$ for the $H^1$ error. The discretization error of the bulk pressure $p^\Omega$ behaves even better. The $L^2$ error decays like
$\mathcal O (h^3)$ and the $H^1$ error as $\mathcal O (h^{\frac 32})$. This may be due to the fact that for this particular example the matrix pressure is continuous over the fracture and that the singularity of the fluid pressure gradient at the crack tip is reproduced by the second enrichment function $G_2$. The $L^2$ error of the fracture pressure behaves like $\mathcal O (h^2)$, which is better than expected as well. The reason for this is unclear. We conjecture that the higher regularity of the matrix pressure leads to a higher regularity in the coupling terms as well, and hence a stabilization effect for the fracture pressure gradient is induced. 

\subsection{A Three-Dimensional Example}
Finally, we demonstrate that our approach and implementation also work for three-dimensional problems.
We consider the domain $\Omega= \left[0,1\right]\times\left[ -\nicefrac 12, \nicefrac 12 \right]\times\left[ -\nicefrac 12, \nicefrac 12 \right]$, again in kilometers. The midsurface
of the fracture is defined as a half disc with radius $R= \frac 14 \si{\kilo \metre}$ and center $(0,0,0)^T$. The material properties are set as in the two-dimensional example,
\ie, the bulk permeability is $K = \SI{0.1}{\milli \text{D}}$, the fracture permeabilities are $K^\nu = K^\tau = \SI{100}{\text{D}}$ and for the matrix stiffness we set $E=\SI{1}{\text{G}\pascal}$ and $\nu = 0.3$.

\begin{figure}[!ht]
\begin{subfigure}{.32\textwidth}
  \centering
  \begin{tikzpicture}[scale=0.6]
\pgfmathsetmacro{\cubex}{5}
\pgfmathsetmacro{\cubey}{5}
\pgfmathsetmacro{\cubez}{5}

\pgfmathsetmacro{\circlex}{-\cubex}
\pgfmathsetmacro{\circley}{-0.5*\cubey}
\pgfmathsetmacro{\circlez}{-0.5*\cubez}

\draw[thick] (0,0,0) -- ++(-\cubex,0,0) -- ++(0,-\cubey,0) -- ++(\cubex,0,0) -- cycle;
\draw[thick] (0,0,0) -- ++(0,0,-\cubez) -- ++(0,-\cubey,0) -- ++(0,0,\cubez) -- cycle;
\draw[thick] (0,0,0) -- ++(-\cubex,0,0) -- ++(0,0,-\cubez) -- ++(\cubex,0,0) -- cycle;
\draw[thick, dashed] (-\cubex,-\cubey,0) -- ++(0,0,-\cubez) --++(\cubex,0,0);
\draw[thick, dashed] (-\cubex,-\cubey,-\cubez) -- ++(0,\cubey,0);


\draw [thick,domain=-90:90,variable=\x, fill=lightgray!30!white] plot ({\circlex + 1.25*cos(\x)}, {\circley},{\circlez+1.25*sin(\x)});
\draw[thick] (-\cubex,-0.5*\cubey,-0.25*\cubez) -- ++ (0,0,-0.5*\cubez);

\end{tikzpicture}
\end{subfigure} \hfill
\begin{subfigure}{.32\textwidth}
	\centering
  \begin{tikzpicture}[scale=0.6]
\pgfmathsetmacro{\cubex}{5}
\pgfmathsetmacro{\cubey}{5}
\pgfmathsetmacro{\cubez}{5}

\draw[fill=blue!20!white,draw=none] (-\cubex,-\cubey,0) -- ++(0,0,-\cubez) --++(\cubex,0,0) --++ (0,0,\cubez) -- cycle;

\node[font=\selectfont] at (-0.5*\cubex,-\cubey,-0.5*\cubez) {$\Gamma_D$};

\draw[thick] (0,0,0) -- ++(-\cubex,0,0) -- ++(0,-\cubey,0) -- ++(\cubex,0,0) -- cycle;
\draw[thick] (0,0,0) -- ++(0,0,-\cubez) -- ++(0,-\cubey,0) -- ++(0,0,\cubez) -- cycle;
\draw[thick] (0,0,0) -- ++(-\cubex,0,0) -- ++(0,0,-\cubez) -- ++(\cubex,0,0) -- cycle;
\draw[thick, dashed] (-\cubex,-\cubey,0) -- ++(0,0,-\cubez) --++(\cubex,0,0);
\draw[thick, dashed] (-\cubex,-\cubey,-\cubez) -- ++(0,\cubey,0);
\end{tikzpicture}
\end{subfigure} \hfill
\begin{subfigure}{.32\textwidth}
	\centering
  \begin{tikzpicture}[scale=0.6]
\pgfmathsetmacro{\cubex}{5}
\pgfmathsetmacro{\cubey}{5}
\pgfmathsetmacro{\cubez}{5}

\pgfmathsetmacro{\circlex}{-\cubex}
\pgfmathsetmacro{\circley}{-0.5*\cubey}
\pgfmathsetmacro{\circlez}{-0.5*\cubez}

\draw [thick,domain=-90:90,variable=\x, fill=lightgray!30!white,draw=yellow!60!black] plot ({\circlex + 2.5*cos(\x)}, {\circley},{\circlez+2.5*sin(\x)});
\draw[thick,green!60!black] (-\cubex,-0.5*\cubey,-\cubez) -- ++ (0,0,\cubez);

\node[font=\selectfont,yellow!60!black] 
 at (-0.3*\cubex,-0.5*\cubey,-0.5*\cubez) {$\gamma_N$};
\node[font=\selectfont,green!60!black] 
at (\circlex-0.2*\cubex,\circley,\circlez) {$\gamma_D$};
\end{tikzpicture}
\end{subfigure}
  \caption{Boundary Conditions for the three-dimensional example}
  \label{fig:setting_3d}
\end{figure}
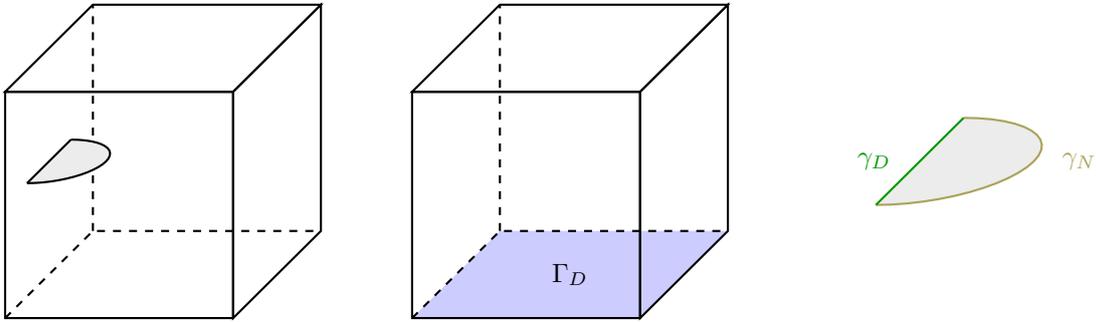

We apply zero Dirichlet boundary conditions
on the bottom and zero Neumann boundary conditions on the remaining part of the cube for the fluid bulk pressure and the displacement. We prescribe
a fracture pressure $p_0^\Sigma = \SI{0.5}{\mega\pascal}$ at the boundary part
$\gamma_D =\left\lbrace 0 \right\rbrace \times \left[-\nicefrac 14, \nicefrac 14 \right]\times\left\lbrace 0 \right\rbrace$ and zero Neumann conditions everywhere else on the boundary of the midsurface (see Figure \ref{fig:setting_3d}).

\begin{figure}[!ht]
	\centering
	\begin{subfigure}{.49\textwidth}
		\includegraphics[width=\linewidth]{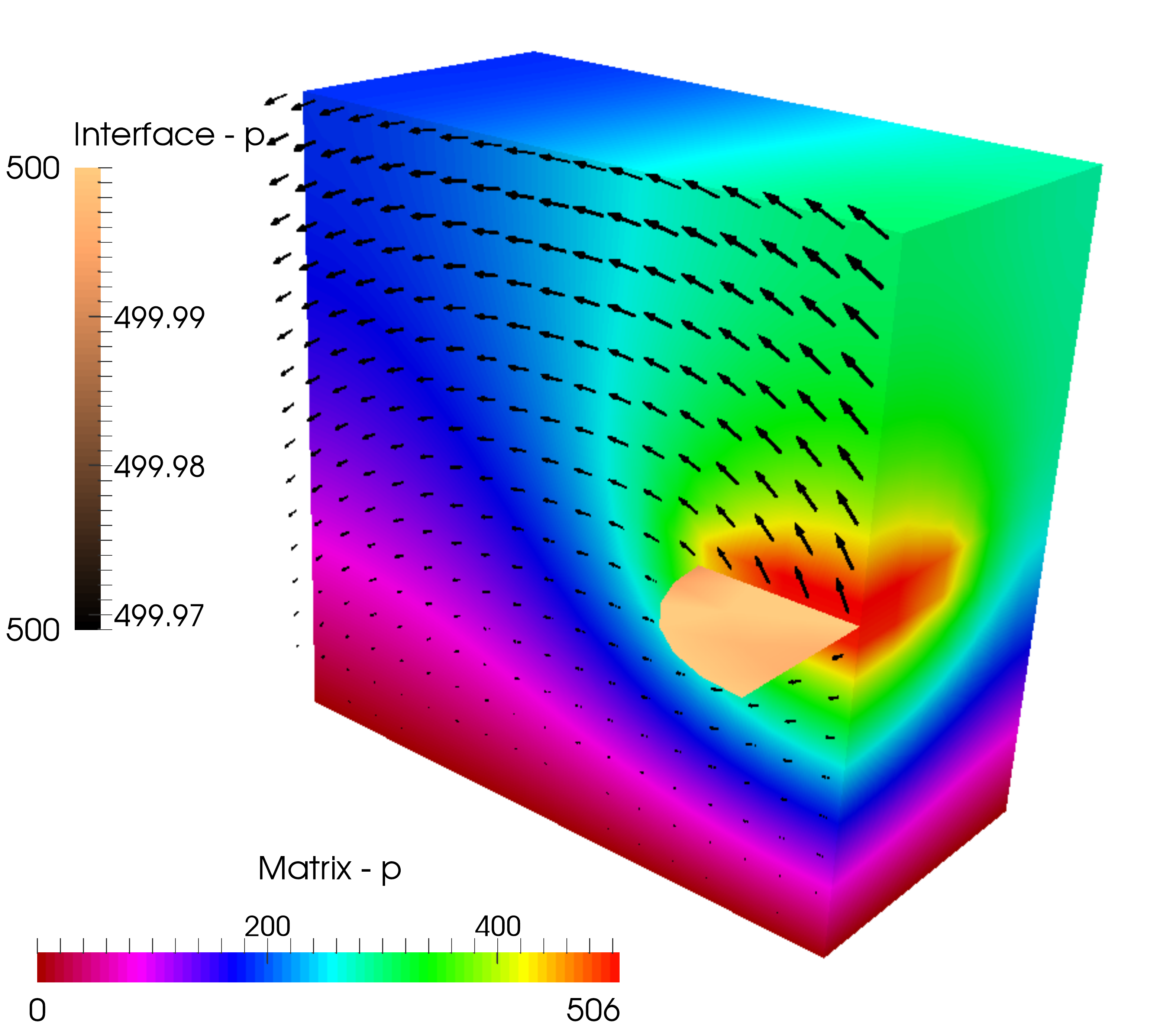}
	\end{subfigure}
	\begin{subfigure}{.49\textwidth}
		\includegraphics[width=\linewidth]{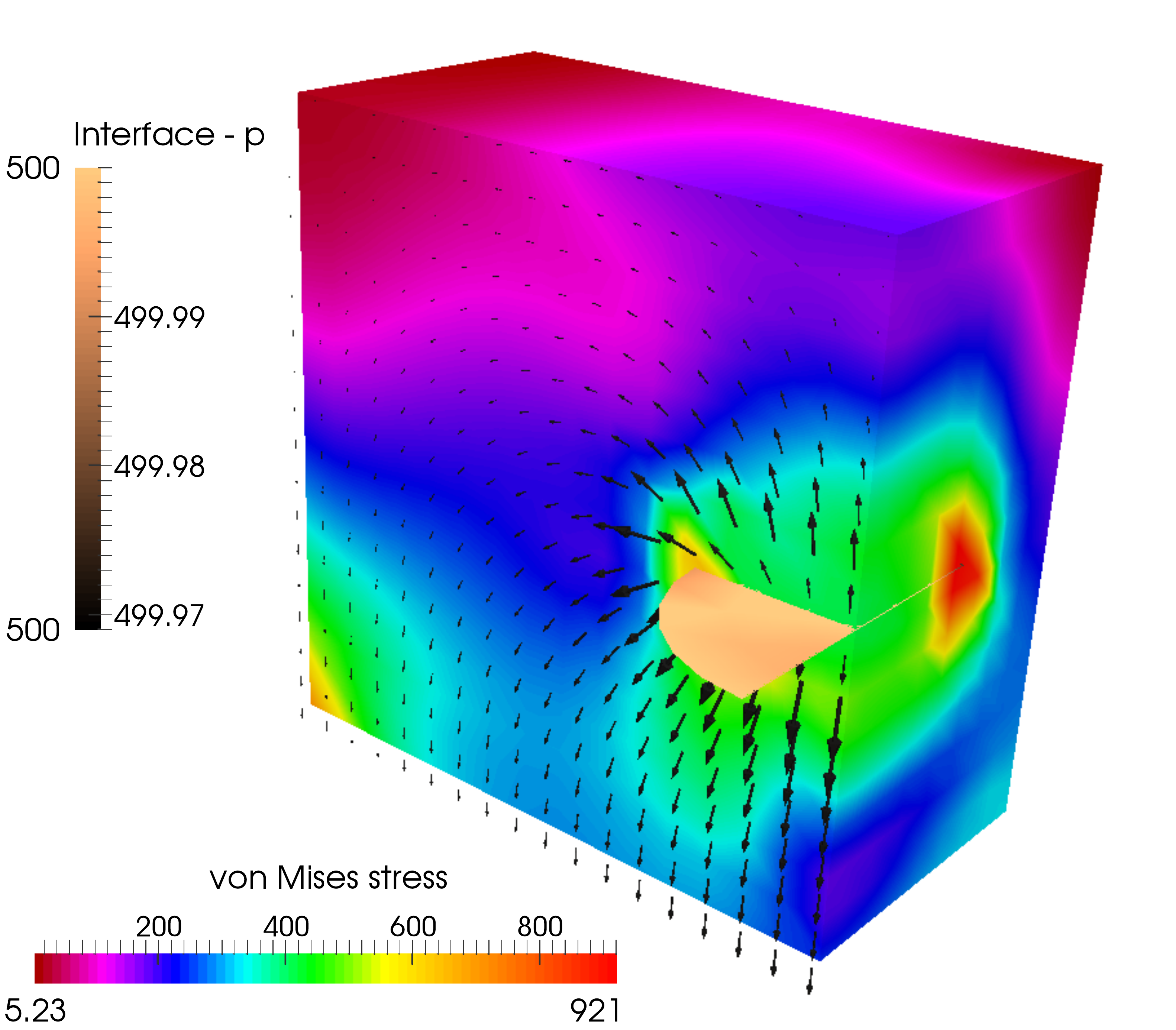}
	\end{subfigure}
	\caption{Results of the three-dimensional example; Left: Bulk pressure, displacement (scaled by $100$) and interface pressure; Right: Von Mises stress of the displacement field, bulk flow (scaled by $5 \cdot 10^{3}$) and interface pressure.}
		\label{fig:3dresults}
\end{figure}

The bulk domain is discretized using a structured tetrahedral grid with $20250$ elements in total.
The interface grid is an unstructured triangle grid with $90$ elements.
The XFEM enrichment radius is set to $R=0.125$, $\beta$ to one, and the  initial crack width
function is set to $b_{h,0} = \sqrt{r} \cdot \SI{e-2}{\meter}$.
We compute the von Mises stress of the displacement by evaluating the stress in the center of the elements. Nodal values are computed by averaging the values of the
elements the node is contained in. The result of this computation can be seen in Figure \ref{fig:3dresults}.

For this example, the substructuring solver again convergences in very few iterations. In the outcome, we observe that	
the fracture pressure is nearly constant. The fluid flows out of the fracture and thus a force is applied onto the fracture boundaries inducing a deformation, which opens the fracture.
The von Mises stress induced by the deformation has a maximum at the fracture front, which would lead to an enlargement of the fracture if fracture growth was included in the model.






\begin{thebibliography}{10}

\bibitem{aliabadi}
M.H. Aliabadi.
\newblock Boundary element formulations in fracture mechanics.
\newblock {\em Applied Mechanics Reviews}, 50:83--96, 1997.

\bibitem{angot-boyer-hubert}
P.~Angot, F.~Boyer, and F.~Hubert.
\newblock Asymptotic and numerical modelling of flows in fractured porous
  media.
\newblock {\em ESAIM: Mathematical Modelling and Numerical Analysis},
  43:239--275, 2009.

\bibitem{bastian-et-al}
P.~Bastian, M.~Blatt, A.~Dedner, C.~Engwer, R.~Klöfkorn, R.~Kornhuber,
  M.~Ohlberger, and O.~Sander.
\newblock A generic grid interface for parallel and adaptive scientific
  computing. part ii: implementation and tests in dune.
\newblock {\em Computing}, 82(2-3):121--138, 2008.


\bibitem{bourdin-francfort-marigo}
B.~Bourdin, G.A. Francfort, and J.-J. Marigo.
\newblock The variational approach to fracture.
\newblock {\em Journal of Elasticity}, 91(1-3):5--148, 2008.

\bibitem{brewster-mitrea}
Kevin Brewster and Marius Mitrea.
\newblock In weighted sobolev spaces on lipschitz manifolds.
\newblock {\em Memoirs on Differential Equations and Mathematical Physics},
  60:15--55, 2013.

\bibitem{costabel-dauge-edge-asymptotics}
M.~Costabel and M.~Dauge.
\newblock General edge asymptotics of solutions of second-order elliptic
  boundary value problems i.
\newblock {\em Proceedings of the Royal Society of Edinburgh: Section A
  Mathematics}, 123(01):109--155, 1993.
%

\bibitem{dangelo-scotti}
C.~D'Angelo and A.~Scotti.
\newblock A mixed finite element method for darcy flow in fractured porous
  media with non-matching grids.
\newblock {\em ESAIM: Mathematical Modelling and Numerical Analysis},
  46:465--489, 2012.

\bibitem{suitesparse-project}
T.~Davis.
\newblock Suitesparse project homepage. \texttt{URL:} \textsf{http://faculty.cse.tamu.edu/davis/welcome} (visited on 05/03/2015)
%

\bibitem{duran-lopez}
R.~G. Dur{\'a}n and F.~L{\'o}pez~Garc{\'i}a.
\newblock Solutions of the divergence and analysis of the stokes equations in
  planar h{\"o}lder-$\alpha$ domains.
\newblock {\em Mathematical Models and Methods in Applied Sciences},
  20(01):95--120, 2010.

\bibitem{edmunds-opic}
D.~E. Edmunds and B.~Opic.
\newblock Weighted poincaré and friedrichs inequalities.
\newblock {\em Journal of the London Mathematical Society}, s2-47(1):79--96,
  1993.
%

\bibitem{formaggia-et-al}
L.~Formaggia, A.~Fumagalli, A.~Scotti, and P.~Ruffo.
\newblock A reduced model for darcy’s problem in networks of fractures.
\newblock {\em ESAIM: Mathematical Modelling and Numerical Analysis},
  48:1089--1116, 2014.

\bibitem{fries-schaetzer-weber}
T.-P. Fries, M.~Sch\"atzer, and N.~Weber.
\newblock XFEM-Simulation of Hydraulic Fracturing in 3D with Emphasis on Stress Intensity Factors.
\newblock In: \textit{11th World Congress on Computational Mechanics (WCCM XI) and 5th European Congress on Computational Mechanics (ECCM V) and 6th European Congress on Computational Fluid Dynamics (ECFD VI)}. (Barcelona, Spain). Ed. by E.
Oñate, X. Oliver, and A. Huerta, CIMNE, pages 3282--3293, 2014.

\bibitem{fumagalli-scotti-13}
A.~Fumagalli and A.~Scotti.
\newblock A numerical method for two-phase flow in fractured porous media with
  non-matching grids.
\newblock {\em Advances in Water Resources}, 62, Part C:454 -- 464, 2013.
\newblock Computational Methods in Geologic CO2 Sequestration.

\bibitem{fumagalli-scotti-11}
A.~Fumagalli and A.~Scotti.
\newblock A reduced model for flow and transport in fractured porous media with
  non-matching grids.
\newblock In A.~Cangiani, R.L. Davidchack, E.~Georgoulis, A.N. Gorban,
  J.~Levesley, and M.V. Tretyakov, editors, {\em Numerical Mathematics and
  Advanced Applications 2011}, pages 499--507. Springer,
  2013.

\bibitem{garagash-detournay-adachi}
D.I. Garagash, E.~Detournay, and J.I. Adachi.
\newblock Multiscale tip asymptotics in hydraulic fracture with leak-off.
\newblock {\em Journal of Fluid Mechanics}, 669:260--297, 2011.

\bibitem{girault-kumar-wheeler}
V.~Girault, K.~Kumar, and M.~F. Wheeler.
\newblock Convergence of iterative coupling of geomechanics with flow in a fractured poroelastic medium.
\newblock {\em Computational Geosciences}, 2016.

\bibitem{girault-wheeler-ganis-mear}
V.~Girault, M.~F. Wheeler, B.~Ganis, and M.~E. Mear.
\newblock A lubrication fracture model in a poro-elastic medium.
\newblock {\em Mathematical Models and Methods in Applied Sciences},
  25(04):587--645, 2015.

\bibitem{inglis}
C.E. Inglis.
\newblock Stresses in plates due to the presence of cracks and sharp corners.
\newblock {\em Transactions of the Institute of Naval Architects}, 55:219--241,
  1913.

\bibitem{jaffre-mnejja-roberts}
J.~Jaffr\'e, M.~Mnejja, and J.E. Roberts.
\newblock A discrete fracture model for two-phase flow with matrix-fracture
  interaction.
\newblock {\em Procedia Computer Science}, 4(0):967--973, 2011.
\newblock Proceedings of the International Conference on Computational Science,
  \{ICCS\} 2011.

\bibitem{khludnev}
A.M. Khludnev and V.A. Kovtunenko.
\newblock {\em Analysis of Cracks in Solids}.
\newblock International series on advances in fracture. WIT Press, 2000.

\bibitem{mazia-rossmann-kozlov}
V.~Kozlov, V.G. Mazia, and J.~Rossmann.
\newblock {\em Spectral Problems Associated with Corner Singularities of
  Solutions to Elliptic Equations}.
\newblock Mathematical surveys and monographs. American Mathematical Society,
  2001.
%

\bibitem{kufner-opic-84}
A.~Kufner and B.~Opic.
\newblock How to define reasonably weighted sobolev spaces.
\newblock {\em Commentationes Mathematicae Universitatis Carolinae},
  025(3):537--554, 1984.

\bibitem{martin-jaffre-roberts}
Vincent Martin, J\'er\^ome Jaffré, and Jean~E. Roberts.
\newblock Modeling fractures and barriers as interfaces for flow in porous
  media.
\newblock {\em SIAM Journal on Scientific Computing}, 26(5):1667--1691, 2005.

\bibitem{mikelic-wheeler-wick}
A.~Mikeli\'c, M.~F. Wheeler, and T.~Wick.
\newblock A Phase-Field Method for Propagating Fluid-Filled Fractures Coupled to a Surrounding Porous Medium.
\newblock {\em Multiscale Modeling \& Simulation}, 13(1):367--398, 2015.

\bibitem{miller}
N.~Miller.
\newblock Weighted sobolev spaces and pseudodifferential operators with smooth
  symbols.
\newblock {\em Transactions of the American Mathematical Society},
  269(1):91--109, 1982.

\bibitem{moes-dolbow-belytschko}
N.~Mo\"es, J.~Dolbow, and T.~Belytschko.
\newblock A finite element method for crack growth without remeshing.
\newblock {\em International Journal for Numerical Methods in Engineering},
  46(1):131--150, 1999.

\bibitem{nicaise-renard-chahine}
Serge Nicaise, Yves Renard, and Elie Chahine.
\newblock Optimal convergence analysis for the extended finite element method.
\newblock {\em International Journal for Numerical Methods in Engineering},
  86(4-5):528--548, 2011.
%

\bibitem{kufner-opic-90}
B.~Opic and A.~Kufner.
\newblock {\em Hardy-type Inequalities}.
\newblock Pitman research notes in mathematics series. Longman Scientific \&
  Technical, 1990.
%

\bibitem{pommier-et-al}
S.~Pommier, A.~Gravouil, N.~Mo\"es, and A.~Combescure.
\newblock {\em Extended Finite Element Method for Crack Propagation}.
\newblock Iste Series. Wiley, 2011.

\bibitem{quarteroni_valli:1999}
Alfio Quarteroni and Alberto Valli.
\newblock {\em Domain Decomposition Methods for Partial Differential
  Equations}.
\newblock Oxford Science Publications, 1999.

\bibitem{rabczuk-et-al}
T.~Rabczuk, S.~Bordas, and G.~Zi.
\newblock Computational methods for fracture.
\newblock {\em Mathematical Problems in Engineering}, 2014, 2014.

\bibitem{rice}
J.~R. Rice.
\newblock Mathematical analysis in the mechanics of fracture.
\newblock {\em Fracture: an advanced treatise}, 2:191--311, 1968.

\bibitem{turesson}
B.O. Turesson.
\newblock {\em Nonlinear Potential Theory and Weighted Sobolev Spaces}.
\newblock Number Nr. 1736 in Lecture Notes in Mathematics. Springer, 2000.

\bibitem{watanabe-et-al}
N.~Watanabe, W.~Wang, J.~Taron, U.~J. Görke, and O.~Kolditz.
\newblock Lower-dimensional interface elements with local enrichment:
  application to coupled hydro-mechanical problems in discretely fractured
  porous media.
\newblock {\em International Journal for Numerical Methods in Engineering},
  90(8):1010--1034, 2012.

\bibitem{wheeler-wick-wollner}
M.~Wheeler, T.~Wick, and W.~Wollner.
\newblock An augmented-lagrangian method for the phase-field approach for
  pressurized fractures.
\newblock {\em Computer Methods in Applied Mechanics and Engineering},
  271(0):69--85, 2014.

\bibitem{wloka}
J.~Wloka.
\newblock {\em Partial Differential Equations}.
\newblock Cambridge University Press, 1987.

\bibitem{yazid-abdelkader-abdelmadjid}
A.~Yazid, N.~Abdelkader, and H.~Abdelmadjid.
\newblock A state-of-the-art review of the x-fem for computational fracture
  mechanics.
\newblock {\em Applied Mathematical Modelling}, 33(12):4269--4282, 2009.

\end{thebibliography}
\end{document}